\documentclass{itogi2016}

 \currentyear{2019}
 \currentvolume{???}

\usepackage[russian]{babel}
\usepackage[cp1251]{inputenc}

\usepackage{graphicx}

\usepackage{enumerate} 
\usepackage{verbatim}
\usepackage[mathscr]{eucal} 

\usepackage{xcolor} 

\usepackage[
unicode,colorlinks,linkcolor=blue,citecolor=red,bookmarksopen,pdfhighlight=/N]{hyperref}

\theoremstyle{plain} 
\newtheorem{theorem}{Теорема} 
 
\newtheorem*{theoremA}{Теорема A}

\newtheorem{problem}{Задача}

\newtheorem*{thHB}{Теорема Хана\,--\,Банаха}

\newtheorem*{thAL}{Теорема Энгера\,--\,Лембке}
\newtheorem*{thH}{Теорема Хёрмандера}
\newtheorem*{thHBK}{Теорема Хана\,--\,Банаха\,--\,Канторовича}
\newtheorem{lemma}{Лемма}[section]
\newtheorem{propos}{Предложение}[section] 
\newtheorem{corollary}{Следствие}[section]

\theoremstyle{definition}
\newtheorem{definition}{Определение}[section] 
\newtheorem{remark}{Замечание}[section] 

\newtheorem*{examples}{Примеры}

\newcommand{\bal}{\rm {bal}}
\newcommand{\RR}{\mathbb{R}} 
\newcommand{\CC}{\mathbb{C}} 
\newcommand{\NN}{\mathbb{N}}

\newcommand{\KK}{\mathbb{K}} 
\newcommand{\const}{{\rm const}} 
\newcommand{\tmin}{{\text{-}}} 
\newcommand{\card}{{\rm card\,}} 
\newcommand{\pr}{{\rm pr}} 
\newcommand{\dd}{\,{\rm d}}

\DeclareMathOperator{\dist}{dist} 
 
\DeclareMathOperator{\clos}{clos} 
\DeclareMathOperator{\Int}{int}
\DeclareMathOperator{\Meas}{Meas}

\DeclareMathOperator{\har}{har} 
\DeclareMathOperator{\Hol}{Hol}
\DeclareMathOperator{\Mer}{Mer}

\DeclareMathOperator{\comp}{cmp}

\DeclareMathOperator{\Zero}{Zero} 
\DeclareMathOperator{\sbh}{sbh}
\DeclareMathOperator{\psbh}{plsbh}
\DeclareMathOperator{\sbl}{sbl}
\DeclareMathOperator{\spl}{spl}
\DeclareMathOperator{\spf}{spf}
\DeclareMathOperator{\Spr}{spr}
\DeclareMathOperator{\ph}{ph}
\DeclareMathOperator{\phar}{plhar}
\DeclareMathOperator{\lin}{lin}
\DeclareMathOperator{\aff}{aff}
\DeclareMathOperator{\add}{add}
\DeclareMathOperator{\sba}{sba}
\DeclareMathOperator{\conv}{conv}   
\DeclareMathOperator{\conc}{conc}  
 
\DeclareMathOperator{\supp}{supp} 
\DeclareMathOperator{\loc}{loc}

\DeclareMathOperator{\incr}{incr}
\DeclareMathOperator{\reg}{reg}

\DeclareMathOperator{\lE}{lenv}
\DeclareMathOperator{\uE}{uenv} 

\renewcommand{\Re}{{\rm Re \,}}

\renewcommand{\theequation}{\thesection.\arabic{equation}}
\begin{document}

\title[Порядковые версии Теоремы Хана\,--\,Банаха и огибающие. II. Применения \dots]{Порядковые версии Теоремы Хана\,--\,Банаха и огибающие. II. Применения в теории функций}

\author[Хабибуллин Булат Нурмиевич]{Б. Н. Хабибуллин}
\address{Башкирский государственный университет
 (БашГУ)}
\email{khabib-bulat@mail.ru}

\author[Розит Алексей Петрович]{А. П. Розит}
\address{МБОУ <<Лицей №~60>> г.~Уфы}
\email{rozit@mail.ru}

\author[Хабибуллина Энже Булатовна]{Э. Б. Хабибуллина}
\address{Башкирский государственный университет
	(БашГУ)}
\email{khabibullinae@gmail.com}




\UDC{517.982, 517.57, 517.55, 517.53}

\begin{abstract}
В главе 1 рассматривается проблема существования верхней/нижней огибающей из выпуклого конуса или, более общ\'о, выпуклого множества для функций на проективном пределе векторных решёток  со значениями в пополнении пространства Канторовича или на расширенной вещественной прямой. Даются векторные, порядковые и топологические двойственные трактовки условий существования такой огибающей и метода её построения.   В главе 2 приводятся применения к существованию нетривиальной (плюри)субгармонической и/или (плюри)гармонической миноранты для функций в областях из конечномерного вещественного или комплексного пространства. Указываются общие подходы к задачам о нетривиальности весовых классов голоморфных функций, к описанию нулевых (под)множеств для таких классов голоморфных функций, к задаче представления мероморфной функции как частного голоморфных функций из заданного весового класса.
\end{abstract}

\thanks{
Исследование выполнено при поддержке РФФИ (проект № 16-01-00024а).}


\maketitle

\tableofcontents

\chapter{Огибающая}\label{ch_en}

\section{Введение. Истоки}
\setcounter{equation}{0}
\subsection{Алгебраические версии}\label{ist} Одна  из основ функционального анализа ---
\begin{thHB}[{\rm об огибающей \cite[II.1.3--4]{AK}, \cite{Buskes}, \cite{Narici}, \cite{FL}, \cite{NB}, \cite{S08}}] Пусть $X$ ---  ве\-к\-т\-о\-р\-ное пространство над полем вещественных чисел\/ $\RR$ 
и $f\colon X\to \RR$ --- функция, т.\,е. $f\in \RR^X$.
\begin{enumerate}[{\bf I.}]
	\item\label{I} Попарно эквивалентны следующие  три утверждения.     
\begin{enumerate}[{\bf 1.}]
\item\label{HBi} Функция $f$ сублинейная\/ {\rm (пишем $f\in \sbl \RR^X$)}, что здесь означает  одновременно
\begin{description}
\item[{\bf (ph)}] положительную  однородность\/ {\rm (пишем $f\in 
\ph \RR^X$):}
\begin{equation}\label{ph}
f(tx)=tf(x)\quad\text{для всех $x\in X$ и $t\in \RR^+:=\{t\in \RR\colon t\geq 0\}$},
\end{equation} 
где\/ $\RR^+$ можно заменить на\/ $\RR_*^+:=\RR\setminus \{0\}$,
\item[{\bf (sa)}] и субаддитивность\/  {\rm (пишем $f\in \sba \RR^X$):}  
\begin{equation}\label{sa}
f(x_1+x_2)\leq f(x_1)+f(x_2)\quad\text
 {для всех $x_1,x_2 \in X$}.
\end{equation}
\end{description}

\item\label{HBii} Для векторного пространства\/ $\lin \RR^X$ над\/ $\RR$  всех линейных функций $l \in \RR^X$:
\begin{equation}\label{lf}
l(t_1x_1+t_2x_2)=t_1l(x_1)+t_2l(x_2)\quad\text{для всех $x_1,x_2\in X$ и $t_1,t_2\in \RR$},
\end{equation} 
функция $f$ совпадает со своей  нижней огибающей  по $\lin \RR^X$: 
\begin{equation}\label{lenvf}
f(x)=\sup \bigl\{l(x)\colon \, l\in \lin \RR^X, \; 
l\bigm|_X \leq f\bigm|_X \bigr\}=:
{\lE}_{\lin \RR^X}^f (x)\quad\text{для всех $x\in X$},
\end{equation}
где $l\bigm|_X\leq f\bigm|_X$ означает, что
$l (x')\leq f(x')$  для всех  $x'\in X$.
\item\label{HBiii}   Значения $f(x)$ тождественно равны 
\begin{equation}\label{repsbl}
\inf\left\{\sum_{k=1}^n t_kf(x_k)\colon
\sum_{k=1}^n t_k x_k=x, \; x_k\in X, \; t_k\in \RR^+, \;
n\in \NN:=\{1,2,\dots\}  \right\}
\text{ при всех $x\in X$}.
\end{equation}
\end{enumerate}
\item\label{II}
Следующие  три утверждения попарно  эквивалентны.
\begin{enumerate}[{\bf 1.}]
	\item\label{HBci} Функция $f$ выпуклая\/ {\rm (пишем $f\in \conv \RR^X$)}, т.\,е. удовлетворяет неравенству Йенсена:
\begin{equation}\label{cof}
f\bigl(tx_1+(1-t)x_2\bigr)\leq tf(x_1)+(1-t)f(x_2)
\quad\text{для всех $x_1,x_2\in X$ и $t\in (0,1)$},
\end{equation}
где в данном случае интервал $(0,1)\subset \RR_*^+$ можно заменить на отрезок $[0,1]\subset \RR^+$. 
	
\item\label{HBcii}
Для векторного пространства\/ $\aff \RR^X$ над\/ $\RR$  всех аффинных функций $a \in \RR^X$:
\begin{equation}\label{af}
a\bigl(tx_1+(1-t)x_2\bigr)= ta(x_1)+(1-t)a(x_2)
\quad\text{для всех $x_1,x_2\in X$ и $t\in (0,1)$},
\end{equation} 
где $(0,1)$ можно заменить на $\RR$,
$f$ совпадает со своей  нижней огибающей  по $\aff \RR^X$:
\begin{equation}\label{aenvf}
f(x)=\sup \bigl\{a(x)\colon \, a\in \aff \RR^X, \; a\bigm|_X\leq f\bigm|_X \bigr\}=:
{\lE}_{\aff \RR^X}^f (x)\quad\text{для всех $x\in X$}.
\end{equation} 
\item\label{HBciii} Значения $f(x)$ тождественно равны 
\begin{equation}\label{repco}
\inf\left\{\sum_{k=1}^n t_kf(x_k)\colon
\sum_{k=1}^n t_k x_k=x, \; x_k\in X, \; t_k\in \RR^+, \;
\sum_{k=1}^n t_k=1, \; n\in \NN  \right\}
\quad \text{ при всех $x\in X$}.
\end{equation}
\end{enumerate}
\item\label{III} Если формально дополнить выпуклые конусы $\sbl \RR^X\subset \conv \RR^X$ функцией\footnote{тождественная $-\infty$} $\boldsymbol{-\infty}\colon x\mapsto -\infty$,
$x\in X$, то точные нижние границы  $\inf$ из   \eqref{repsbl} и \eqref{repco} тождественно равны
 соотв.\footnote{сокращение для <<соответственно>>} \begin{enumerate}[{\bf 1)}]
\item\label{III1} наибольшей сублинейной миноранте функции $f$:
\begin{equation}\label{gslm}
{\lE}_{\sbl \RR^X}^f (x):=
\sup \bigl\{ p(x)\colon p\in \sbl \RR^X , \; p\bigm|_X\leq f\bigm|_X\bigr\}=\text{\eqref{repsbl}}, \quad x\in X,
\end{equation} 
где для пустого подмножества $\varnothing\subset \RR^X$ полагаем $\sup \varnothing :=\boldsymbol{-\infty}$,
\item\label{III2} и наибольшей выпуклой миноранте функции $f\in \RR^X$:
\begin{equation}\label{gclm}
{\lE}_{\conv \RR^X}^f (x):=
\sup \bigl\{ v(x)\colon v\in \conv  \RR^X , \; v\bigm|_X\leq f\bigm|_X\bigr\}=\eqref{repco}, \quad x\in X,
\end{equation}
\end{enumerate}
\end{enumerate}
\end{thHB}  
Чаще всего формулируют только часть \ref{I} и только  импликацию \ref{HBi}$\Rightarrow$\ref{HBii}, но обратная \ref{HBii}$\Rightarrow$\ref{HBi} --- элементарное следствие  определений линейной \eqref{lf} и сублинейной {\bf (ph)}+{\bf (sa)} функций. Часть \ref{II} достаточно просто следует из части \ref{I} (см., например,  \cite{W60}, \cite[основная теорема, часть (2)]{KhFE13}).

Эквивалентности \ref{HBi}$\Leftrightarrow$\ref{HBiii}
и \ref{HBci}$\Leftrightarrow$\ref{HBciii} элементарны и в неявной форме содержатся, например, в  \cite[доказательство теоремы 1]{Rev09}. По схеме из \cite{Rev09} и определению 
сублинейности \eqref{ph}--\eqref{sa} без труда выводится \ref{III1}. По аналогии с этим выводом по определению  выпуклости \eqref{cof} легко устанавливается наверняка известная и сформулированная где-нибудь ранее часть \ref{III2}.

В <<Математической энциклопедии>> \cite[Хана--Банаха теорема]{ME} формулировка теоремы Хана--Банаха для векторного пространства некорректна: {\it <<В случае действительного пространства $X$ полунорму можно заменить положительно однородным функционалом, \dots>>}, т.\,е. опущено требование субаддитивности {\bf (sa)}.  Тем не менее, основной ориентир в выборе терминологии и обозначений, где это возможно,  --- именно 
\cite{ME}, а также монография \cite{AK} и первая часть \cite{KhRI} настоящей работы, в которой рассмотрен простой модельный случай  лишь однородных функций {\bf (ph)}+\eqref{ph} и их обобщений.
В различных источниках и у разных авторов терминология  зачастую существенно разнится, поэтому в нашем изложении по возможности все, даже элементарные, определения, понятия и утверждения, встречавшиеся нам в литературе хотя бы раз  в различных смыслах и трактовках, приводятся полностью во избежание  разночтений. Нас в первом приближении  будет интересовать актуальный в теории оптимизации \cite{DELV} и в теории функций \cite{Kh01_1}--\cite{Kh14} случай 
\begin{description}
\item[(sbl)]
{\it {\rm \large(}$f\colon X\to \RR_{+\infty}:=\RR\cup \{+\infty\}$ --- сублинейная функция{\rm \large)}
$\overset{{\rm def}}{\Longleftrightarrow}$
{\rm \large(}выполнено {\bf (ph)} с $\RR_*^+$ и {\bf (sa)}{\rm \large)}, где $+\infty :=\sup \RR$, $t+(+\infty):=+\infty=:(+\infty)+t$ при $t\in \RR_{+\infty}$, $t\cdot (+\infty):=+\infty$ при $t\in \RR_*^+$.}
\end{description}
Такие функции называют ещё {\it гиполинейными\/} (hypolinear \cite{AL}, \cite{DELV}) или {\it расширенными cублинейными\/}
(extended sublinear \cite{S68}, \cite[замечание 1.4]{S08}) функционалами при одном из трёх соглашений

\begin{equation}\label{v0}
f(0)=0 \quad {\Longleftrightarrow}
\quad f(0)\in \RR \quad {\Longleftrightarrow} \quad
\bigl(\text{{\bf (ph)} с $\RR^+$ и $0\cdot (+\infty):=0$}\bigr).
\end{equation}
Примеры функций $\boldsymbol{+\infty}\colon x\mapsto +\infty$  и $x\mapsto \begin{cases}
+\infty, \; &x\in \RR^+,\\
0, \; &x\notin \RR^+,
\end{cases}$ $x\in X$, показывают, что {\bf (sbl)}
$\nRightarrow$ \eqref{v0}. 

В такой постановке задача описания нижней огибающей существенно усложняется и далеко не тривиальна. 
Представляется уместной пространная цитата из 
\cite[1. Введение]{DELV}:
``A surprising fact occurs when, as requested in many constraint optimization problem, $p$ is allowed to take the value $+\infty$. When the dimension of the underlying space $X$ is infinite, S. Simons provided (see\footnote{Наше дополнение к цитате --- см. также \cite[замечание 1.4]{S08}, \cite[(1.6)]{AL}.} the paragraph of the article \cite{S68} entitled ‘Counterexample
to 4’, at p.114) a highly counterintuitive example: a hypolinear function (that is a convex ph-function) 
$p \colon X \to \RR_{+\infty}$ such that the inequality $g \leq p$ is false for all the linear mappings $g \colon  X \to \RR$.
Arguably, the better illustration of the difficulty to address this case is the unusually
large number of flawed Hahn\,--\,Banach type theorems for hypolinear functions which can be
found in the mathematical literature; in the articles,
\cite{AL}, \cite{Za1}--\cite{Za2} the reader can find examples and criticism of as much as seven such incorrect results published between 1969 and 2005.'' Эту цитату можно дополнить еще одним более поздним восьмым некорректным результатом \cite[теорема 2.4]{GPP} 2013 г., где вместо сублинейных 
{\bf (sbl)+\eqref{v0}} рассматриваются противоположные им суперлинейные  функционалы со значениями в $\RR_{-\infty}:=\RR\cup \{-\infty\}$. Этот результат 
по существу используется для доказательства \cite[теоремы 2.7, 2.8, следствия 2.5, 2.9]{GPP}, но опровергается контрпримерами \cite[Counterexample to (3), p.~113]{S68}, \cite[замечание 2.3]{S08}, \cite[пример (1.6)]{AL}.  

\subsection{Топологические версии} По-видимому, первая  {\it топологическая версия\/} теоремы Хана\,--\,Банаха об огибающей для расширенных сублинейных {\bf (sbl)}+\eqref{v0} функционалов была дана  Л.~Хёрмандером \cite{Hor} ещё в 1955 г., хотя, как отмечено в \cite[0.1, теорема 8 и далее]{Tik87}, как минимум, в конечномерном случае $X=\RR^n$ для расширенных выпуклых функционалов она была известна и Г.~Минковскому.    

\begin{thH}[{\rm об огибающей \cite{Hor}, \cite[(3.6), (4.8)]{AL},  \cite[${\bf 2}^{\circ}$]{KR}, \cite[4.7.3(3)]{KK}}] 
Пусть $X$ --- локально выпуклое пространстве, $f\colon X\to \RR_{+\infty}$ или $f=\boldsymbol{-\infty}$ --- тожественная $-\infty$ на $X$; $C(X)\subset \RR^X$ --- векторное пространство над\/ $\RR$ непрерывных функций.
\begin{enumerate}[{\bf I.}]
\item\label{HI} Следующие два утверждения эквивалентны.
\begin{enumerate}[{\bf 1.}]
\item\label{HI1} Функция $f$ сублинейная в смысле\/ {\bf (sbl)} полунепрерывная снизу на $X$ или  $f=\boldsymbol{-\infty}$. 
\item\label{HI2} $f$ совпадает со своей нижней огибающей по выпуклому конусу $C(X)\cap \lin \RR^X$:
\begin{equation}\label{gelc}
f(x)=\sup \bigl\{ l(x)\colon l\in C(X)\cap \lin \RR^X , \; l\bigm|_X\leq f\bigm|_X\bigr\}=:{\lE}_{C(X)\cap \lin \RR^X}^f  \quad\text{для всех  $x\in X$}.
\end{equation}
\end{enumerate}
\item\label{HII} Эквивалентны следующие два утверждения.
\begin{enumerate}[{\bf 1.}]
\item\label{HII1} $f\in (\RR_{+\infty})^X$ --- выпуклая 
в смысле \eqref{cof} полунепрерывная снизу или  $f=\boldsymbol{-\infty}$ на $X$.
\item\label{HII2} $f$ совпадает со своей нижней огибающей по выпуклому конусу $C(X)\cap \lin \RR^X$:
\begin{equation}\label{glc}
f(x)=\sup \bigl\{ a(x)\colon a\in C(X)\cap \aff \RR^X , \; a\bigm|_X\leq f\bigm|_X\bigr\}=:{\lE}_{C(X)\cap \aff \RR^X}^f  \quad\text{для всех  $x\in X$}.
\end{equation}
\end{enumerate}
\end{enumerate}
\end{thH}

<<Индивидульное>> развитие теоремы Хёрмандера  для отдельных $x$  даёт
\begin{thAL}[{\rm об огибающей \cite[теоремы (3.4), (4.7)]{AL}}] Пусть $X$ --- локально выпуклое пространство и $P\subset X$ --- выпуклый конус, $K\subset X$ --- выпуклое множество. 
\begin{enumerate}[{\bf I.}]
\item\label{ALHI} Пусть $f\colon P\to \RR_{+\infty}$ --- гиполинейный функционал на $P$, полунепрерывный снизу в нуле.
Тогда следующие два утверждения эквивалентны.
\begin{enumerate}[{\bf 1.}]
\item\label{ALHI1} $f$ --- полунепрерывная снизу в  $x\in P$. 
\item\label{ALHI2} $f(x)=\sup \bigl\{l(x)\colon
l\in C(X)\cap \lin \RR^X, \; l\bigm|_P\leq f\bigm|_P\bigr\}$.
\end{enumerate}
\item\label{ALHII} Пусть $f\colon K\to \RR_{+\infty}$ --- расширенная выпуклая функция, локально ограниченная снизу.
Тогда следующие два утверждения эквивалентны.
\begin{enumerate}[{\bf 1.}]
\item\label{ALHII1} $f$ --- полунепрерывная снизу в  $x\in P$. 
\item\label{ALHII2} $f(x)=\sup \bigl\{a(x)\colon
a\in C(X)\cap \aff \RR^X, \; a\bigm|_K\leq f\bigm|_K\bigr\}$.
\end{enumerate}
\end{enumerate}
\end{thAL}

Введение топологии позволяет дать и законченную векторную версию теоремы Хана\,--\,Банаха о мажорируемом продолжении линейных функционалов \cite[введение, стр. 127--128]{AL}: {\it для гиполинейного функционала $f$ на векторном пространстве $X$ и для $l_0\in \lin \RR^{X_0}$, удовлетворяющего ограничению $l_0\leq f$ на векторном подпространстве $X_0\subset X$ линейное продолжение $l\in \lin \RR^X$, удовлетворяющее условиям $l=l_0$ на $X_0$ и  $l\leq f$ на $X$, возможно тогда и только тогда, когда функция 
$x\mapsto \inf_{x'\in X_0}\bigl(f(x+x')-l_0(x')\bigr)$, $x\in X$, ограничена снизу в некоторой окрестности нуля в сильнейшей локально выпуклой топологии на $X$, в которой непрерывны все функции из $\lin \RR^X$}.

{\it Теорема Хана\,--\,Банаха\,--\,Канторовича\/} \cite[1.4.13--14]{KK}  в {\it порядковой версии}\/ 
в духе предыдущих теорем об огибающей  формулируется в подразделе \ref{ord_ver}, где она непосредственно потребуется. 

Ниже, в подразделе \ref{lue_sp}, даются возможные общие постановки рассматриваемых в работе задач, мотивированные для нас прежде всего  предшествующими применениями 
в теории функций \cite{Kha92_1}--\cite{Kh14}, \cite{Kh01_0}. При этом мы вынуждены повторить часть основных определений и соглашений из \cite{KhRI}. 

\section{Основные определения и постановки задач}\label{ssp}
\setcounter{equation}{0}

\subsection{Упорядоченные множества. Пополнение}\label{os_c} Пусть ${\rm S}$  ---  (частично) упорядоченное множество \cite{AK} с отношением порядка (рефлексивным, транзитивным, антисимметричным) $\leqslant$, т.\,е. пара $({\rm S},\leqslant )$; 
$\,\geqslant\,$ и $\,>\,$ ---  соотв.  {\it обратные к\/} $\,\leqslant\,$ и {\it строгому порядку\/}  $\,<\;:=\;\leqslant\, \cap \,\neq \,$. 

Пара $({\rm S},\leqslant)$, или множество ${\rm S}$, {\it полное снизу} (соотв. {\it сверху\/}), если для каждого непустого подмножества
${\rm S}_0\subset {\rm S}$ существует  {\it  точная нижняя\/} (соотв. {\it верхняя\/}) {\it граница\/} $\inf {\rm S}_0$ (соотв. $\sup {\rm S}_0$). Множество   ${\rm S}$ {\it полное,\/} 
если  ${\rm S}$ полное и снизу, и сверху.       Подмножество ${\rm S}_0\subset {\rm S}$  {\it ограничено снизу\/} (соотв. {\it сверху\/}), если 
существует элемент  $s_0 \in {\rm S}$, для которого  $s_0\leqslant s$ (соотв. $s\leqslant s_0$) для всех  $s\in {\rm S}_0$.  
Множество  ${\rm S}$  {\it порядково полное} {\it снизу\/} (соотв. {\it сверху\/}) \cite{AK}, \cite{KK}, \cite{KR}, \cite{BV}\footnote{{\it lower\/} ({\it upper\/} resp.) {\it order-complete}}, если для каждого непустого ограниченного снизу (соотв. сверху) подмножества  ${\rm S}_0\subset {\rm S}$ существует $\inf {\rm S}_0\in {\rm S}$ (соотв. $\sup {\rm S}_0\in {\rm S}$). Множество ${\rm S}$  {\it порядково полное\/,} если ${\rm S}$ порядково полное и снизу, и сверху.  

Пусть  ${\rm S}$ --- {\it порядково полное.\/} 
Если   $\inf {\rm S} $ и/или $\sup {\rm S} $ не существуют, то часто удобна и полезна операция  {\it (полу-)пополнения\/}  
порядково полного  ${\rm S} $ до полного  пут\"ем добавления  {\it  \underline{символов}\/}  $ \inf {\rm S} $ и/или $ \sup {\rm S} $, если таких элементов  первоначально  в  ${\rm S}$ нет.
Конкретнее\footnote{В  \cite[1.3.3]{KK} использованы иные обозначения, например, ${\rm S}^{\bullet}$ вместо ${{\rm S}}^{\uparrow}$.}, 
\begin{enumerate}
	\item[{[$\downarrow$]}] ${\rm S}_{\downarrow} :=\{\inf {\rm S}\}\cup {\rm S}$ --- {\it полупополнение\/}, или полурасширение, множества ${\rm S}$ {\it вниз,\/} или влево;
	\item[{[$\uparrow$]}] ${{\rm S}}^{\uparrow}:={\rm S}\cup \{\sup {\rm S}\}$  --- {\it полупополнение\/}, или полурасширение, множества ${\rm S}$  {\it вверх,\/} или вправо;
	\item[{[$\updownarrow$]}] ${{\rm S}^{_\uparrow}_{^\downarrow}}:={\rm S}_{\downarrow}\cup {\rm S}^{\uparrow}$ --- {\it пополнение\/}, или расширение,  множества  ${\rm S}$ в порядковом смысле, 
\end{enumerate}
где порядок  $\leqslant\,$ продолжен  естественным путём на эти пополнения, т.\,е.  $\inf {\rm S}\leqslant s\leqslant \sup {\rm S}$ {\it для всех элементов\/} 
$s\in {{\rm S}^{_\uparrow}_{^\downarrow}}$.  Очевидно, пополнения  $ {{\rm S}^{_\uparrow}_{^\downarrow}}, {{\rm S}_\downarrow}, 
{\rm S}^{\uparrow}$  с таким отношением порядка соотв. полное, полное снизу, полное сверху упорядоченные множества.
Для пустого множества  $\varnothing$ полагаем
\begin{equation}\label{varnoth}
\sup \varnothing :=\inf {\rm S} \quad \text{для $\varnothing\subset {\rm S}_{\downarrow}\subset {{\rm S}^{_\uparrow}_{^\downarrow}}$} ;
\quad \inf \varnothing :=\sup {\rm S} \quad\text{для  $\varnothing\subset {{\rm S}}^{\uparrow}$} \subset {{\rm S}^{_\uparrow}_{^\downarrow}}.
\end{equation}

\subsection{Верхняя и нижняя огибающие}\label{lue_sp}
Для множеств $X,Y$ традиционно через $Y^X$ обозначаем множество всех {\it функций} 
(отображений, операторов, функционалов, форм и проч.)\footnote{
В основном будем использовать термин {\it функция,\/} индифферентный к природе множеств $X,Y$
\cite[Функция]{ME}.}   $f\colon X\rightarrow Y$,  $f\colon x\mapsto f(x)$ или $x\mapsto f(x)$, $x\in X$, определённых на $X$.  

Для $X_0\subset X$ через $f\bigm|_{X_0}$ обозначаем  {\it сужение\/ $f$ на\/} $X_0$. Пишем  $\varphi\bigm|_{X_0} =f\bigm|_{X_0}$ и говорим, что <<{\it $\varphi=f$ на\/ $X_0$}>>,  если $\varphi (x)=f (x)$ {\it для всех\/} $x\in X_0$.  В противном случае $\varphi\bigm|_{X_0} \neq f\bigm|_{X_0}$. Пишем $ \varphi\bigm|_{X_0} \leqslant f\bigm|_{X_0}$  и говорим, что <<{\it $\varphi\leq f$ на\/ $X_0$}>>, или $\varphi$ {\it минорирует\/} $f$,
или $f$ {\it мажорирует\/} $\varphi$, {\it на\/}  $X_0$, если $\varphi (x) \leqslant  f(x)$ {\it для всех\/ $x\in X_0$}.  
Пусть $Y=\rm  {\rm S}^{_\uparrow}_{^\downarrow}$  --- пополнение  порядково полного $({\rm S}, \leqslant\,)$.  
Отношение $f\bigm|_{X}\leqslant \varphi\bigm|_{X}$  определяет {\it отношение поточечного порядка на множестве функций\/}
$\bigl({\rm S}^{_\uparrow}_{^\downarrow}\bigr)^X$. Очевидно, множество $\bigl({\rm S}^{_\uparrow}_{^\downarrow}\bigr)^X$ с отношением поточечного порядка, обозначаемого тем же символом $\,\leqslant\,$, {\it полное,\/} а именно:
{\it для произвольного   $F \subset \bigl({\rm S}^{_\uparrow}_{^\downarrow}\bigr)^X$ всегда существуют функции 
	\begin{equation*}\label{siF} 
	\sup  F\colon x\mapsto \sup_{f\in  F}f(x)\in {{\rm S}^{_\uparrow}_{^\downarrow}}, \quad
	\inf F\colon x\mapsto \inf_{f\in F}f(x)\in {{\rm S}^{_\uparrow}_{^\downarrow}}, \quad x\in X,
	\end{equation*}}
когда на  $X$ рассматриваются и постоянные  функции
\begin{equation}\label{conf}
\inf{\bigl({\rm S}^{_\uparrow}_{^\downarrow}\bigr)^X} 			\colon x\mapsto \inf S\in 	{\rm S}^{_\uparrow}_{^\downarrow}, \quad 
\sup {\bigl({\rm S}^{_\uparrow}_{^\downarrow}\bigr)^X} \colon x\mapsto \sup S \in { {\rm S}^{_\uparrow}_{^\downarrow}}	, \quad x\in X,
\end{equation}
а для  $\varnothing\subset \bigl({\rm S}^{_\uparrow}_{^\downarrow}\bigr)^X$  
в соответствии с соглашением  \eqref{varnoth} определены точные  границы
\begin{equation}\label{df:pmi}
\sup \varnothing :=\inf{\bigl({\rm S}^{_\uparrow}_{^\downarrow}\bigr)^X} \in 
\bigl({\rm S}^{_\uparrow}_{^\downarrow}\bigr)^X, \quad 
\inf \varnothing :=\sup {\bigl({\rm S}^{_\uparrow}_{^\downarrow}\bigr)^X} \in 
\bigl({\rm S}^{_\uparrow}_{^\downarrow}\bigr)^X .
\end{equation}

\begin{definition}\label{df:c}
Пусть\/  $({\rm S},\leqslant\,)$  порядково полное,   $f\in \bigl({\rm S}^{_\uparrow}_{^\downarrow}\bigr)^X$,  $X_0\subset X$,  $L \subset \bigl({\rm S}^{_\uparrow}_{^\downarrow}\bigr)^{X_0}$. {\it Нижнюю\/} (соотв. {\it верхнюю\/})  {\it ${L}$-огибающую\/}, или {\it огибающую по\/} ${L}$  {\it  для\/}  $f$ 
{\it на\/ $X_0$}, определяем как функцию
\begin{subequations}\label{luE}
\begin{align}	 
{\lE}_{{L}}^f&\colon	x \mapsto \sup  \bigl\{ {l} (x) \colon  {L}  \ni {l}\bigm|_{X_0}\leqslant f\bigm|_{X_0}\bigr\}\in 
{{\rm S}^{_\uparrow}_{^\downarrow}}, \quad x\in {X_0}  
\tag{\ref{luE}l}\label{luEl}\\
\Bigl(\text{соотв. }
{\uE}^{{L}}_f&\colon x\mapsto \inf \bigl\{ {l} (x) \colon f\bigm|_{X_0} \leqslant {l}\bigm|_{X_0} \in  {L}  \bigr\}\in {{\rm S}^{_\uparrow}_{^\downarrow}}, \quad x\in X_0{\Bigr)}.
\tag{\ref{luE}u}\label{luEu}
\end{align}
\end{subequations}

Очевидно, всегда
\begin{equation}\label{lxu}
{\lE}_{{L}}^f\bigm|_{X_0}\leq f\bigm|_{X_0}\leq {\uE}^{{L}}_f\bigm|_{X_0}.
\end{equation}
\end{definition}
Функция $f\colon X\rightarrow Y$ с упорядоченными  $(X, \leq )$, $(Y, \leq )$ {\it возрастающая} на $X$, если для любых $x_1,x_2\in X$ из $x_1\leq  x_2$ следует $f(x_1)\leq  f(x_2)$. При этом, если $F\subset Y^X$, то множество всех возрастающих функций $f\in F$ обозначаем как $\incr F$. Функция $f\in Y^X$ {\it строго возрастающая\/} на $X$, если для любых $x_1,x_2\in X$ из $x_1< x_2$ следует $f(x_1)< f(x_2)$. 
Аналогично для убывания. Функция  (строго) возрастающая или  убывающая --- ({\it строго\/}) {\it монотонная.} Функции 
$f\mapsto {\lE}_{{L}}^f$ и $f\mapsto {\uE}^{{L}}_f$, $f\in \bigl({\rm S}^{_\uparrow}_{^\downarrow}\bigr)^X$, 
возрастающие на $\bigl({\rm S}^{_\uparrow}_{^\downarrow}\bigr)^X$, но, вообще говоря, не строго возрастающие. 

\subsection{Общие постановки задач}\label{stpr}
В  приложениях роль   $X$ из определения \ref{df:c}  часто играет некоторый  класс функций и ${\rm S}=\mathbb R$ \cite{Kh01_2}--\cite{Kh10}.
Основные проблемы, диктуемые определением \ref{df:c} в свете теорем Хана--Банаха, Хёрмандера, Энгера\,--\,Лембке, ---  

\begin{problem}\label{pr1} 
Описать функции $f\in \bigl({\rm S}^{_\uparrow}_{^\downarrow}\bigr)^X$,   равные своей нижней (верхней)\/   ${L}$-огибающей на $X_0$. 
\end{problem}

\begin{problem}\label{pr2}
Указать метод(ы)  конструктивного построения  ${L}$-огибающих\/ 
для $f\in \bigl({\rm S}^{_\uparrow}_{^\downarrow}\bigr)^X$.
\end{problem}
Именно эти задачи \ref{pr1}, \ref{pr2} решают для определенных классов ${L}$ и функций $f$ результаты подраздела~\ref{ist} --- см. и ср.  
\eqref{lenvf}, \eqref{aenvf}, \eqref{gslm}, \eqref{gclm}, \eqref{gelc}, \eqref{glc} и  пп. \ref{ALHI2} и \ref{ALHII2}
теоремы Энгера\,--\,Лембке.

В связи с приложениями не меньший, а для некоторых применений  (см., например, \cite{Kh96}--\cite{Kh14}) даже  б\'ольший интерес, чем задача \ref{pr1},   представляет собой менее требовательная

\begin{problem}\label{pr_3} 
Описать те $f\in \bigl({\rm S}^{_\uparrow}_{^\downarrow}\bigr)^X$, для которых\/  
${\lE}_{{L}}^f  \neq   \inf{\bigl({\rm S}^{_\uparrow}_{^\downarrow}\bigr)^X}$    или\/ 
${\uE}^{{L}}_f  \neq   \sup{\bigl({\rm S}^{_\uparrow}_{^\downarrow}\bigr)^X}$. 
\end{problem}
Если порядково полное  ${\rm S}$ изначально дополнительно снабжено какими-либо алгебраическими операциями, согласованными с отношением порядка $\,\leqslant\,$, то продолжения этих операций на пополнения, или (полу)расширения, вообще говоря, с ограничениями, могут определяется в каждом конкретном случае в зависимости от  проблематики  (см. и ср. \cite[\S~4]{L}, \cite[1.3.1]{KK}).   

Для пары добавляемых символов $\inf {\rm S}\notin {\rm S}$ и/или $\sup {\rm S}\notin {\rm S}$ часто, в особенности, если на ${\rm S}$ используется  бинарная операция в аддитивной форме как основная, по аналогии с расширениями вещественной прямой  $\mathbb  R$ (c операцией сложения) вверх и/или вниз, используются соотв. обозначения 
$-\infty:=\inf {\rm S}$ и/или $+\infty:=\sup {\rm S}$,  а также ряд других в иных ситуациях \cite[1.2]{KhRI}.

Для множества $S$ через ${\card} S$ обозначаем число элементов в  $S$ или мощность множества $S$.

Векторные пространства рассматриваются над полем $\RR$, если не оговорено противное.
Нулевой вектор различных векторных пространств обозначаем одним и  тем же символом $0$.
Задачи \ref{pr1}--\ref{pr_3} будут решаться на проективных пределах векторных решёток по аналогии с \cite{Kh96}--\cite{Kh10}, \cite{KKh09} для функций (функционалов) со значениями в пространстве Канторовича ${\KK}$ и в ${\KK}_{\pm\infty}$.

Пусть $X,Y$ --- упорядоченные векторные пространства. Через $X^+$ обозначаем множество всех положительных векторов в $X$. В частности,  для подмножества $F\subset Y^X$ тогда $F^+$ --- подмножество всех {\it положительных\/} функций $f$ из $F$, для которых по определению  $f(X^+)\subset Y^+$.
Порядково полная векторная решётка $\mathbb K$, $\card \KK>1$, --- {\it пространство Канторовича\/} \cite{AK}, \cite{KK}, \cite{KR}. При  этом $\KK$ не полное ни снизу, ни сверху и\footnote{Здесь и далее ссылки над знаками бинарных отношений означают, что эти соотношения как-то связаны с приведёнными ссылками, например, вытекают из них.} 
\begin{equation}\label{Kinfty}
\begin{split}
\inf {\KK}:=-\infty, \quad \sup {\KK}:=+\infty,\quad
\sup \varnothing &\overset{\eqref{varnoth}}{:=}-\infty, \quad \inf \varnothing \overset{\eqref{varnoth}}{:=}+\infty \quad \text{для $\varnothing \subset {\KK}$}; 
\\
-(-\infty):=+\infty, \quad -(+\infty):=-\infty,
\quad &-\infty\leq k\leq +\infty\quad\text{для всех $k\in \KK_{\pm \infty}$};
\\
{\KK}_{-\infty}:={\KK}_{\downarrow},\quad {\KK}_{+\infty}&:={\KK}^{\uparrow},
\quad {\KK}_{\pm\infty}:={\KK}_{\downarrow}^{\uparrow}; \\
(-\infty)+k:=:k+(-\infty)&:=-\infty, 
\text{ для всех $k\in {\KK}_{-\infty}$},\\ 
(+\infty)+k:=:k+(+\infty)&:=+\infty \text{ для всех $k\in {\KK}_{+\infty}$};
\\
t(-\infty):=:(-t)(+\infty):=-\infty, 
\quad t(+\infty)&:=:(-t)(-\infty):=+\infty
\text{ для всех $t\in \RR_*^+$};
\\
{\KK}_{-\infty}^S:=({\KK}_{-\infty})^S,\quad 
{\KK}_{+\infty}^S&:=({\KK}_{+\infty})^S,\quad
{\KK}_{\pm \infty}^S:=({\KK}_{\pm \infty})^S; 
\\
\inf {\KK}_{\pm \infty}^S:=:\inf {\KK}_{-\infty}^S
\overset{\eqref{conf}}{:=}\boldsymbol{-\infty},
\quad &\sup{\KK}_{\pm \infty}^S:=:\sup{\KK}_{+\infty}^S
\overset{\eqref{conf}}{:=}\boldsymbol{+\infty}
\quad \text{при $S\neq \varnothing$};\\
\sup \varnothing \overset{\eqref{df:pmi}}{:=}\boldsymbol{-\infty}, \quad \inf \varnothing \overset{\eqref{df:pmi}}{:=}\boldsymbol{+\infty} \quad &\text{для пустого подмножества $\varnothing \subset {\KK}^S\subset ({\KK}_{\pm \infty})^S$}.
\end{split}
\end{equation}

\section{Проективный предел упорядоченных векторных пространств}\label{pros}
\setcounter{equation}{0}

Для векторных пространств $X,Y$ через $\lin Y^X$ обозначаем векторное пространство всех линейных функций $l$ из $Y^X$, удовлетворяющих \eqref{lf}. В случае упорядоченных векторных пространств $X,Y$ полагаем
$\lin^+Y^X:=(\lin Y^X)^+$. Очевидно, в этом случае $\lin^+Y^X=\incr \lin Y^X$.

Пусть $(X_n)_{n\in \NN}$ ---
последовательность упорядоченных векторных пространств  над $\RR$ и $\leq_n$ --- отношение порядка  на $X_n$. Проекцию вектора $x=(x_n)$ из векторного пространства-произведения $\prod_{n} X_n$ на пространство $X_n$ обозначаем через ${\pr}_n x=x_n$.
На $\prod_{n} X_n$ вводится естественный порядок ${\leqslant}\,$,
а именно: $x=(x_n) \leqslant x'=(x_n')$ в $\prod_{n} X_n$, если и только если $x_n\leq_n x_n'$ для всех $n\in{\NN}$.

\begin{definition}\label{dfprp}
Пусть $p_n \in \lin^+ X_n^{X_{n+1}}$. Векторное подпространство $X$ в $\prod_{n} X_n$, c условием
\begin{equation}\label{dfX}
(x\in X)\;\overset{\rm def}{\Longleftrightarrow}\; \bigl({\pr}_n x = p_n ({\pr}_{n+1} x)\quad\text{для всех $n\in \NN$}\bigr), 
\end{equation}
снабженное индуцированным с  $\prod_{n} X_n$ отношением порядка  $\leqslant\,$, называем {\it проективным пределом последовательности $(X_n)$ относительно $(p_n)$\/} и обозначаем
$X=\projlim_n X_np_n$. 
\end{definition}
В записи $\projlim_n$ нижний индекс $n$ часто опускаем. Из \eqref{dfX} сразу следует, что 
\begin{equation}\label{1_2}
p_n({\pr}_{n+1}B)={\pr}_nB\quad\text{для любого подмножества $B\subset X=\projlim X_np_n$}. 
\end{equation}

\subsection{Проективный предел векторных решёток}\label{prlat} Упорядоченное векторное пространство --- {\it векторная решётка\/}, если для всякой пары векторов в нём существует точная  верхняя граница $\sup$, а значит, и точная нижняя граница $\inf$. Операцию $\sup$ в упорядоченном множестве $S$ для её идентификации обозначаем как $S\text{-}\sup$. В частности, для  $(X_n, \leq_n)$ наряду с обозначением  $X_n$-$\sup$ иногда, для краткости, используем и обозначение $n$-$\sup$. Аналогично для $\inf$.
Проективный предел $X=\projlim X_np_n$ векторных решёток
$X_n$ называем {\it правильным проективным пределом\/}, если 
$p_n$ {\ сохраняют $\sup$ для конечных множеств\/}, т.\,е.
\begin{equation}\label{nsupf}
{X_n\text{-}\sup\,} p_n(B_{n+1})=p_n(X_{n+1}\text{-}\sup\,B_{n+1})\quad\text{при всех $n\in \NN$ и $B_{n+1}\subset X_{n+1}$ с $\card B_{n+1}<\infty$}. 
\end{equation}

\begin{propos}\label{prpr} Для правильного  предела $X=\projlim X_np_n$ векторных решёток $X_n$
\begin{enumerate}[{\rm (i)}]
\item\label{pri} если $B\subset X$ и  $\card {\pr}_nB<\infty$ для каждого $n\in \NN$, то существует  
\begin{equation}\label{supB}
v=X\text{-}\sup B=\Bigl(\prod_{n} X_n\Bigr)\text{-}\sup B, \quad
\pr_n v=X_n\text{-}\sup \pr_n B;
\end{equation}
\item\label{prii} $X$ и  все  ${\pr}_nX$ с индуцированным с $X_n$
отношением порядка ${\leq}_n$ ---  векторные решётки;
\item\label{priii}  для любого $n\in \NN$ сужения  $p_n\bigm|_{{\pr}_{n+1} X}$ сохраняют\/ $\sup$ для конечных множеств.
\end{enumerate}
\end{propos}
\begin{proof} \eqref{pri}. Положим
$v=(v_n)\in \prod_{n} X_n$, где $v_n=X_n\text{-}\sup \pr_n B$.
Очевидно, $v=\bigl(\prod_{n} X_n\bigr)\text{-}\sup B$. Поскольку 
$p_n$ сохраняют $\sup$, используя \eqref{nsupf}  и \eqref{1_2},
получаем
\begin{equation*}
p_n(v_{n+1})=p_n(X_{n+1}\text{-}\sup {\pr_{n+1}}B)
\overset{\eqref{nsupf}}{=}X_n\text{-}\sup p_n({\pr}_{n+1}B)
\overset{\eqref{1_2}}{=} X_n\text{-}\sup {\pr}_nB=v_n.
\end{equation*}
Значит для $x=v$ имеет место правая часть \eqref{dfX} и $v\in X$ по  определению \ref{dfprp}. Так как порядок на $X$ индуцирован с  $\prod_{n} X_n$, то $v=X\text{-}\sup B\in X$ и  \eqref{pri}
доказано. 

\eqref{prii}. Пространство-произведение векторных решёток --- векторная решётка \cite{AK}, \cite[гл.II, \S~1]{Bour}. По п. \eqref{pri} ввиду \eqref{supB} с $\card B=2$ его векторное подпространство $X$ --- тоже векторная решётка. 
Пусть теперь $x_n, x'_n\in {\pr}_nX$, то есть
 существуют $x,x'\in X$, для которых ${\pr}_nx=x_n$,
${\pr}_nx'=x'_n$. Поскольку $X$ --- векторная решётка,
существует вектор $v=X\text{-}\sup\, \{x,x'\}\in X$,
для которого $ X_n\text{-}\sup\, \{ x_n, x'_n\} ={\pr}_nv\in {\pr}_nX$. Таким образом, 
\begin{equation}\label{prnXs}
{\pr}_nv=(\pr_nX)\text{-}\sup\, \{x_n,x'_n\}
=X_n\text{-}\sup\, \{x_n,x'_n\}
\end{equation} 
и ${\pr}_nX$ --- векторная решётка, что завершает доказательство п.~\eqref{prii}. 

\eqref{priii}. Из \eqref{prnXs} для  $x_n,x_n'\in \pr_n X$ при $n\in \NN\setminus \{1\}$ имеем
\begin{multline*}
 p_{n-1}\bigm|_{\pr_n X}\bigl((\pr_nX) \text{-}\sup\, \{ x_n, x'_n\} \bigr)
\overset{\eqref{prnXs}}{=}p_{n-1}\bigl(X_n\text{-}\sup\,\{ x_n, x'_n\}\bigr)\\
\overset{\eqref{nsupf}}{=} X_{n-1}\text{-}\sup\,\bigl\{ p_{n-1}( x_n),\, p_{n-1}(x'_n)\bigr\}
\overset{\eqref{prnXs}}{=}
(\pr_{n-1}X)\text{-}\sup\,\Bigl\{ p_{n-1}\bigm|_{\pr_nX}( x_n),\, p_{n-1}\bigm|_{\pr_nX}( x'_n)\Bigr\}.
\end{multline*}
Отсюда все  $p_{n-1}\bigm|_{\pr_n X}$ при $n\in \NN\setminus \{1\}$ сохраняют $\sup$ для конечных множеств из $\pr_nX$.
\end{proof}

Предел $X=\projlim X_np_n$ называем {\it
приведённым\/}, если  ${\pr}_nX=X_n$ для каждого $n\in \NN$ (ср. с \cite[гл.IV, \S~4]{Schef}). 
Не умаляя общности общности,  любой проективный
предел можно считать приведённым, поскольку  
 $X=\projlim X_np_n\overset{\eqref{1_2}}{=}\projlim\, ({\pr}_nX)\,\bigl(p_n\bigm|_{X_{n+1}}\bigr)$. 
При этом все отображения $p_n\bigm|_{X_{n+1}}\colon \pr_{n+1}X\to \pr_nX$ в силу \eqref{1_2} 
сюръективные. Более того, если первоначально рассматривался
правильный проективный предел векторных решёток, то по предложению \ref{prpr} после перехода к приведённому проективному пределу
мы по-прежнему имеем дело с правильным проективным пределом векторных решёток $\pr_nX$ относительно $p_n\bigm|_{X_{n+1}}$.

\begin{remark}\label{prl_subs} Векторные решётки $X$ и $Y$ решёточно изоморфны, если существует линейная биекция из $X$ на $Y$, сохраняющая порядок. Каждый приведённый правильный проективный предел $X=\projlim X_np_n$  векторных решёток для произвольной возрастающей подпоследовательности $(n_k)_{k\in \NN}\subset \NN$ решёточно
изоморфен такому же пределу $X=\projlim_k X_{n_k}q_k$ относительно 
$q_k =p_{n_{k}}\circ p_{n_{k}+1}\circ
\dots \circ p_{n_{k+1}-1}\; :\; X_{n_{k+1}} \to  X_{n_k}$, где $\circ$ --- операция суперпозиции.
\end{remark}

\subsection{Примеры}\label{examp} 
\paragraph{\bf 1. Векторные решётки-произведения и решётки из последовательностей} Рассмотрим  $(X_k, \le_k')_{k\in \NN}$ --- последовательность векторных решеток c отношением порядка $\leq_k'$. Тогда каждое конечное произведение $\prod_{k=1}^n X_k$ --- векторная решетка при  естественном отношении порядка $(x_1,\dots, x_n)\leq_n (x'_1,\dots, x'_n)$, $x_k\in X_k$, при $x_k\leq_k' x_k'$ для всех  $k=1, \dots, n$. Для 
\begin{equation}\label{prodprp}
p_n\colon \prod_{k=1}^{n+1}X_k \to \prod_{k=1}^{n}X_k, \quad
p_n\colon (x_1,\dots, x_n, x_{n+1})\mapsto (x_1,\dots, x_n),
\end{equation} 
пространство-произведение $\prod_{k\in \NN} X_k$ векторных решеток решёточно изоморфно  приведённому правильному проективному пределу $\projlim \, (\prod_{k=1}^nX_n)p_n$ векторных решёток $\prod_{k=1}^nX_n$ относительно $p_n$. Таким образом, каждое пространство-произведение счетного числа векторных решёток можно рассматривать как приведённый правильный проективный предел.

В частном случае совпадающих при всех $k$ векторных решёток $(X_k, \leq_k') =(X,\leq')$, $k\in \NN$, пространство-произведение $X^{\NN}$ --- это векторная решётка всех последовательностей $(x_k)_{k\in \NN}$ с векторами $x_k\in X$,  которую можно решёточно изоморфно отождествить с приведённым правильным проективным пределом $\projlim \, (X^n)p_n$    
относительно $p_n$ из \eqref{prodprp}. 

\paragraph{\bf 2. Векторные решётки функций и мер} Пусть $D$ --- топологическое хаусдорфово пространство. Для $S\subset D$ через $\clos S$, $\Int S$ и $\partial S$ обозначаем соотв. {\it замыкание, внутренность} и {\it границу\/} $S$ в $D$. Подмножество $S$ {\it предкомпактно в\/} $D$ (пишем $S\Subset D$), если $\clos S$ --- компакт в $D$. 

Предполагаем, что для $D$  {\it существует  исчерпание последовательностью\/ $(D_n)_{n\in \NN}$  предкомпактных открытых  подмножеств из\/}  $D$:
\begin{equation}\label{Dnexh}
\Int D_n=D_n\Subset D_{n+1},\; n\in \NN, \quad \bigcup_{n\in \NN} D_n=D.
\end{equation}
Различные {\it достаточные условия для существования исчерпания\/} \eqref{Dnexh} давно известны и не раз отмечались и использовались \cite[\S~1, перед (1.3)]{Kh01_1}, \cite[лемма 2.1]{GPP}:
\begin{enumerate}[{(a)}]
\item\label{aD} $D$ локально компактное  и  представимо в виде объединения  счётного числа компактных подмножеств, т.\,е. $\sigma$-компактно \cite[следствие 2.77]{AB}. 
\item\label{bD} $D$ удовлетворяет первой аксиоме счётности и  полукомпактное (hemicompact), т.\,е. существует последовательность компактных подмножеств со свойством: каждое компактное подмножество содержится в некотором компактном подмножестве из этой последовательности. 
\end{enumerate}
По каждому из условий \eqref{aD}$\Leftarrow$\eqref{bD} пространство $D$ нормальное, в силу чего любая непрерывная  функция из  $\RR^{\clos D_n}$ --- сужение какой-либо непрерывной функции из $\RR^D$ на $\clos D_n$. 

\paragraph{\bf 3. Векторная решётка $C(D)$} Для $S\subset D$ через $C(S)\subset \RR^{S}$ обозначаем векторную решётку непрерывных на $S$ функций с отношением поточечного порядка на $S$.  При исчерпании 
\eqref{Dnexh} для $p_n\colon f\mapsto f\bigm|_{\clos D_n}$, $f\in C(\clos D_{n+1})$, определён проективный предел векторных решеток $\projlim_n C(\clos D_n)p_n$, который для нормального $D$ приведённый правильный  и его можно решёточно изоморфно отождествить с $C(D)$.
\paragraph{\bf 4. Векторная решётка $L_{\loc}^p(D,\mu)$} Пусть $p\in \RR_*^+$, $D$ локально компактное с исчерпанием \eqref{Dnexh} и $\mu$ --- вещественная мера Радона на $D$; $L^p \bigl(\clos D_n, \mu\bigm|_{\clos D_n}\bigr)$ --- векторная решётка над $\RR$ функций-классов эквивалентности на $\clos D_n$ по отношению $=$ $\mu$-п.\,в. и отношением порядка $\leq$ $\mu$-п.\,в. поточечно \cite[гл.~IV, \S~2, 6]{Bour}. Тогда определён приведённый правильный проективный предел $L^p_{\loc}(D,\mu )=\projlim_n L^p\bigl(\clos D_n,\mu\bigm|_{\clos D_n}\bigr)p_n$ над $\RR$ локально интегрируемых в $p$-ой степени модулем по мере Радона $\mu$ функций-классов эквивалентности, где $p_n\colon f\mapsto f\bigm|_{\clos D_n}$, $f\in L^p \bigl(\clos D_{n+1}, \mu\bigm|_{\clos D_{n+1}}\bigr)$.
\paragraph{\bf 5. Векторная решётка мер ${\Meas}(D)$} Пусть $D$ локально компактное с исчерпанием \eqref{Dnexh}, ${\Meas}(\clos D_n)$, $n\in \NN$, --- векторные решётки вещественных мер Радона на компактах $\clos D_n$, $p_n\colon \mu\mapsto \mu\bigm|_{D_n}$ --- сужение на $\clos D_n$ меры $\mu \in {\Meas}(\clos D_{n+1})$ \cite[гл.~III]{Bour}. Тогда векторная решетка всех вещественных мер Радона ${\Meas}(D)$ на $D$ \cite[гл.~III, \S~2]{Bour} можно рассматривать как приведённый проективный предел $\projlim_n {\Meas}(\clos D_n)p_n$.

\subsection{Положительные аддитивные функции на проективном пределе}\label{+lf}

\begin{definition}\label{dfadd} Пусть $(X,+)$ и $(Y,+)$ ---  полугруппы по сложению. Функция $a\colon X\to Y$ {\it аддитивная,\/} если $a(x_1+ x_2)=a(x_1)+ a(x_2)$ для любых $x_1, x_2 \in X$. Множество всех таких аддитивных функций обозначаем как $\add Y^X$. Когда $X,Y$ --- упорядоченные векторные пространства, то полагаем $\add^+Y^X := (\add Y^X)^+$, $\add^{\reg} Y^X:=\add^+Y^X-\add^+Y^X$ (ср. с \cite[I.1.4.III]{AK}). 
\end{definition}
Элементарно устанавливается
\begin{propos}\label{pr_addgr}
Пусть $(X,+)$ и $(Y,+)$ --- аддитивные упорядоченные группы. Тогда для любой функции $a\in \add Y^X$ имеем $a(0)=0\in Y$, $a(-x)=-a(x)$. Когда $X,Y$ --- упорядоченные векторные пространства, то $\add^+Y^X := (\add Y^X)^+=\incr \add Y^X$.
\end{propos}

\begin{propos}[{\rm ср. с \cite[предложение 2.1]{Kh01_1}}]\label{Tf+}
Пусть $X=\projlim X_np_n$ --- приведённый правильный  проективный предел векторных решёток $(X_n,\leq_n)$.
Для любой функции  $a\in \add^+\RR^X$ существуют номер $n_a\in \NN$ и последовательность  функций $a_n\in \add^+\RR^{X_n}$, $n\geqslant n_a$, для которых  
\begin{equation}\label{l+}
a(x)=a_n({\pr}_nx)\quad\text{для всех $x\in X$}
\end{equation}
при всех $n\geq n_a$. Обратно, если $a_n\in \add^+\RR^{X_n}$, 
то \eqref{l+} определяет $a\in \add^+\RR^X$.

Очевидно, $\add^+\RR^X$, $\add^+\RR^{X_n}$ здесь можно заменить соотв. на $\add^{\reg}\RR^X$, $\add^{\reg}\RR^{X_n}$.
\end{propos} 
\begin{proof} Часть {\it <<Обратно, \dots>>\/} очевидна, поскольку суперпозиция $a\circ \pr_n$  аддитивной положительной функции $a$ с линейной положительной функцией  $\pr_n$ аддитивна и положительна.  

Покажем сначала, что найдётся номер $n\in \NN$, для которого  $a_n(x)=a_n(x')$, как только ${\pr}_nx={\pr}_nx'$, т.\,е. $a(x-x')= 0$ при ${\pr}_n(x-x')=0$. Предположим, что такого номера $n$ не существует. Тогда по предложению \ref{pr_addgr} для каждого $m\in \NN$ найдётся элемент $ x^{(m)}\in X$, для которого
\begin{equation}\label{2_3}
\quad {\pr}_m x^{(m)}=0, \quad a( x^{(m)})> 0\quad\text{при всех $m\in \NN$}.
\end{equation}
Из последнего строгого неравенства, переходя от $x^{(m)}$ к $k_mx^{(m)}$ с достаточно быстро растущей последовательностью $(k_m)\subset \NN$, по аксиоме Архимеда можно добиться того, что последовательность\/  
$\bigl(a( k_mx^{(m)})\bigr)_{m\in \NN}\subset \RR$, 
где $a( k_mx^{(m)})=k_m a(x^{(m)})$, {\it не ограничена сверху\/} в $\RR$ и  ${\pr}_m (k_mx^{(m)})=k_m{\pr}_m x^{(m)}\overset{\eqref{2_3}}{=}0$.
В силу этих равенств  и линейности  $p_1, p_2, \dots ,p_{m-1}$ для каждого $ k_mx^{(m)}$ имеем ${\pr}_n(k_mx^{(m)})=0$ при $n\leqslant m$, то есть при фиксированном $n$ лишь конечное число элементов последовательности $({\pr}_n k_mx^{(m)})_{m\in \NN}$ отлично от нуля в векторной решётке $X_n$. Следовательно, по предложению 
\ref{prpr}\eqref{pri} последовательность  $(k_mx^{(m)})_{m\in \NN}$  {\it ограничена сверху\/} в $X$. 
По предложению \ref{pr_addgr} функция $a\in \add^+\RR^X$  возрастающая. Отсюда и последовательность $\bigl(a( k_mx^{(m)})\bigr)_{m\in \NN}\subset \RR$ {\it ограничена сверху\/} в $\KK$.  Противоречие!

Пусть $n$ --- номер, для которого  $a(x)=a(x')$, как только ${\pr}_nx={\pr}_nx'$. Тогда $a_n$ однозначно определяется по 
\eqref{l+}. Аддитивность $a_n$ следует из  линейности ${\pr}_n$
и аддитивности $a$. Кроме того, если $x_n\in {\pr}_nX$, $0{\leq}_nx_n$,  то существует $x\in X$, для которого ${\pr}_nx=x_n$. По предложению \ref{prpr}\eqref{pri}
имеем $x_n\overset{\eqref{supB}}{=}{\pr}_n{x^+}$, где $x^+=\sup \{ x,0\} \geqslant 0$ в $X$, и $a_n(x_n)=a_n({\pr}_n{x^+})=a(x^+)\geqslant 0$. Установлена положительность $a_n$. Теперь можно положить $n_a=n$.
Остальные $a_m\in \add^+(\RR^{X_m})$ однозначно определяются при $m\geqslant n$ по правилу $a_m(x_m):=a_n\bigl(p_n\circ p_{n+1}\circ \dots \circ p_{m-1}(x_m)\bigr)$,  $m\geqslant n\geqslant n_a$.
\end{proof}

\begin{remark}\label{r_supp} Предложение \ref{Tf+} обобщает тот известный факт, что линейные (положительные) функционалы на $C(D)$ и $L_{\loc}^p(D,\mu)$ из подраздела \ref{examp}
имеют компактный носитель.
\end{remark}

\begin{remark}\label{r_arch} Предложение \ref{Tf+} остается в силе и для $a\in \add^+ {\mathbb G}^X$, когда $(\mathbb G,+)$ --- 
аддитивная линейно упорядоченная  архимедова группа. Доказательство почти дословно такое же.
\end{remark}

\begin{corollary}[{\rm \cite[предложение 2.1]{Kh01_1}}]\label{Tf+c}
В условиях предложения\/ {\rm \ref{Tf+}}
для каждой функции  $l\in \lin^+\RR^X$ существуют $n_l\in \NN$ и последовательность  $l_n\in \lin^+\RR^{X_n}$, $n\geqslant n_l$, для которых  
\begin{equation}\label{ll+}
l(x)=l_n({\pr}_nx)\quad\text{для всех $x\in X$}
\end{equation}
при всех $n\geq n_l$. Обратно, если $l_n\in \lin^+\RR^{X_n}$, 
то \eqref{ll+} определяет $l\in \lin^+\RR^X$.

Здесь $\lin^+\RR^X$, $\lin^+\RR^{X_n}$ можно заменить соотв. на $\lin^{\reg}\RR^X:=\lin^+\RR^X-\lin^+\RR^X$, $\lin^{\reg}\RR^{X_n}$.
\end{corollary} 
\begin{proof} В дополнение к предложению \ref{Tf+} достаточно добавить, что суперпозиция двух линейных функций линейна.
\end{proof}

\begin{definition}\label{df_aff} Пусть $X,Y$ --- векторные пространства. Функция  $a\colon X\to Y$ {\it аффинная,\/} если 
имеет место \eqref{af}. Множество всех аффинных функций из $Y^X$ обозначаем как $\aff Y^X$. 
Когда $X,Y$ --- упорядоченные векторные пространства, то $\aff^+Y^X:=(\aff Y^X)^+$.
\end{definition}
\begin{propos}\label{afflin}
 $a\in \aff X^Y$, если и только если имеем единственное  представление 
\begin{equation}\label{affl}
a=l_a+\boldsymbol{1} a(0), \quad\text{где $l_a\in \lin Y^X$ и  
$\boldsymbol{1}\colon x\mapsto 1\in \RR$ --- тождественная	единица на $X$}.
\end{equation}  
Пусть $X,Y$ --- упорядоченные векторные пространства. Тогда
\begin{enumerate}[{\rm (i)}]
\item\label{alin} $a\in \incr \aff Y^X$, если и только если в представлении  \eqref{affl} $l_a\in \lin^+ Y^X$;
\item\label{ali} если в  \eqref{affl} одновременно $l_a\in \lin^+ Y^X$ и  $a(0)\in Y^+$, то $a\in \aff^+ Y^X$;
\item\label{aliK} если $a\in \aff^+ \RR^X$, то одновременно $a(0)\in \RR^+$ и $l_a\in \lin^+ \RR^X$.
\end{enumerate}
\end{propos}
Нам не удалось найти доказательства этого тривиального предложения. Поэтому даём полное
\begin{proof}
Достаточность очевидна. Единственность представления \eqref{affl} следует из $l_a(0)=0\in Y$. Для доказательства необходимости положим  $l_a:=a-\boldsymbol{1} a(0)$. Тогда $l_a(0)=0\in Y$ и $l_a$ --- аффинная функция. Следовательно,  при любом $t\in [0,1]$ имеем $l_a(tx)=
a\bigl(tx+(1-t)0\bigr)-a(0)=ta(x)+(1-t)a(0)-a(0)=
t\bigl(a(x)-a(0)\bigr)=tl_a(x)$. При $\RR\ni t> 1$ имеем   
$1/t\in (0,1)$ и $l_a(x)=l_a\bigl((1/t)tx\bigr)=(1/t)l_a(tx)$, откуда $l_a(tx)=tl_a(x)$ при $t\in \RR^+$, т.\,е. $l_a$ --- положительно однородная функция. При этом $l_a(x)+l_a(-x)=\frac12l_a(2x)+
\frac12 l_a(-2x)=l_a(0)=0$ и $l_a(x)=-l_a(-x)$ для любого $x\in X$. Это доказывает однородность $l_a$ на $X$. Аддитивность $l_a$ следует из цепочки равенств $l_a(x_1+x_2)= l_a\bigl(\frac12 2x_1+\frac12 2x_2\bigr)=\frac12 l_a(2x_1)+\frac12 l_a(2x_2)=
\frac12 2l_a(x_1)+\frac12 2l_a(x_2)=l_a(x_1)+l_a(x_2)$.

Утверждения \eqref{alin} и \eqref{ali} 
 очевидны. Если $a\in \aff^+Y^X$, то, очевидно, $a(0)\in Y^+$.

\eqref{aliK}. Пусть теперь $Y=\RR$ и $a\in \aff^+\RR^X$. Тогда $a(0)\in \RR^+$. Предположим, что $l_a(x)<0$ для некоторого $x\in X^+$. Тогда для любого $n\in \NN$ имеем 
$0\leq a(nx)=l_a(nx)+a(0)=nl_a(x)+a(0)$ и $0>nl_a(x)>-a(0)$, что невозможно по аксиоме Архимеда для $\RR$. 
\end{proof}

\begin{corollary}\label{cor_aff} В условиях предложения\/ {\rm \ref{Tf+}}
для любых $a\in \incr \aff \RR^X$  или $a\in \aff^+\RR^X$ существуют $n_l\in \NN$ и последовательность  $l_n\in \lin^+\RR^{X_n}$, $n\geqslant n_l$, для которых  
\begin{equation}\label{la+}
a(x)=l_n({\pr}_nx)+a_0\quad\text{для всех $x\in X$,}
\end{equation}
при всех $n\geq n_l$, где $a_0=a(0)\in \RR$, а при $a\in \aff^+\RR^X$ ещё и $a_0\in \RR^+$. Обратно, если $l_n\in \lin^+\RR^{X_n}$, 
то \eqref{la+} при $a_0\in \RR$ определяет $a\in \incr \aff \RR^X$, а при $a_0\in \RR^+$ определяет $a\in \aff^+ \RR^X$.

Очевидно, здесь множества $\incr \aff \RR^X$, $\aff^+\RR^X$, $\lin^+\RR^{X_n}$ можно заменить соотв. на множества $\incr^{\reg}\aff \RR^X:=\incr\aff \RR^X-\incr\aff \RR^X$,
$\aff^{\reg}\RR^X:=\aff^+\RR^X-\aff^+\RR^X$, $\lin^{\reg}\RR^{X_n}$, но без каких-либо дополнительных ограничений на $a_0\in \RR$.
\end{corollary} 
\begin{proof} Следует из сочетания предложений \ref{afflin} и \ref{Tf+}.
\end{proof}

\begin{remark} По общей схеме из \cite[гл.~3]{BD}
 можно определить проективный предел несчётного семейства упорядоченных векторных пространств, как это сделано в \cite{Lu05a}. В то же время нам не известно естественное обобщение утверждений подраздела \ref{+lf}, в частности, следствия \ref{Tf+c}, в таком варианте, хотя ранее и была предпринята такая попытка в \cite[\S~2]{Lu05a}, оказавшаяся неудачной, поскольку некорректна \cite[лемма~2]{Lu05a}.
\end{remark}

\begin{remark}
Если $X$ --- векторная решётка, то по теореме Рисса--Канторовича $\lin^{\reg} \RR^X$ --- векторная решётка и даже пространство Канторовича с положительным конусом $\lin^+ \RR^X$  \cite[теорема 6(3.1)]{AK}. Следствие \ref{Tf+c} означает, что
упорядоченное векторное пространство $\lin^{\reg} \RR^X$, где  $X=\projlim_n X_np_n$  --- приведённый  правильный проективный предел векторных решёток  $X_n$, решёточно изоморфно индуктивному пределу $\injlim_n i_n (\lin^{\reg} \RR^{X_n})$ относительно положительных линейных функций $i_n\colon \lin^{\reg} \RR^{X_n}\rightarrow \lin^{\reg} \RR^{X_{n+1}}$,
сопряжённых к линейным отображениям $p_n$. Не останавливаясь подробно на соответствующих определениях, отметим, что они легко конструируются по общей схеме  \cite[гл.~3, \S~2]{BD}
и по аналогии с индуктивными пределами топологических
векторных пространств \cite{SiS}, \cite[гл.~II и гл.~IV, \S~4]{Schef}.
\end{remark}

\section{Верхняя огибающая функции на проективном пределе}\label{ueL}
\setcounter{equation}{0}

Всюду далее $\KK$ --- пространство Канторовича с отношением порядка $\leq$. При этом отношения порядка и на иных упорядоченных множествах могут  обозначаться тем же символом $\leq$.

\subsection{Случай произвольной функции}\label{ueLar} 
Неравенство из правой части \eqref{lxu} даёт 
\begin{propos}\label{pr3_1}
 Пусть $f\overset{\eqref{Kinfty}}{\in}  \KK_{\pm\infty}^X$.
Тогда
\begin{equation}\label{eq3_1}
 f(x)\leqslant {\uE}^{L}_f(x)\overset{\eqref{luEu}}{=}
\inf \bigl\{ l(x)\colon l\in L, \;
f\bigm|_{X_0}\leqslant l\bigm|_{X_0} \bigr\}\quad
\text{для любых $x\in X_0\subset X$ и $L\subset \KK^{X_0}$}.
\end{equation}
Если правая часть \eqref{eq3_1} $=-\infty$, то в
\eqref{eq3_1} имеет место равенство.
\end{propos}

В свете предложения \ref{pr3_1} и в связи с задачами \ref{pr1}--\ref{pr_3} будем исследовать условия,  
при которых  правая часть в неравенстве из \eqref{eq3_1} принадлежит $\KK$ и в \eqref{eq3_1} возможно равенство, или противоположное неравенство, в случае проективного предела векторных решёток $X$.

\begin{definition}\label{def3_1}
Пусть $X=\projlim X_np_n$ ---  приведённый правильный проективный предел векторных решёток $X_n$, $f\in \KK_{\pm\infty}^X$. Функцию
\begin{equation}\label{eq3_3}
\sup\tmin\pr_n f \colon x_n \mapsto \sup \bigl\{ f(x):  x\in X,  \; {\pr}_nx=x_n\bigr\},\quad
x_n\in X_n,
\end{equation}
называем {\it $\sup$-проекцией $f$ на $X_n$}.
Функции $\sup\tmin \pr_n f$ могут принимать значение
$\pm\infty$.
\end{definition}

\begin{propos}\label{pr_fn} В обозначениях определения\/ {\rm \ref{def3_1}} при $f_n:=\sup\tmin \pr_n f\colon X_n\to \KK_{\pm\infty}$
\begin{enumerate}[{\rm (i)}]
\item\label{fi}  $f_n(\pr_nx)\geqslant f_{n+1}(\pr_{n+1}x)\geqslant f(x)$ 
для любых $x\in X$ и $n\in {\NN}$;
\item\label{fii} для любых  $X_0\subset X$, $L\subset \KK^{X_0}$, $n\in \NN$, $L_n\subset \KK^{\pr_n X_0}$ при условии включения
\begin{equation}\label{LnL}
L_n\circ \pr_n:=\{ l_n\circ \pr_n\colon l_n\in L_n\}\subset L
\end{equation}
имеет место неравенство
\begin{multline}\label{eq3_4}
{\uE}^L_f(x)\overset{\eqref{luEu}}{=}\inf\bigl\{ l(x)\colon l \in L, \; f\bigm|_{X_0} \leqslant l\bigm|_{X_0} \bigr\} 
\\
\leqslant
\inf\bigl\{ l_n (\pr_n x)\colon 
l_n \in L_n, \; f_n\bigm|_{\pr_nX_0}\leqslant l_n\bigm|_{\pr_nX_0}\bigr\}
\overset{\eqref{luEu}}{=}{\uE}^{L_n}_{f_n}(\pr_n x) \quad\text{для всех $x\in X_0$}.
\end{multline}
\end{enumerate}
\end{propos}
\begin{proof} \eqref{fi}. Последнее неравенство сразу следует из определения \ref{def3_1},
а первое --- из \eqref{eq3_3} и соотношения \eqref{dfX} определения \ref{dfprp} элементов проективного предела.

\eqref{fii}. Пусть правая часть \eqref{LnL} $<+\infty$ и  
$f_n\bigm|_{\pr_nX_0} \leqslant l_n\bigm|_{\pr_nX_0}$, $l_n\in L_n$. Для  $l:=l_n\circ {\pr}_n\overset{\eqref{LnL}}{\in} L$ по определению \ref{def3_1}
при всех $x\in X_0$ имеем 
$f(x)\overset{\eqref{eq3_3}}{\leqslant} f_n(\pr_nx)\leqslant l_n(\pr_n x)=l(x)$. Отсюда $f\bigm|_{X_0}\leqslant l\bigm|_{X_0}$ и $l_n(\pr_nx)$ не меньше левой части \eqref{eq3_4}. Переходя к $\inf$ по всем $l_n\bigm|_{\pr_nX_0}\geqslant f_n\bigm|_{\pr_nX_0}$, получаем требуемое неравенство \eqref{eq3_4}.
\end{proof}

\begin{definition}\label{df3_2}
Пусть $(x^{(k)})_{k\in \NN}$ --- последовательность элементов упорядоченного множества $X$. Её {\it верхний предел\/}
определяется как
\begin{equation}\label{lsxk}
\limsup_{k\to \infty} x^{(k)} :=\inf_{n\in \NN}\sup_{k\geqslant n} x^{(k)},
\end{equation}
если все $\sup$ и $\inf$ в правой части существуют.
Последовательность $(x^{(k)})$ из  $X=\projlim_n X_np_n$  {\it стабилизируется к\/} $x\in X$,
если ${\pr}_nx^{(k)}={\pr}_nx$ при всех $k\geqslant n$.
\end{definition}

\begin{propos}\label{pr_stab}
Пусть $(x^{(k)})$ --- последовательность из  $X=\projlim_n X_np_n$. Тогда
\begin{enumerate}[{\rm (i)}]
\item\label{sti}
если ${\pr}_k x^{(k)}={\pr}_kx$ для  бесконечной последовательности номеров $k$, то последовательность  $(x^{(k)})$  стабилизируется к $x\in X$;
\item\label{stii} если $X$ --- приведённый правильный  проективный предел векторных решёток и\/ $(x^{(k)})$ стабилизируется к\/ $x\in X$, то она ограничена и существует\/ $\limsup_{k\to \infty} x^{(k)} \overset{\eqref{lsxk}}{=}x$.
\end{enumerate}
\end{propos}

\begin{proof} \eqref{sti}. Сразу следует из \eqref{dfX}--\eqref{1_2}.

\eqref{stii}. Пусть $x=(x_n)$. При фиксированном $n$ среди ${\pr}_nx^{(k)}$, $k\in {\NN}$,
лишь конечное число проекций отличается от $x_n$.
По предложению  \ref{prpr}\eqref{pri} 
и определению  \eqref{lsxk} получаем требуемое.
\end{proof}

\begin{theorem}\label{th3_1} 
Пусть $X=\projlim X_np_n$ --- приведённый правильный проективный предел векторных решёток $X_n$, ${f}\in \KK_{\pm\infty}^X$,  
$ {f}_n\overset{\eqref{eq3_3}}{:=}\sup\tmin \pr_n {f}$, $n\in \NN$, а также при некотором $n\in \NN$ имеем ${f}_n(X_n)\subset \KK_{-\infty}$,  $X_0\subset X$,
$L\subset \KK^{X_0}$, $L_n \subset \KK^{\pr_n X_0}$, $n\in \NN$,  существует $n_{f}(X_0)\in \NN$, для которого
\begin{equation}\label{drqn}
{f}_n(\pr_nx)\overset{\eqref{eq3_4}}{=} {\uE}^{L_n}_{{f}_n}(\pr_nx) \quad\text{для всех $x\in X_0$ и $n\geq n_{f}(X_0)$}.
\end{equation}
Для некоторого $x\in X_0$ предполагаем, что выполнено одно из следующих двух условий: 

\begin{enumerate}[{\rm (i)}]
\item\label{stabi}   для любой стабилизирующейся к $x$
последовательности $(x^{(k)})$ из $X$ 
выполнено
\begin{equation}\label{eq3_6}
\limsup_{k\to \infty} {f}\bigl(x^{(k)}\bigr)\leqslant
{f}\Bigl(\limsup_{k\to \infty} x^{(k)}\Bigr),
\end{equation}
\item\label{stabii} 
$f$ --- возрастающая функция на $X$
и для любой убывающей и стабилизирующейся
к $x\in X_0$ последовательности  $(x^{(k)})$ из  $X$
имеем
\begin{equation}\label{eq3_7}
\inf_k {f}\bigl(x^{(k)}\bigr)\leqslant {f}\Bigl(\inf_k x^{(k)}\Bigr).
\end{equation}
\end{enumerate}
Если  для $L$ и $L_n$ выполнено \eqref{LnL} при  всех $n\geq  n_f(X_0)$, то  
\begin{equation}\label{drqnq}
{f}(x)\overset{\eqref{eq3_4}}{=} {\uE}^{L}_{{f}}(x)\quad\text{для этого $x\in X_0$}.
\end{equation}
\end{theorem}
\begin{proof} Пусть для некоторого $c\in \KK$ выполнено неравенство 
\begin{equation}\label{eq3_2}
c<\inf\bigl\{ l(x)\colon l\in L=\lin \KK^{X_0}, \; {f}\bigm|_{X_0}
\leqslant l\bigm|_{X_0} \bigr \}\overset{\eqref{luEu}}{=}{\uE}^{L}_{f}(x) .
\end{equation}
Ввиду \eqref{LnL} по предложению \ref{pr_fn}\eqref{fii} согласно  \eqref{eq3_4} и \eqref{drqn} при $n\geq n_{f}(X_0)\geq n_{f}$ имеем
\begin{equation*}
c\overset{\eqref{eq3_4}}{<}
\inf\bigl\{ l_n(x_n)\colon l_n\in L_n,\;
 {f}_n\bigm|_{\pr_n X_0}\leqslant l_n\bigm|_{\pr_n X_0}\bigr\}
\overset{\eqref{luEu}}{=} {\uE}^{L_n}_{{f}_n}(x_n)\overset{\eqref{drqn}}{=}{f}_n(x_n).
\end{equation*}
Отсюда по определению \ref{def3_1} для $f_n$
при любом $n\geq n_{f}(X_0)$ найдётся $x^{(n)}\in X$, для которого
\begin{equation}\label{eq3_9}
c<{f}\bigl(x^{(n)}\bigr),\quad {\pr}_n x^{(n)}=x_n.
\end{equation}
Ввиду последнего равенства последовательность $(x^{(n)})$ по предложению \ref{pr_stab} стабилизируется к $x$  и $\limsup\limits_{k\to \infty} x^{(k)}=x$.
Из \eqref{eq3_9} по условию \eqref{eq3_6} из п.~\eqref{stabi}   получаем
\begin{equation*}
c\overset{\eqref{eq3_9}}{\leqslant} \limsup_{k\to \infty} {f}\bigl(x^{(k)}\bigr)\overset{\eqref{eq3_6}}{\leqslant}
{f}\Bigl(\limsup_{k\to \infty} x^{(k)}\Bigr)={f}(x).  
\end{equation*}
Отсюда и из \eqref{eq3_2} следует неравенство, противоположное 
\eqref{eq3_1}, т.\,е. имеем \eqref{drqnq}.

Покажем, что из \eqref{stabii} следует \eqref{stabi}, т.\,е.
\eqref{eq3_7} влечёт за собой \eqref{eq3_6}. Последовательность $(x^{(k)})$, стабилизирующаяся к $x$, по предложению \ref{pr_stab}\eqref{stii} ограничена. Значит существуют
${\hat x}^{(k)}=\sup_{n\geqslant k} x^{(n)}$.  
По построению $({\hat x}^{(k)})$ ---  убывающая последовательность,
которая по построению и предложению \eqref{prpr} также стабилизируется к $x$ и $\inf {\hat x}^{(k)}=x$. Если выполнено условие \eqref{eq3_7}, то в силу возрастания  ${f}$ по определению верхнего предела \eqref{lsxk} имеем
\begin{equation*}
{f}\Bigl(\limsup_{n\to \infty} x^{(n)}\Bigr)= 
{f}\Bigl(\, \inf_k {\hat x}^{(k)}\Bigr) \overset{\eqref{eq3_7}}{\geqslant} \inf_k {f}\bigl({\hat x}^{(k)}\bigr)=\inf_k  {f}\Bigl(\,\sup_{n\geqslant k} x^{(n)}\Bigr) \geqslant \inf_k \sup_{n\geqslant k} {f}\bigl(x^{(n)}\bigr) =\limsup_{n\to \infty} {f}\bigl(x^{(n)}\bigr).
\end{equation*}
Таким образом,  выполнено \eqref{eq3_6} и доказательство завершено. 
\end{proof}
\begin{remark}\label{uenseq} Согласно замечанию \ref{prl_subs} в теореме \ref{th3_1} достаточно требовать выполнения условий \eqref{drqn} и \eqref{LnL} не при всех  $n\geq  n_f(X_0)$, а только для какой-либо  бесконечной последовательности  номеров $n\geq  n_f(X_0)$.
\end{remark}

\subsection{Случай суперлинейной функции}\label{oespr}
$C$ --- {\it конус в векторном пространстве\/} $X$, если для любого $t\in \RR$ имеем $tC:=\{tc\colon c\in C\}\subset C$ при любом $t\in \RR_*^+$. Конус $C$ {\it выпуклый}, если для любых $x_1,x_2\in X$ имеем $x_1+x_2\in C$. 
Конус $C$ {\it с вершиной\/ $0$,\/} если $0\in C$.

\begin{definition}\label{qspl}
Функция $f\colon C\rightarrow \KK_{-\infty}$ на выпуклом конусе $C$ 
из векторного пространства $X$ {\it суперлинейная\/}, если она {\it супераддитивная}:
$f(x_1+x_2)\geqslant f(x_1)+f(x_2)$, {\it для всех\/}
$x_1, x_2\in C$,
и {\it строго положительно однородная}:
$f(tx)=t f(x)$, {\it для всех\/ 
$x\in C$ и\/  $t \in \RR_*^+$.} Если, в дополнение, $C$ c вершиной в $0$, то требуем  $f(0)\in \KK$, откуда $f(0)\overset{\eqref{v0}}{=}0$. 
Множество всех таких суперлинейных функций обозначаем как 
$\spl \KK_{-\infty}^C$, $\KK_{-\infty}^{C}:=(\KK_{-\infty})^{C}$,
$\spl^{\uparrow} \KK_{-\infty}^C:=\spl \KK_{-\infty}^C \cup \{\boldsymbol{+\infty}\}$, где 
функция $\boldsymbol{+\infty}$
--- тождественная $+\infty \in \KK_{+\infty}$ из \eqref{Kinfty}.  
Каждая функция $f\in \spl \KK_{-\infty}^C$ может быть канонически продолжена на всё $X$ как функция из $f\in \spl \KK_{-\infty}^X$ значениями $-\infty$ вне $C$, но функция $\boldsymbol{+\infty}\in \spl^{\uparrow} \KK_{-\infty}^C$ канонически продолжается  на $X$ как функция $\boldsymbol{+\infty}\in \spl^{\uparrow} \KK_{-\infty}^X$.
Аналогично определяются множества всех сублинейных функций
$\sbl \KK_{+\infty}^C:=-\spl \KK_{-\infty}^C$, 
$\sbl_{\downarrow} \KK_{+\infty}^C:=-\spl^{\uparrow} \KK_{-\infty}^C\ni \boldsymbol{-\infty}$ и их канонические продолжения на $X$: значением $+\infty$ на $X\setminus C$
для функций из $\sbl \KK_{+\infty}^X$, но 
$\boldsymbol{-\infty}\in \sbl_{\downarrow} \KK_{+\infty}^X$
для $\boldsymbol{-\infty}\in \sbl_{\downarrow} \KK_{+\infty}^C$.
\end{definition}
В связи с приложениями к теории функций здесь и далее, в отличие от подраздела \ref{ist} из введения, нам удобнее иметь дело не с сублинейными, а с <<противоположными>> им суперлинейными функциями и их верхними огибающими по пространству линейных функций.

\begin{propos}\label{spl_n}
В обозначениях определений\/ {\rm \ref{def3_1}}--{\rm  \ref{qspl}} пусть $f\in \spl\KK_{-\infty}^X$, а также $f_n\overset{\eqref{eq3_3}}{:=}\sup\tmin\pr_n f$. Если $f_n(X_n)\subset \KK_{-\infty}$, то  $f_n\in \spl \KK_{-\infty}^{X_n}$.
\end{propos}
\begin{proof} Пусть $x_n,x_n'\in X_n$.
Супераддитивность следует из цепочки (не)равенств
\begin{multline*}
f_n(x_n+x_n')=\sup\bigl\{ f(x'') \colon {\pr}_nx''=x_n+x_n'\bigr\}
 \geqslant \sup\bigl\{ f(x+x')\colon  {\pr}_nx=x_n,\; {\pr}_nx'=x_n'\bigr\} \\
 \geqslant \sup\bigl\{ f(x)+f(x')\colon  {\pr}_nx=x_n,\;  {\pr}_nx'=x_n' \bigr\}  \\
 =\sup\bigl\{ f(x)\colon {\pr}_nx=x_n\bigr\}  +
 \sup \bigl\{ f(x) \colon {\pr}_nx'=x_n' \bigr\}= f_n(x_n)+f_n(x_n').
\end{multline*}

Из определения \ref{def3_1} в силу положительной однородности проекции  ${\pr}_n$ легко следует строго положительная однородность $f_n$. Кроме того, по определению \ref{def3_1} и условию $f_n(X_n)\subset \KK_{-\infty}$ имеем $0=f(0)\leqslant f_n(0)< +\infty$, т.\,е. $f_n(0)\in \KK$. Таким образом, $f_n\in \spl \KK_{-\infty}^{X_n}$.
\end{proof}
\begin{propos}\label{prLlin} Пусть $X:=\projlim X_np_n$ приведённый правильный проективный предел векторных решёток, $X_0$ --- векторное подпространство в $X$. Тогда $\pr_n X_0$ --- векторное подпространство в $X_n$ и для   $L:= \lin \KK^{X_0}$ и 
$L_n:=\lin \KK^{\pr_n X_0}$ при  всех $n\in \NN$ выполнено \eqref{LnL}.
\end{propos}
\begin{proof} Из линейности проекции $\pr_n$ образ $\pr_n X_0$ векторного подпространства $X_0$ --- векторное подпространство в $X_n$, а для любой  функции $l_n\in \lin \KK^{\pr_n X_0}$ суперпозиция двух линейных функций  $l:=l_n\circ \pr_n\in \lin \KK^{X_0}$, что означает  $L_n\circ \pr_n\overset{\eqref{LnL}}{\subset} L$. 
\end{proof}
Из теоремы \ref{th3_1} и предложения \ref{prLlin}  получаем
\begin{corollary}\label{corspl} Пусть выполнены условия теоремы\/ 
{\rm \ref{th3_1}} до п.~\eqref{stabii}
включительно, а также дополнительно $f\in  \spl^{\uparrow}\KK_{-\infty}^X$ и $L,L_n$ выбраны как в предложении\/ {\rm \ref{prLlin}}. Тогда
имеем \eqref{drqnq}.
\end{corollary}
\begin{proof}  При $f\in \spl \KK_{-\infty}^X$ достаточное для \eqref{drqnq} условие  \eqref{LnL} --- в предложении \ref{prLlin}. При $f=\boldsymbol{+\infty}$ равенство \eqref{drqnq} ввиду соглашений \eqref{Kinfty}  очевидно.
\end{proof}
\begin{remark}\label{uenseqsq} Замечание \ref{uenseq} остаётся в силе и для  следствия \ref{corspl}.
\end{remark}

\subsection{Случай вогнутой функции}\label{oecon}
\begin{definition}\label{qcon}
Функция $f\colon K\rightarrow \KK_{-\infty}$ {\it вогнутая\/}
на выпуклом  подмножестве $K$ векторного пространства  $X$, 
если <<противоположная>>  функция $-f \colon K\to \KK_{+\infty}$, $(-f)(x)= -f(x)$ при $x\in K$, --- {\it выпуклая} в смысле \eqref{cof} c заменой $X$ на $K$. 
Множество всех таких вогнутых функция обозначаем как $\conc \KK_{-\infty}^{K}$, $\conc^{\uparrow} \KK_{-\infty}^{K}
\overset{\eqref{Kinfty}}{:=}(\conc \KK_{-\infty}^{K})\cup \{\boldsymbol{+\infty}\}$.
Каждая функция $f\in \conc \KK_{-\infty}^K$ может быть канонически продолжена на всё $X$ как функция из $f\in \conc \KK_{-\infty}^X$ значениями $-\infty$ вне $K$, но функция $\boldsymbol{+\infty}\in \conc^{\uparrow} \KK_{-\infty}^K$ канонически продолжается  на $X$ как функция $\boldsymbol{+\infty}\in \conc^{\uparrow} \KK_{-\infty}^X$. Аналогично $\conv \KK_{+\infty}^{K}$ ---
множество всех выпуклых функций на $K$ c каноническими выпуклым продолжения значением $+\infty$ на $X\setminus K$, но функция 
$\boldsymbol{-\infty}\overset{\eqref{Kinfty}}{\in} \conv_{\downarrow} \KK_{+\infty}^K:=
(\conv \KK_{+\infty}^K)\cup \{\boldsymbol{-\infty}\}$ канонически продолжается  на $X$ как $\boldsymbol{-\infty}\in \conv_{\downarrow} \KK_{+\infty}^X$. В частности, 
$\aff \KK^X=\conc \KK^{X}\cap \conv \KK^{X}$.
\end{definition}

Также в связи с приложениями к теории функций здесь, в отличие от подраздела \ref{ist} из введения, нам удобнее иметь дело не с выпуклыми, а с вогнутыми функциями и их верхними огибающими по пространству  аффинных функций.

\begin{propos}\label{conc_n}
В обозначениях определений\/ {\rm \ref{def3_1}, \ref{pr_fn}} и\/ {\rm  \ref{qspl}} пусть $f\in \conc \KK_{-\infty}^X$, а также $f_n\overset{\eqref{eq3_3}}{:=}\sup\tmin\pr_n f$. Если $f_n(X_n)\subset \KK_{-\infty}$, то  $f_n\in \conc \KK_{-\infty}^{X_n}$.
\end{propos}
\begin{proof} Пусть $x_n,x_n'\in X_n$, $t\in (0,1)$.
Выпуклость следует из цепочки (не)равенств 
\begin{multline*}
f_n\bigl(tx_n+(1-t)x_n'\bigr)=\sup\bigl\{ f(x'') \colon {\pr}_nx''=tx_n+(1-t)x_n'\bigr\}
 \\
\geqslant \sup\bigl\{ f\bigl(x_t+x_t'\bigr)\colon  {\pr}_n x_t=t x_n,\; {\pr}_n x_t'= (1-t)x_n'\bigr\} \\
=\sup\bigl\{ f\bigl(tx+(1-t)x'\bigr)\colon  {\pr}_n x= x_n,\; {\pr}_n x'= x_n'\bigr\}\\
 \geqslant \sup\bigl\{ tf(x)+(1-t)f(x')\colon  {\pr}_nx=x_n,\;  {\pr}_nx'=x_n' \bigr\}  \\
 =t\sup\bigl\{ f(x)\colon {\pr}_nx=x_n\bigr\}  +
(1-t) \sup \bigl\{ f(x)\colon {\pr}_nx'=x_n' \bigr\}= tf_n(x_n)+(1-t)f_n(x_n'),
\end{multline*}
использующей линейность проекции  $\pr_n\colon X\to X_n$. 
\end{proof}

\begin{propos}\label{prLaff} Пусть выполнены условия предложения\/ {\rm \ref{prLlin}}.
Тогда  для   $L:= \aff \KK^{X_0}$ и 
$L_n:=\aff \KK^{\pr_n X_0}$ при  всех $n\in \NN$ выполнено \eqref{LnL}.
\end{propos}
\begin{proof} В предложении \ref{prLlin} уже отмечалось, что  $\pr_nX_0$ --- векторное подпространство в $\pr_nX_0$, значит корректно определено $L_n$. 
Для любой  функции $a_n\in \aff \KK^{\pr_n X_0}$ суперпозиция её с линейной функцией $\pr_n$ даёт $a:=a_n\circ \pr_n\in \aff \KK^{X_0}$, что означает  $L_n\circ \pr_n\overset{\eqref{LnL}}{\subset} L$. 
\end{proof}

Из теоремы \ref{th3_1} и предложения \ref{prLaff}  получаем
\begin{corollary}\label{coraff} Пусть выполнены условия теоремы\/ {\rm \ref{th3_1}} до п.~\eqref{stabii}
включительно, а также дополнительно $f\in  \conc^{\uparrow} \KK_{-\infty}^X$ и $L,L_n$ выбраны как в предложении\/ {\rm \ref{prLaff}}. Тогда имеем \eqref{drqnq}.
\end{corollary}
\begin{proof} 
При $f\in \conc \KK_{-\infty}^X$ достаточное для \eqref{drqnq} условие  \eqref{LnL} --- в предложении \ref{prLaff}. При $f=\boldsymbol{+\infty}$ равенство \eqref{drqnq} ввиду соглашений \eqref{Kinfty}  очевидно.
\end{proof}
\begin{remark}\label{uenseqaf} Замечание \ref{uenseq} остаётся в силе и для  следствия \ref{coraff}.
\end{remark}

\section{Супремальные функции 
и выметание}\label{suprf}
\setcounter{equation}{0}

\subsection{Супремальная функция}\label{sprf}

\begin{definition}\label{defsprf}
{\it Cупремальной функцией, порождённой функцией\/
$q\in \KK^H$ относительно\/ $H\subset X$\/} на 
упорядоченном множестве $X$ с отношением порядка $\leq$, называем  функцию 
\begin{equation}\label{eq4_1}
\spf_{H,q}\colon 
x\mapsto \sup \bigl\{ q(h) \colon H\ni h\leqslant x\bigr\}
\overset{\eqref{Kinfty}}{\in} \KK_{\pm \infty},\quad x\in X. 
\end{equation}
\end{definition}

\begin{propos}\label{pr_spf} Cупремальная функция \eqref{eq4_1} обладает следующими свойствами.
\begin{enumerate}[{\rm 1.}]
\item\label{fHpI} $\spf_{H,q}(x)\in \KK_{+\infty}$, если и только если существует такой элемент $h\in H$, что $h\leq x$.
\item\label{fHpII}  $\spf_{H,q}$ возрастающая на $X$, по $H$ относительно отношения включения и по  $q$ в том смысле, что 
при $x\leq x'\in X$, $H\subset H'\subset X$ и для $q'\in \KK^{H'}$ 
при  $q\bigm|_{H}\leq q'\bigm|_H$ имеем $\spf_{H,q}(x)\leq \spf_{H',q'}(x')$.
\item\label{fHpIII} Пусть $X$ --- упорядоченное векторное пространство. Тогда
\begin{enumerate}[{\rm (i)}]
\item\label{Xosi} если $0\in H$ и $q(0)\geq 0$, то $\spf_{H,q}$ --- положительная функция;
\item\label{Xosii} если $H$ --- конус в $X$ и $q$ --- строго положительно однородная функция на $H$, то $\spf_{H,q}$ --- строго положительно однородная функция на $X$; если, в дополнение, $H$ содержит отрицательный вектор, т.\,е.
$H\cap (-X^+)\neq \varnothing$, то $\spf_{H,q}(0)=+\infty$ или $\spf_{H,q}(0)=0$, а также $\spf_{H,q}$ --- положительная функция {\rm (пишем $\spf_{H,q}\in (\KK_{-\infty}^X)^+$);}  
\item\label{Xosiii} если $H$ --- выпуклое множество, $\spf_{H,q}(X)\subset \KK_{-\infty}$ и 
$q\in \conc \KK^H$, то $\spf_{H,q}\in \conc \KK_{-\infty}^X$; 

\item\label{Xosiv} если $H$ --- выпуклый конус, $\spf_{H,q}(X)\subset \KK_{-\infty}$ и  $q$ --- суперлинейная функция из  $\spl \KK^H$, то
$\spf_{H,q}$ --- супераддитивная строго положительно однородная функция на $X$; если, в дополнение, $H\cap (-X^+)\neq \varnothing$, то $\spf_{H,q}\in \spl^+ \KK_{-\infty}^X:=
 (\KK_{-\infty}^X)^+ \cap \spl \KK_{-\infty}^X$. 
\end{enumerate}
\end{enumerate}
\end{propos}
\begin{proof} Пункты \ref{fHpI}, \ref{fHpII} и 
\ref{Xosi} очевидным образом следуют из определения \ref{defsprf}, \eqref{eq4_1}.

\ref{Xosii}. Пусть $t\in \RR_*^+$. Первая часть о строгой положительной однородности следует из равенств
\begin{equation*}
\spf_{H,q}(tx)=\sup \bigl\{ q(h) \colon H\ni h\leqslant tx\bigr\}=\sup \left\{ tq\Bigl(\frac{1}{t}h\Bigr) \colon H\ni \frac{1}{t} h\leqslant x\right\}=t\spf_{H,q}(x).
\end{equation*}
Если при этом $H$ содержит отрицательный вектор, то $-\infty\overset{\eqref{eq4_1}}{<}\spf_{H,q}(0)$, откуда либо 
$\spf_{H,q}(0)=+\infty$, либо, в противном случае, ввиду 
$\spf_{H,q}(0)=2\spf_{H,q}(0)$, имеем $\spf_{H,q}(0)=0$. B любом из случаев по п.~\ref{fHpII} возрастающая на $X$ функция  $\spf_{H,q}$ положительна.

\ref{Xosiii}. Пусть $t_1,t_2\in \RR_*^+$ и $t_1+t_2=1$. Тогда 
для любых $x_1,x_2\in X$ имеет место цепочка (не)равенств
\begin{multline}\label{cpfconc}
\spf_{H,q}(t_1x_1+t_2x_2)\overset{\eqref{eq4_1}}{=}
\sup \bigl\{ q(h) \colon H\ni h\leqslant t_1x_1+t_2x_2\bigr\}
\\
\geq \sup \bigl\{ q(t_1h_1+t_2h_2) \colon H\ni  h_1\leq x_1,\;  H\ni h_2\leq x_2\bigr\}
\\
\geq \sup \bigl\{ t_1q(h_1)+t_2q(h_2) \colon H\ni h_1\leq x_1,\;  H\ni h_2\leq x_2\bigr\}
\\
=\sup \bigl\{ t_1q(h_1)\colon H\ni h_1\leq x_1\bigr\}
+\sup \bigl\{t_2q(h_2) \colon  H\ni h_2\leq x_2\bigr\}
\overset{\eqref{eq4_1}}{=}t_1 \spf_{H,q}(x_1) + t_2 \spf_{H,q}(x_2).
\end{multline}
\ref{Xosiv}. Ввиду $+\infty \notin \spf_{H,q}(X)$ и  $q\in \spl \KK^H$ по-прежнему имеет место  цепочка (не)равенств 
\eqref{cpfconc}, но уже для всех $t_1,t_2\in \RR_*^+$ и $x_1,x_2\in X$, что доказывает супераддитивность и строго положительную однородность $\spf_{H,q}$ на $X$.  
При $H\cap (-X^+)\neq \varnothing$ и $+\infty \notin \spf_{H,q}(X)$ согласно  п.~\ref{Xosii} получаем  $\spf_{H,q}(0)=0$ и $\spf_{H,q}\in \spl^+\KK_{-\infty}^X$.
\end{proof}

\begin{propos}[{\rm ср. с \cite[предложение 4.2]{Kh01_1}}]\label{prqn} 
Пусть $X=\projlim_n X_np_n$
--- приведённый правильный проективный предел векторных решёток, $n\in \NN$,
$H\subset X$ и $H_n:=\pr_n H$, $n\in \NN$, 
\begin{equation}\label{q1qn}
q_1\in \KK^{H_1},\quad  q:=q_1 \circ \pr_1 \colon H\to \KK, \quad q_n:=q_1\circ p_1\circ \dots \circ p_{n-1} 
\colon H_n\to \KK.  
\end{equation} 
Тогда  в обозначении \eqref{eq4_1} определения\/ {\rm \ref{defsprf}} супремальная функция $\spf_{H,q}$, 
порождённая функцией $q$ относительно $H$,
через  её $\sup$-проекцию\/ $\sup\tmin\pr_n$ из определения\/ {\rm \ref{def3_1}} в обозначении \eqref{eq3_3} на $X_n$ связана с супремальными функциями $\spf_{H_n,q_n}$ равенствами 
\begin{equation}\label{spflln}
\sup\tmin\pr_n \spf_{H,q}=\spf_{H_n,q_n}\quad\text{на $X_n$ при любых $n\in \NN$}.
\end{equation}
\end{propos}
\begin{proof} При  $x_n\in X_n$ по определениям  \ref{def3_1}  и \ref{defsprf}  
для левой части \eqref{spflln} имеем 
\begin{equation}
\begin{split}\label{eq4_6}
\bigl(\sup\tmin\pr_n \spf_{H,q}\bigr) (x_n)
&\overset{\eqref{eq3_3}}{=}
\sup \bigl\{\spf_{H,q}(x)\colon x\in X,\; \pr_n x=x_n \bigr\}\\ 
&\overset{\eqref{eq4_1}}{=}
\sup \Bigl\{\sup\bigl\{q(h)\colon H\ni h\leq x\bigr\}
\colon x\in X,\; \pr_n x=x_n \Bigr\}\\
&=\sup \bigl\{q(h)\colon H\ni h\leq x,\;
\pr_n x=x_n \bigr\}\\
&\overset{\eqref{q1qn}}{=} \sup \bigl\{q_n(\pr_n h)\colon H\ni h\leq x,\; \pr_n x=x_n \bigr\}\\
&=\sup \bigl\{q_n(h_n)\colon h\in H, \; h\leq x,\;
h_n=\pr_n h, \;  \pr_n x=x_n \bigr\}.
\end{split}
\end{equation} 
В то же время для правой части \eqref{spflln} по определению \ref{defsprf} имеем
\begin{equation}\label{spfln}
\spf_{H_n,q_n}(x_n)\overset{\eqref{eq4_1}}{=}
\sup\bigl\{q_n(h_n)\colon h_n\in H_n,\; h_n\leq_n x_n \bigr\},
\end{equation}
где $\leq_n$ --- отношение порядка на $X_n$.
Если покажем, что множество
\begin{equation}\label{eq4_7}
\bigl\{h_n\colon h\in H, \; h\leq x,\;
h_n=\pr_n h, \;  \pr_n x=x_n \bigr\},
\end{equation}
по которому берётся $\sup$ в правой части \eqref{eq4_6},
совпадает с множеством
\begin{equation}\label{eq4_8}
\bigl\{h_n\colon h_n\in H_n,\; h_n\leq_n x_n \bigr\},
\end{equation}
по которому берётся $\sup$ в правой части \eqref{spfln},
то \eqref{eq4_6} совпадает с \eqref{spfln} при всех $x_n\in X_n$ и равенство \eqref{spflln} будет доказано.

Если $h_n$ из \eqref{eq4_7}, то $h_n{\leq}_n x_n$, следовательно, множество \eqref{eq4_7} включается в \eqref{eq4_8}.

Обратно, пусть $h_n$ из \eqref{eq4_8}. Существуют
элементы $h\in H$ и $x'\in X$, для которых 
$h_n={\pr}_n h$, ${\pr}_nx'=x_n$. Полагаем $x=\sup \{ h,x'\}$.
Тогда ${\pr}_nx=x_n\,$ и, очевидно, $h\leqslant x$. Следовательно, \eqref{eq4_8} содержится в \eqref{eq4_7}, что и требовалось.
\end{proof}

\subsection{Выметание и росток функции}\label{Sbal} 
\begin{definition}[{\rm обобщение \cite[определение 1]{Kh06}, \cite[III.1.3]{AK}, \cite[гл.~XI, \S~3, О41]{M}}]\label{balpb}
Пусть $X$ --- упорядоченное множество с отношением порядка $\leq$,  $H\subset X$, $q\in \KK^H$, $X_0\subset X$. {\it Функцию\/} $b\colon  X_0 \to \KK_{\pm\infty}$
называем {\it выметанием  функции\/ $q$ относительно\/ пары $(H,X_0)$} и пишем $q {\prec_H^{X_0}} b$,  если 
для любых $h\in H$ и $x\in X_0$ при $h\leq x$ имеем $q(h)\leq b(x)$.
Пусть $L\subset \KK_{\pm\infty}^{X_0}$.    <<Луч>> 
\begin{equation}\label{sprqb}
\Spr(q; H,X_0,L):=\bigl\{b\in L
\colon q {\prec_H^{X_0}} b\bigr\},
\end{equation}
следуя \cite[III.1.3]{AK} и \cite[\S~4]{Kh01_1}, будем называть {\it ростком функции\/ $q$ относительно тройки\/ $(H,X_0,L)$}. 
\end{definition}
Элементарную взаимосвязь между супремальной функцией и выметанием даёт 
\begin{propos}[{\rm ср. с \cite[предложение 1]{Kh06}}]\label{prsprn}
В обозначениях определений\/ {\rm \ref{defsprf}}  и\/ 
{\rm \ref{balpb}} 
\begin{multline}\label{s_b}
\spf_{H,q}(x) \overset{\eqref{eq4_1}}{=}
 \sup \bigl\{ q(h) \colon H\ni h\leqslant x\bigr\}
\leq \inf \bigl\{b(x)\colon q {\prec_H^{X_0}} b\bigr\}\\
\overset{\eqref{sprqb}}{\leq}
\inf \bigl\{b(x)\colon  b\in \Spr(q; H,X_0,L) \bigr\}
=\inf \bigl(\Spr(q; H,X_0,L) (x)\bigr)
\quad\text{для всех $x\in X_0$}.
\end{multline}
Так, для $x\in X_0$ существование элемента $h\in H$, удовлетворяющего неравенству $h\leq x$, влечёт за собой строгое  неравенство $-\infty < \inf \bigl\{b(x)\colon q {\prec_H^{X_0}} b\in L\bigr\}$, т.\,е. $\inf \bigl(\Spr(q; H,X_0,L) (x)\bigr)\neq -\infty$.
\end{propos}
\begin{proof} 
Если в \eqref{s_b} левая часть $=-\infty$ или правая часть равна $=+\infty$, то неравенство \eqref{s_b} очевидно. Если  левая часть $\neq -\infty$, то существуют элемент $h\in H$, удовлетворяющий неравенству $h\leq x$. Для любых  $H\ni h\leq x\in x_0$ при $q\prec_H^{X_0} b$ имеем $q(h)\leq b(x)$, откуда 
$\sup_{h\leq x} q(h)\leq b(x)$. Применяя к правой части операцию $\inf$ по всем выметаниям $b$ функции $q$ относительно пары $(H,X_0)$, получаем первое неравенство в  \eqref{s_b}. Следующее неравенство в \eqref{s_b} и заключительное утверждение --- очевидное следствие из \eqref{s_b} ввиду $L\subset \KK_{\pm\infty}^{X_0}$.
\end{proof}
\begin{remark}\label{rmspr} В свете предложения \ref{prsprn} ключевой интерес представляют условия на $(X,\leq)$, $H$, $q$, $X_0$, $L$, при которых в \eqref{s_b} имеет место равенство левой и правой частей --- задачи \ref{pr1} и \ref{pr2} из подраздела \ref{stpr} введения. Для менее требовательной задачи \ref{pr_3} --- условия на те же объекты и на $x$, при которых из $\inf \bigl(\Spr(q; H,X_0,L) (x)\bigr)\neq -\infty$ следует $\spf_{H,q}(x)\neq -\infty$. 
\end{remark}

\begin{propos} В обозначениях определения\/ 
{\rm \ref{balpb}} имеют место следующие свойства.
\begin{enumerate}[{\rm 1.}]
\item\label{spr1}  При $L\subset L'\subset \KK_{\pm\infty}^{X_0}$ для всех $x\in X_0$ справедливы соотношения 
\begin{equation}\label{sprL}
\Spr(q; H,X_0,L)\subset \Spr(q; H,X_0,L'),\quad
\inf \bigl(\Spr(q; H,X_0,L)(x)\bigr)
\geq \inf \bigl(\Spr(q; H,X_0,L')(x)\bigr).
\end{equation}
\item\label{spr2}
Если $L\subset \incr \KK_{\pm\infty}^{X_0}$, то функция  
$\inf \bigl(\Spr(q; H,X_0,L)(\cdot)\bigr)$ возрастающая на $X_0$. 
\item\label{spr3} Если $H\subset H'\subset X$, то для всех 
$x\in X_0$ справедливы соотношения
\begin{equation}\label{sprH}
\Spr(q; H,X_0,L)\supset \Spr(q; H',X_0,L),\quad
\inf \bigl(\Spr(q; H,X_0,L)(x)\bigr)
\leq \inf \bigl(\Spr(q; H',X_0,L)(x)\bigr).
\end{equation}
\item\label{spr4} Если $q'\in \KK^H$ и $q\bigm|_H\leq q'\bigm|_{H}$, то 
для всех $x\in X_0$ справедливы соотношения
\begin{equation}\label{sprq}
\Spr(q; H,X_0,L)\supset \Spr(q'; H,X_0,L),\quad
\inf \bigl(\Spr(q; H,X_0,L)(x)\bigr)
\leq \inf \bigl(\Spr(q; H',X_0,L)(x)\bigr).
\end{equation}
\item\label{spr5}  
Если $X_0\subset X_0'\subset X$ и  для множества функций $L'\subset \KK_{\pm\infty}^{X_0'}$ множество их сужений $L'\bigm|_{X_0}$ на $X_0$ включено в $L$, то
$\inf \bigl(\Spr(q; H,X_0,L) (x)\bigr)\leq
\inf \bigl(\Spr(q; H,X_0',L') (x)\bigr)$ для $x\in X_0$.
\item\label{sprO} Пусть $X$ --- упорядоченное векторное пространство. Тогда 
\begin{enumerate}[{\rm (i)}]
\item\label{bXosi} если $0\in H$ и $q(0)\geq 0$, то $\inf \bigl(\Spr(q; H,X_0,L) (\cdot)\bigr)$ --- положительная функция на $X_0$;
\item\label{bXosii} если $X_0$ --- конус в $X$, 
а все функции
из $L\subset \KK_{\pm \infty}^{X_0}$ строго положительно однородны,  то $\inf \bigl(\Spr(q; H,X_0,L)\bigr)\subset \KK_{\pm \infty}^{X_0}$ --- строго положительно однородная функция на $X_0$; если, в дополнение, 
$H\cap (-X^+)\neq \varnothing$ и $0\in X_0$, то $\inf \bigl(\Spr(q; H,X_0,L) (0)\bigr)=+\infty$ или $=0$, 
а при $L\subset \incr \KK_{\pm\infty}^{X_0}$ имеем $\inf \bigl(\Spr(q; H,X_0,L) (\cdot)\bigr)\in (\KK_{-\infty}^{X_0})^+$;  
\item\label{bXosiii} если подмножество $X_0\subset X$ --- \underline{полугруппа по сложению} (соотв. \underline{выпуклое множество} или \underline{выпуклый конус c вершиной и $H\cap (-X^+)\neq \varnothing$}) и все функции из $L\subset \KK_{-\infty}$ \underline{супераддитивны} (соотв. \underline{вогнуты} или \underline{суперлинейны}) и $\Spr(q; H,X_0,L)(X_0)\subset \KK_{-\infty}^{X_0}$, то функция $\inf \bigl(\Spr(q; H,X_0,L) \bigr)$ \underline{супераддитивна} (соотв. \underline{вогнута} или \underline{суперлинейна}) на $X_0$; 
\end{enumerate}
\end{enumerate} 
\end{propos}
\begin{proof} Свойства из пп.~\ref{spr1}--\ref{spr5} и \ref{bXosi} легко  следуют из определения \ref{balpb}. 

\ref{bXosii}. Первая часть утверждения следует из равенств с $t\in \RR_*^+$: 
\begin{equation}
\inf \bigl\{b(tx)\colon  b\in \Spr(q; H,X_0,L) \bigr\}
=t\inf \bigl\{b(x)\colon  b\in \Spr(q; H,X_0,L) \bigr\}.
\end{equation}
Если  в $H\cap (-X^+)\neq \varnothing$, то $\spf_{H,q}(0)>-\infty$. По предложению \ref{prsprn} 
$\inf \bigl(\Spr(q; H,X_0,L) (0)\bigr)>-\infty$. При 
$\inf \bigl(\Spr(q; H,X_0,L) (0)\bigr)\neq +\infty$ имеем 
$\inf \bigl(\Spr(q; H,X_0,L) (0)\bigr)\in \KK$, откуда в силу строгой положительной однородности $\inf \bigl(\Spr(q; H,X_0,L) (0)\bigr)=0$. Это согласно п.~\ref{spr2}
при $L\subset \incr \KK_{\pm\infty}^{X_0}$ даёт положительность функции $\inf \bigl(\Spr(q; H,X_0,L) (\cdot)\bigr)$ на $X_0$. 

\ref{bXosiii}. Для краткости положим $S:=\Spr(q; H,X_0,L)$.
Для $t_1=t_2=1$ (соотв. для $t_1,t_2\in \RR_*^+, t_1+t_2=1$,
или для $t_1,t_2\in \RR_*^+$) имеем 
\begin{multline*}
\inf\bigl\{b(t_1x_1+t_2x_2)\colon b\in S\bigr\}\geq
\inf\bigl\{b(t_1x_1)+ b(t_2x_2)\colon b\in S\bigr\}\\
\geq \inf\bigl\{t_1b(x_1)+t_2 b(x_2)\colon b\in S\bigr\}
\geq t_1\inf\bigl\{b_1(x_1)\colon b_1\in S\bigr\}
+t_2\inf\bigl\{b_2(x_2)\colon b_2\in S\bigr\}.
\end{multline*}
Для доказательства суперлинейности осталось отметить, что по п.~\ref{bXosii} имеем  $(\inf S)(0)=0$.
\end{proof}

\begin{propos}[{\rm обобщение \cite[предложение 4.1]{Kh01_1}, \cite[III.1.3.VI]{AK}}]\label{pr_pf}
В обозначениях определений\/ {\rm \ref{defsprf}}  и\/ 
{\rm \ref{balpb}} имеет место равенство
\begin{equation}\label{spr_spf}
\Spr(q; H,X_0,L)=\bigl\{b\in L\colon \spf_{H,q}\bigm|_{X_0}\leq b\bigm|_{X_0} \bigr\} 
\end{equation}
\end{propos}
\begin{proof} Покажем, что правая часть \eqref{spr_spf} включена в левую. Пусть 
\begin{equation}\label{spf_b}
\spf_{H,q}\bigm|_{X_0}\leq b\bigm|_{X_0},\quad\text{где $b\in L$}.
\end{equation} 
Допустим, что $h\in H$ и $h\leq x\in X_0$. Тогда по определению \ref{defsprf}\eqref{eq4_1} супремальной функции, в силу её возрастания на $X$ 
по предложению \ref{pr_spf}, п.~\ref{fHpII}, имеем
\begin{equation}\label{qspfb}
q(h)=\spf_{H,q}(h)\leq \spf_{H,q}(x)\overset{\eqref{spf_b}}{\leq} b(x).
\end{equation}
Это по определению \ref{balpb} означает, что $b$ --- выметание
$q$ относительно пары $(H,X_0)$, т.\,е. $q\prec_H^{X_0}b$, и
$b\overset{\eqref{sprqb}}{\in} \Spr(q;H,X_0,L)$. Таким образом, правая часть \eqref{spr_spf} включена в левую.

Обратно, пусть $b\in \Spr(q; H,X_0,L)$. Тогда по определениям супремальной функции \eqref{eq4_1} и ростка функции $q$ относительно тройки $(H,X_0,L)$ в определении 
\ref{balpb} при $x\in X_0$  имеем     
\begin{equation*}
\spf_{H,q}(x)\overset{\eqref{eq4_1}}{:=}
\sup \bigl\{ q(h) \colon H\ni h\leqslant x\bigr\}\overset{\eqref{sprqb}}{\leq}
b(x), \quad\text{при всех $x\in X_0$.} 
\end{equation*}
Таким образом, левая часть \eqref{spr_spf} включена в правую.
\end{proof}

\section{Представление верхней огибающей на проективном пределе}\label{repue}
\setcounter{equation}{0}

\subsection{Порядково-линейная версия для выпуклого конуса $H$}\label{ord_ver} В этом подразделе $\KK$ --- пространство Канторовича. Нам потребуется (ср. с теоремами об огибающей из введения)
\begin{thHBK}[{\rm об огибающей \cite[1.4.14(2)]{KK}}]
Пусть $X$ --- векторное пространство и $f\in \KK^X$. Попарно эквивалентны следующие  три утверждения.     
\begin{enumerate}[{\bf 1.}]
\item\label{HBiK} $f$ --- суперлинейная функция на $X$, т.\,е. $f\in \spl \KK^X$.
\item\label{HBiiK} $f(x)\overset{\eqref{luEu}}{=}{\uE}^{\lin \KK^X}_f (x)$ для всех $x\in X$.
\item\label{HBiiiK}   Значения $f(x)$ тождественно равны 
\begin{equation}\label{repsblK}
\sup\left\{\sum_{k=1}^n t_kf(x_k)\colon
\sum_{k=1}^n t_k x_k=x, \; x_k\in X, \; t_k\in \RR^+, \;
n\in \NN\right\}\;
\text{ при всех $x\in X$}.
\end{equation}
\end{enumerate}
\end{thHBK}
\begin{definition}Пусть $(X,\leq)$ --- упорядоченное множество, $X_0\subset X$. Подмножество $H\subset X$ {\it минорирующее\/} для $X_0$, если для любого $x\in X_0$ существует такой элемент 
$h\in H$, что $h\leqslant x$.
\end{definition}

\begin{theorem}\label{thrdq}
 Пусть $X=\projlim X_np_n$ ---     приведённый правильный
 проективный предел векторных решёток  $X_n$,
$n\in \NN$, $q_1\in \incr \spl\KK^{X_1}$ и $q\overset{\eqref{q1qn}}{:=}q_1\circ \pr_1$,  $X_0\subset X$ --- векторное подпространство, $H\subset X$ --- выпуклый конус, 
$H_n:=\pr_n H\subset X_n$. Пусть выполнены следующие условия:
\begin{enumerate}[{\rm (i)}]
\item\label{qHii} каждая проекция $H_n$ содержит отрицательный элемент из $X_n$, т.\,е. $H_n\cap (-X_n^+)\neq \varnothing$;
\item\label{qHi} для каждого $n\in \NN$ проекция $H_n$  минорирующая для $\pr_n X_0$;
\item\label{qHiii-}  для любого  ограниченного сверху подмножества $H'\subset H$ имеем $\sup H'\in H$;
\item\label{qHiv}  если последовательность
\begin{equation}\label{seqhn}
(h^{(k)})_{k\in \NN} \subset H
\end{equation} 
убывающая ограниченная снизу в $X$, то
$\inf_{k\in \NN} h^{(k)} \in H$;
\item\label{qHv}  если последовательность
\eqref{seqhn}  убывающая  в $X$ и 
\begin{equation}\label{qh1}
q(h^{(1)})\geq \inf_{k\in \NN} q(h^{(k)})\in \KK,
\end{equation} 
то последовательность  \eqref{seqhn}
ограничена снизу в $X$, а также\footnote{Из неравенства в \eqref{qhk} следует равенство, поскольку противоположное неравенство для $q\in \incr \KK^X$ очевидно.} 
\begin{equation}\label{qhk}
q\bigl(\inf_{k\in \NN} h^{(k)}\bigr)\geqslant \inf_{k\in \NN} q(h^{(k)}), \quad \inf_{k\in \NN} h^{(k)}
\overset{\eqref{qHiv}}{\in} H.
\end{equation}
\end{enumerate}

Тогда  для $L:=\lin \KK^{X_0}$ супремальная функция  $\spf_{H,q}$ со значениями в $\KK_{-\infty}$ является супераддитивной и положительно однородной, а также допускает представление 
\begin{multline}\label{eq5_1}
\spf_{H,q}(x)\overset{\eqref{eq4_1}}{:=}
\sup \bigl\{ q(h) \colon H\ni h\leqslant x\bigr\}
\\=\inf \bigl\{l(x) \colon l\in L; \;  
q(h)\leq l(x') \text{ при всех $h\in H$, $h\leq x'$, $x'\in X_0$}\bigr\}
\\
=\inf \bigl\{l(x) \colon l\in L, \;  
q {\prec_H^{X_0}} l\bigr\}
\overset{\eqref{sprqb}}{=:}\inf
\bigl(\Spr(q; H,X_0,L) (x)\bigr)
\quad \text{для всех $x\in X_0$}.
\end{multline}
Если усилить\/ \eqref{qHii} до $H\cap (-X^+)\neq \varnothing$, то $\spf_{H,q}\in \spl^+\KK_{-\infty}^X$. 
\end{theorem}
\begin{proof} Будет использована теорема \ref{th3_1} в виде следствия \ref{corspl} из неё. 

 Положим $f:=\spf_{H,q}$ и 
$f_n=\sup\tmin\pr_n f$, $L:=\lin \KK^{X_0}$ и  $L_n:= \lin \KK^{\pr_n X_0}$. Тогда по предложению \ref{prLlin} выполнено \eqref{LnL} и   
по предложению \ref{prqn} в обозначениях \eqref{q1qn} 
имеет место равенства 
\begin{equation}\label{spflln+}
f_n=\sup\tmin\pr_n \spf_{H,q}\overset{\eqref{spflln}}{=}\spf_{H_n,q_n}\quad\text{на $X_n$ при любых $n\in \NN$},
\end{equation}
где по предложению \ref{pr_spf}, п.~\ref{Xosiv}, с учётом условия \eqref{qHii}, $\spf_{H_n,q_n}\in \spl^+ \KK_{-\infty}^{X_n}$, т.\,е. $f_n \overset{\eqref{spflln+}}{\in} \spl^+ \KK_{-\infty}^{X_n}$.
Кроме того, по предложению \ref{pr_spf}, п.~\ref{Xosiv}, функция $\spf_{H,q}$ супераддитивная и строго положительно однородная, а при $H\cap (-X^+)\neq \varnothing$
ещё и $\spf_{H,q}(0)=0$,  $\spf_{H,q}\in \spl^+\KK_{-\infty}^X$. 

Напомним, что  $\pr_nX_0$ --- векторное подпространство в $X_n$ по предложению \ref{prLlin}. Рассмотрим сужения $f_n\bigm|_{\pr_nX_0}=\spf_{H_n,q_n}\bigm|_{\pr_nX_0}$. В правой части здесь --- суперлинейная функция на $\pr_nX_0$ со значениями только в $\KK$ по определению
\ref{defsprf}, равенство \eqref{eq4_1}, супремальной функции
$\spf_{H_n,q_n}\bigm|_{\pr_nX_0}$  и по условию
\eqref{qHi} о минорировании. Таким образом, 
$f_n\bigm|_{\pr_nX_0} \in \spl^+ \KK^{\pr_n X_0}$. 
Тогда по теореме Хана\,--\,Банаха\,--\,Канторовича об огибающей, применённой в части импликации \ref{HBiK}$\Rightarrow$\ref{HBiiK} в векторном пространстве $\pr_n X_0$, имеем 
$f_n\bigm|_{\pr_nX_0}\overset{\eqref{luEu}}{=}{\uE}^{\lin \KK^{\pr_nX_0}}_{f_n|_{\pr_nX_0}}$
на $\pr_n X_0$.
Последнее равенство означает, что выполнено условие \eqref{drqn} теоремы \ref{th3_1}.  

Функция $f=\spf_{H,q}$ --- возрастающая на $X$ по п.~\ref{fHpII} предложения \ref{pr_spf}. Проверим выполнение условия \eqref{stabii} теоремы \ref{th3_1}. Для этого требуется, чтобы для любой убывающей и стабилизирующейся
к $x \in X_0$ последовательности  $(x^{(k)})$ из  $X$ было выполнено неравенство \eqref{eq3_7}. По определению  \ref{df3_2} 
\begin{equation}\label{eq5_2}
\inf_{k\in \NN} x^{(k)}=x\in X_0,\;
x^{(k+1)}\leq x^{(k)},\;
x^{(k)}=\sup_{m\geq k} x^{(m)},
 \; \pr_n x^{(k)}=
\pr_n x^{(n)}\text{ при всех $k\geq n\geq 1$}.
\end{equation}
По определению \ref{defsprf} супремальной функции  
\begin{multline}\label{eq5_3}
f(x^{(k)})=\spf_{H,q}(x^{(k)})\overset{\eqref{eq4_1}}{=}\sup \bigl\{ q(h)\colon  h\leqslant x^{(k)}, \; h\in H \bigr\}
\\=\sup \bigl\{ q_1(\pr_1 h)\colon  h\leqslant x^{(k)},
\; h\in H \bigr\}=:a_k \in \KK_{-\infty}.
\end{multline}
Последовательность $ (a_k)_{k\in \NN}$ --- убывающая. При $\inf_{k\in \NN} a_k =-\infty$ условие \eqref{eq3_7} выполнено автоматически. Поэтому далее убывающая последовательность $ (a_k)_{k\in \NN}$ ограничена в $\KK$ и  
\begin{equation}\label{eqa0}
a_0\overset{\eqref{eq5_2}}{:=}\inf\limits_{k\in \NN} a_k 
\overset{\eqref{eq5_3}}{=}\inf f(x^{(k)})\in \KK ,
\quad a_0\leq a_k \leq a_1 \quad \text{для всех $k\in \NN$}.
\end{equation}
По условию \eqref{qHiii-} существует убывающая в $X$ последовательность
\begin{equation}\label{hk}
h^{(k)}:=\sup \{h\in H\colon h\leq x^{(k)}\}\in H,\quad 
k\in \NN, \quad h^{(k)}\leq x^{(k)},
\end{equation}
для которой в силу возрастания $q_1$ и $q=q_1\circ \pr_1$ имеем равенства 
\begin{equation}\label{akhk}
a_k\overset{\eqref{eq5_3}}{=}q(h^{(k)})=q_1(\pr_1 h^{(k)}), \quad k\in \NN.
\end{equation}
Пусть сначала выполнено \eqref{qh1}, т.\,е.  $ \inf_{k\in \NN} q(h^{(k)})\in \KK$. Тогда по условию \eqref{qHv} последовательность \eqref{seqhn}  ограничена снизу в $X$ и по условию  \eqref{qHiv} существует 
\begin{equation}\label{eq5_7}
h=\inf_{k\in \NN} h^{(k)} \in H,
\end{equation}
а  из \eqref{hk} и стабилизации последовательности 
$(x^{(k)})$ к $x\in X_0$ следует
\begin{equation}\label{eq5_8}
h\leqslant x.
\end{equation}
Из условия \eqref{qHv} доказываемой теоремы и из \eqref{hk} по построению \eqref{eq5_7} вектора  $h\in H$ вытекает
\begin{equation*}
q(h)=q\Bigl(\inf_{k\in \NN} h^{(k)}\Bigr)\overset{\eqref{qHv}}{\geqslant} \inf_{k\in \NN} q(h^{(k)})\overset{\eqref{akhk}}{=}
 \inf_{k\in \NN} a_k\overset{\eqref{eqa0}}{=}a_0.
\end{equation*}
Следовательно, согласно \eqref{eq5_3} и \eqref{eq5_8}, для некоторого $h\in H$, $h\overset{\eqref{eq5_8}}{\leqslant} x$, выполнено
\begin{equation*}
q(h)\geqslant \inf_{k\in \NN} \spf_{H,q} (x^{(k)}).
\end{equation*}
По определению \ref{defsprf} супремального функционала $\spf_{H,q}$ и ввиду \eqref{eq5_3} это означает, что
\begin{equation}\label{eq5_9}
\inf f(x^{(k)})=\inf\limits_k \spf_{H,q}(x^{(k)})\leqslant
\spf_{H,q} \Bigl(\inf_{k\in \NN} x^{(k)}\Bigr)=
f\Bigl(\inf_{k\in \NN} x^{(k)}\Bigr).
\end{equation}
Таким образом, в случае ограниченной снизу последовательности $(h^{(k)})_{k\in \NN}$ выполнено \eqref{eq3_7}.

Если же  \eqref{qh1} не выполнено, т.\,е.  $ q(h^{(1)})\geq \inf_{k\in \NN} q(h^{(k)})= -\infty$ и по условию  
\eqref{qHiv} последовательность \eqref{seqhn} 
не ограничена снизу, то согласно \eqref{akhk} имеем $\inf_{k\in \NN} a_k=-\infty$. Как уже отмечалось выше, при $\inf_{k\in \NN} a_k =-\infty$ условие \eqref{eq3_7}, или \eqref{eq5_9}, выполнено автоматически.

По теореме \ref{th3_1} в виде следствия \ref{corspl}
заключаем, что суперлинейная функция $f=\spf_{H,q}\in
\spl \KK_{-\infty}^X$ допускает описание через верхнюю огибающую на $X_0$:
\begin{equation*}
\spf_{H,q}(x)=\inf \bigl\{ l(x)\colon l\in L, \;
\spf_{H,q}\leqslant l \text{ на } X_0 \bigr\}, \quad x\in X_0.
\end{equation*}
По равенству \eqref{spr_spf} предложения \ref{pr_pf} правая часть здесь совпадает с правой частью \eqref{eq5_1}.
\end{proof}

\begin{corollary}[{\rm развитие \cite[теорема 5.1]{Kh01_1}}]\label{cor_RR} Пусть $\KK:=\RR$ и выполнены все условия теоремы\/ {\rm \ref{thrdq}} за исключением условий  \eqref{qHiii-} и \eqref{qHiv}, которые заменяем на одно условие 
\begin{enumerate}
\item[{\rm (iii--iv)}]\label{iii_iv} для любой ограниченной в $X$ последовательности \eqref{seqhn} существует верхний предел
\begin{equation}\label{limsh}
\limsup_{k\to \infty} h^{(k)} \overset{\eqref{lsxk}}{:=}\inf_{n\in \NN}\sup_{k\geqslant n} h^{(k)} \in H \quad\text{\rm (см. определение \ref{df3_2})}.
\end{equation}
\end{enumerate}
Тогда функция $\spf_{H,q}$ такая же, как в теореме {\rm \ref{thrdq}} и выполнено заключение \eqref{eq5_1}, где $L:=\lin \RR^{X_0}$ можно заменить на  
$\lin^+ \RR^{X_0}$, если в $H$ есть отрицательный вектор, т.\,е. $H\cap (-X^+)\neq \varnothing$.
\end{corollary}
\begin{proof} Будет использована 
\begin{lemma}\label{lemiii_iv} Для каких бы то ни было условий на подмножество $H\subset X$ условие\/ {\rm (iii--iv)} эквивалентно сочетанию условий
\eqref{qHiv} теоремы\/ {\rm \ref{thrdq}} и 
\begin{enumerate}
\item[{\rm ($\rm iii'$)}]  {\it для любой  ограниченной сверху в $X$ последовательности\/ \eqref{seqhn} имеем\/ $\sup_{k\in \NN} h^{(k)} \in H$.}
\end{enumerate}
\end{lemma}
\begin{proof}[Доказательство леммы\/ {\rm \ref{lemiii_iv}}]
Достаточность условий \eqref{qHiv} теоремы\/ {\rm \ref{thrdq}} и ($\rm iii'$) для выполнения условия (iii--iv) очевидна по определению верхнего предела в \eqref{limsh}, или в \eqref{lsxk}.

Обратно, пусть выполнено условие условия (iii--iv). Если $h, h'\in H$, то для счётной последовательности  $b:=(h,h',h,h', \dots)$ чередующихся пар $h,h'$ по условию 
(iii--iv) имеем 
\begin{equation*}
\sup\{h,h'\}=\sup b\overset{\eqref{lsxk}}{=}
\limsup b\overset{\eqref{limsh}}{\in} H.  
\end{equation*}
Таким образом, точная верхняя граница любого конечного подмножества из $H$ принадлежит $H$. Отсюда, если последовательность  \eqref{seqhn} ограничена сверху, то новая  последовательность $\bigl(\sup_{k\leq n} h^{(k)} \bigr)_{n\in \NN}\subset H$ ограничена и по условию (iii--iv) 
\begin{equation*}
\sup_{k\in \NN} h^{(k)}=\sup_{n\in \NN}
\Bigl\{\,\sup_{k\leq n} h^{(k)}\Bigr\}
\overset{\eqref{lsxk}}{=}\limsup_{n\to \infty} \Bigl(\sup_{k\leq n} h^{(k)}\Bigr) \overset{\eqref{limsh}}{\in}
H.
\end{equation*}
Это означает, что выполнено условие ($\rm iii'$). 
Кроме того, для любой убывающей ограниченной снизу, а значит и ограниченной, последовательности \eqref{seqhn} её верхний предел, если он существует, равен точной нижней границе этой последовательности. Таким образом, из условия (iii--iv) следует условие \eqref{qHiv} теоремы \ref{thrdq}.  
\end{proof}

Условие ($\rm iii'$) слабее условия \eqref{qHiii-} теоремы \ref{thrdq}, поэтому нам придётся изменить часть доказательства теоремы \ref{thrdq}.
Далее приведём лишь те этапы доказательства, которые 
для $\KK=\RR$ отличаются от соответствующих частей 
доказательства теоремы \ref{thrdq}. Изменения начинаем сразу после \eqref{eqa0}. 

Пусть $\varepsilon \in \RR_*^+$. Ввиду \eqref{eqa0} и \eqref{eq5_3} найдётся вектор $h_{\varepsilon}^{(k)}\in H$, для которого
\begin{equation}\label{eq5_4}
q(h_{\varepsilon}^{(k)})\geq a_k-\varepsilon, \quad h_{\varepsilon}^{(k)}\leq x^{(k)}\leq x^{(1)}, \quad k\in \NN. 
\end{equation}
Так как последовательность $(h_{\varepsilon}^{(k)})_{k\in \NN}$  ограничена сверху, по условию ($\rm iii'$) существуют
векторы 
\begin{equation}\label{an6_10}
h^{(k)} :=\sup_{n\geq k} h_{\varepsilon}^{(n)}\in H, 
\quad h_{\varepsilon}^{(k)}\leq h^{(k)}\overset{\eqref{eq5_4}}{\leq} x^{(k)}\leq x^{(1)}, \quad k\in \NN,
\end{equation}
которые образуют убывающую последовательность  \eqref{seqhn}.
При этом согласно \eqref{eq5_4} ввиду стабилизации  $(x^{(k)})_{k\in \NN}\subset X$ из \eqref{eq5_2} к $x=(x_n)_{n\in \NN}\in X_0$, $x_n:=\pr_n x\in X_n$, имеем 
\begin{equation}\label{eq5_5-6}
\pr_n h^{(k)}\leq_n \pr_n x^{(k)}=x_n \quad\text{при $k\geq n\in \NN$}, \quad q(h^{(k)})\geq q(h^{(k)}_{\varepsilon})
\overset{\eqref{eq5_4}}{\geq} a_k-\varepsilon. 
\end{equation}

Пусть сначала выполнено \eqref{qh1}, т.\,е.  $ \inf_{k\in \NN} q(h^{(k)})\in \RR$. Тогда по условию \eqref{qHv} последовательность \eqref{seqhn}  ограничена снизу в $X$ и по условию  \eqref{qHiv} имеем \eqref{eq5_7},
а  из \eqref{an6_10} и стабилизации последовательности 
$(x^{(k)})$ к $x\in X_0$ следует \eqref{eq5_8}, т.\,е.
$h\leqslant x$. Из условия \eqref{qHv} доказываемой теоремы и из \eqref{eq5_4} по построению \eqref{eq5_7} вектора  $h\in H$ вытекает
\begin{equation*}
q(h)\overset{\eqref{eq5_7}}{=}q\Bigl(\inf_{k\in \NN} h^{(k)}\Bigr)\overset{\eqref{qHv}}{\geqslant} \inf_{k\in \NN} q(h^{(k)})\overset{\eqref{an6_10}}{\geq} \inf_{k\in \NN} q(h^{(k)}_\varepsilon) \overset{\eqref{eq5_4}}{\geq} 
 \inf_{k\in \NN} a_k-\varepsilon  \overset{\eqref{eqa0}}{=}a_0-\varepsilon.
\end{equation*}
Отсюда для произвольного $\varepsilon\in \RR_*^+$ 
при некотором, зависящем от $\varepsilon$, $h\in H$, $h\overset{\eqref{eq5_8}}{\leqslant} x$, выполнено
\begin{equation*}
q(h)\geqslant \inf_{k\in \NN} \spf_{H,q} (x^{(k)})-\varepsilon.
\end{equation*}
По определению \ref{defsprf} супремального функционала $\spf_{H,q}$ согласно  \eqref{eq5_3} это означает, что
\begin{equation}\label{eq5_9e}
\inf f(x^{(k)})=\inf\limits_k \spf_{H,q}(x^{(k)})\leqslant
\spf_{H,q} \Bigl(\inf_{k\in \NN} x^{(k)}\Bigr)+\varepsilon=
f\Bigl(\inf_{k\in \NN} x^{(k)}\Bigr)+\varepsilon.
\end{equation}
Это, ввиду произвола в выборе $\varepsilon\in \RR_*^+$, даёт  
\eqref{eq5_9}. Теперь повторяем рассуждения и выкладки после \eqref{eq5_9} до конца доказательства теоремы  \ref{thrdq}.
Осталось заменить $\lin \KK^{X_0}$ на $\lin^+ \RR^{X_0}$
в \eqref{eq5_1}.
\begin{lemma}\label{lemp} Пусть в условиях и обозначениях теоремы\/ {\rm \ref{thrdq}} до п.~\eqref{qHii} включительно $\KK:=\RR$, $l\in \lin \RR^{X_0}$ и $q\prec_H^{X_0}l$, $H\cap (-X^+)\neq \varnothing$. Тогда $l\in \lin^+\RR^{X_0}$, т.\,е. функция $l$ положительна на $X_0$.
\end{lemma}
\begin{proof}[Доказательство леммы\/ {\rm \ref{lemp}}] Поскольку в $H$ существует 
отрицательный вектор, из строго положительной однородности $\spf_{H,q}$ по предложению \ref{pr_spf}, п.~\ref{Xosii}, 
имеем $\spf_{H,q}(0)=0$. Из равенства \eqref{spr_spf} предложения \ref{pr_pf} и определения \ref{balpb} ростка 
в \eqref{sprqb} имеем $\spf_{H,q}\leq l$ на $X_0$. 
По предложению \ref{pr_spf}, п.~\ref{Xosiv}, ввиду 
условия $H\cap (-X^+)\neq \varnothing$ имеем $\spf_{H,q}\in \spl^+ \RR_{-\infty}^{X}$. Пусть $x\in X_0^+$. Тогда
\begin{equation*}
0=\spf_{H,q}(0)\leq \spf_{H,q}(x)\leq l(x), 
\end{equation*}
что доказывает положительность $l$ на $X_0$.
\end{proof}
Применение леммы \ref{lemp} завершает доказательство следствия
\ref{cor_RR}. 
\end{proof}

\subsection{Векторно-аффинная версия для выпуклого множества $H$}\label{vec_ver}

\begin{theorem}\label{thrdv} 
Пусть $X$ и $X_0\subset X$ те же, что и в теореме\/ 
{\rm \ref{thrdq}}, но $q_1\in \incr \conc \RR^{X_1}$ и $q\overset{\eqref{q1qn}}{:=}q_1\circ \pr_1$, $H\subset X$ --- выпуклое множество, $H_n=\pr_n H\subset X_n$, $\KK:=\RR$. Пусть выполнены условия \eqref{qHi} и \eqref{qHv} из теоремы\/  {\rm \ref{thrdq}}, а также условие\/ {\rm (iii--iv)} из следствия\/ {\rm \ref{cor_RR}}. Тогда  для $L:=\aff \RR^{X_0}$ вогнутая 
функция  $\spf_{H,q}\in \conc \RR_{-\infty}^X$ допускает представление \eqref{eq5_1}. 

Если $0\in H$ и $q_1(0)\geq 0$, то $\spf_{H,q}\in \conc^+\RR_{-\infty}^X$, а  $L=\aff \RR^{X_0}$ в \eqref{eq5_1} можно заменить на \begin{equation}\label{aff+}
L:=\aff^+\RR^{X_0}\overset{\eqref{affl}}{=}\lin^+\RR^{X_0}+\RR^+ \cdot \boldsymbol{1} .
\end{equation} 
\end{theorem}
\begin{proof} Будет использована теорема \ref{th3_1} в виде следствия \ref{coraff} из неё. 

Положим $f:=\spf_{H,q}$ и $f_n=\sup\tmin\pr_n f$, $L:=\aff \KK^{X_0}$ и  $L_n:= \aff \KK^{\pr_n X_0}$. Тогда по предложению \ref{prLaff} выполнено \eqref{LnL} и   
по предложению \ref{prqn} в обозначениях \eqref{q1qn} имеет место равенство \eqref{spflln+}, где по предложению \ref{pr_spf}, п.~\ref{Xosiii}, имеем  $\spf_{H_n,q_n}\in \conc  \KK_{-\infty}^{X_n}$, т.\,е. $f_n \overset{\eqref{spflln+}}{\in} \conc \KK_{-\infty}^{X_n}$, 
а также  $f=\spf_{H,q}\in \conc \RR_{-\infty}^X$ . Напомним, что $\pr_nX_0$ --- векторное подпространство в $X_n$ по предложению \ref{prLlin}.
Рассмотрим сужения  $f_n\bigm|_{\pr_nX_0} =\spf_{H_n,q_n}\bigm|_{\pr_nX_0}$. В правой части здесь --- вогнутая функция на $\pr_nX_0$ со значениями исключительно в $\RR$ по определению
\ref{defsprf}, равенство \eqref{eq4_1}, супремальной функции
$\spf_{H_n,q_n}\bigm|_{\pr_nX_0}$  и по условию
\eqref{qHi} из теоремы  {\rm \ref{thrdq}} о минорировании. Таким образом, 
$f_n\bigm|_{\pr_nX_0} \in \conc \KK^{\pr_n X_0}$. 
Тогда по теореме Хана\,--\,Банаха об огибающей из введения, применённой в части импликации \ref{HBci}$\Rightarrow$\ref{HBcii} в <<вогнутой форме>> в векторном пространстве $\pr_n X_0$, имеем 
$f_n\bigm|_{\pr_nX_0}\overset{\eqref{luEu}}{=}{\uE}^{\aff \KK^{\pr_nX_0}}_{f_n|_{\pr_nX_0}}$
на $\pr_n X_0$, откуда выполнено условие \eqref{drqn} теоремы \ref{th3_1}.

Функция $f=\spf_{H,q}$ --- возрастающая на $X$ по п.~\ref{fHpII} предложения \ref{pr_spf}. Проверим выполнение условия \eqref{stabii} теоремы \ref{th3_1}. Для этого требуется, чтобы для любой убывающей и стабилизирующейся
к $x \in X_0$ последовательности  $(x^{(k)})$ из  $X$ было выполнено неравенство \eqref{eq3_7}. По определению  \ref{df3_2} имеем \eqref{eq5_2}.
По определению \ref{defsprf} супремальной функции  также 
выполнено \eqref{eq5_3}.

Последовательность $ (a_k)_{k\in \NN}$ --- убывающая. При $\inf_{k\in \NN} a_k =-\infty$ условие \eqref{eq3_7} выполнено автоматически. Поэтому далее убывающая последовательность $ (a_k)_{k\in \NN}$ ограничена в $\RR$ и
выполнено  \eqref{eqa0} с $\KK:=\RR$.    
Из условия (iii-iv) по лемме \ref{lemiii_iv} следуют условия \eqref{qHiv} теоремы\/ {\rm \ref{thrdq}} и ($\rm iii'$)
леммы \ref{lemiii_iv}. 

Пусть $\varepsilon \in \RR_*^+$. Ввиду \eqref{eqa0} и \eqref{eq5_3} найдётся вектор $h_{\varepsilon}^{(k)}\in H$, для которого выполнено \eqref{eq5_4}. Далее дословно повторяем рассуждения и выкладки из доказательства следствия 
от \eqref{eq5_4} до \eqref{eq5_9e} включительно. Последнее
ввиду произвола в выборе $\varepsilon\in \RR_*^+$ даёт  
\eqref{eq5_9}.
Теперь остаётся лишь повторить рассуждения и выкладки после \eqref{eq5_9} до конца доказательства теоремы  \ref{thrdq}, но с вогнутой функцией $f=\spf_{H,q}\in 
\conc \RR_{-\infty}^X$, что даёт \eqref{eq5_1} для $L:=\aff \RR^{X_0}$.  

При $0\in H$ и $q_1(0)\geq 0$, т.\,е. $q(0)\geq 0$, по предложению \ref{pr_spf}, п.~\ref{Xosi}, $\spf_{H,q}$ --- положительная функция. Тогда по предложению \ref{pr_pf} из равенства \eqref{spr_spf} следует, что для $l\in
\Spr(q; H,X_0,\aff \RR^{X_0})$ имеем $l\geq \spf_{H,q}$ на $X_0$, откуда $l\in (\RR^{X_0})^+$. Таким образом, $\aff \RR^{X_0}$ можно заменить на $\aff^+\RR^{X_0}$. Последнее равенство в \eqref{aff+} составляет содержание 
пунктов \eqref{ali} и \eqref{aliK} предложения \ref{afflin}.
\end{proof}

\subsection{Топологические версии для выпуклого  $H$}\label{top_ver}

Для топологических вариантов теорем  о   представлении супремальных функций  потребуется понятие решётки Фреше  \cite{Schef}.

Пусть $X$  ---  векторная  решётка.  Для  $x\in  X$,  как
обычно, полагаем  $|x|=\sup \{ -x,  x\}$ --- абсолютная величина вектора $x$. Если  к  тому  же  $X$  ---  полное  метризуемое  локально выпуклое пространство,  т.\,е.  пространство  Фреше,  и обладает таким базисом окрестностей нуля,  что  для  любой  окрестности нуля $V$ из этого базиса при любом $x\in V$ неравенство $|x'|\leqslant |x|$ влечёт за собой $x'\in  V$,  то  $X$  называется {\it решёткой
Фреше.} Решётка   Фреше  ---  упорядоченное
локально выпуклое  пространство,  т.\,е.   конус   положительных
векторов в нём замкнут \cite[гл.~V, 7.2]{Schef}.

Кроме того, если $p$  --- линейное отображение решётки Фреше в решётку Фреше, то $p$   непрерывно   \cite[гл.V, 6.4, 6.1]{Schef}. В частности, непрерывны линейные положительные функционалы на решётке Фреше. Таким образом,   если  $X=\projlim_n  X_np_n$  ---  приведённый правильный проективный предел решёток Фреше как упорядоченных векторных пространств,  то отображения $p_n$ непрерывны  и  порождают  на  $X$ топологию проективного предела \cite[гл.~II]{Schef}, \cite{SiS}. 
Топология    на    $X$ определяется как индуцированная с топологического произведения $\prod_{n} X_n$.  Если   проекции   ${\pr}_n$   рассматривать   как
определённые на  $X$,  то  ${\pr}_n$  ---  линейные  непрерывные
отображения и в качестве базиса  окрестностей  нуля  в  $X$
можно взять    прообразы ${\pr}^{-1}_n (V_n)$ всевозможных   окрестностей   нуля   $V_n$ пространств $X_n$.   Пространство   $X$   как   топологический проективный предел  локально выпуклых хаусдорфовых пространств также будет   локально   выпуклым   и   хаусдорфовым.

Конус положительных элементов в $X$ замкнут в $X$.
Действительно, если $K_n$ --- конусы положительных векторов
в решётках Фреше $X_n$ , то они замкнуты в $X_n$.
Конус положительных элементов $K$ в $X$ --- это пересечение
$\bigcap\limits_n {\pr }_n^{-1}(K_n)$, замкнутых в силу непрерывности отображений ${\pr }_n$ множеств ${\pr }_n^{-1}(K_n)$, т.\,е. конус $K$ замкнут в $X$. Следовательно, $X$ --- упорядоченное локально выпуклое пространство.
Всюду далее проективный предел решёток Фреше рассматривается
одновременно и как упорядоченное локально выпуклое пространство
с топологией проективного предела, описанной выше.

Пусть $S$ --- подмножество в локально выпуклом пространстве. $S$ называется {\it секвенциально замкнутым,\/} если для любой сходящейся последовательности векторов из $S$ предел этой последовательности принадлежит $S$. $S$\/  называется {\it секвенциально предкомпактным,\/} если  любая последовательность векторов из $S$ содержит сходящуюся подпоследовательность.

\begin{theorem}[{\rm вариация \cite[теорема 6.1]{Kh01_1}}]\label{th6_1} 
Пусть $X=\projlim_n X_np_n$ ---  приведённый правильный проективный предел  решёток Фреше $X_n$, $n\in \NN$, и $X_0$ --- векторное подпространство в $X$, $q_1\in \lin^+\RR^{X_1}$ и $q\overset{\eqref{q1qn}}{:=}q_1\circ \pr_1$, выпуклый конус $H\subset X$ секвенциально замкнут в $X$ и
$H\cap (-X^+)\neq \varnothing$.  Допустим, что при каждом $n\in \NN$ выполнено одно из двух условий:
\begin{enumerate}[{\rm (i)}]
\item\label{ki} $H_n={\pr}_n H$ --- минорирующее множество для 
${\pr}_n X_0$,
\item\label{kii} $\sup$-проекция   $f_n\overset{\eqref{eq3_3}}{:=}\sup\tmin\pr_n \spf_{H,q}$  полунепрерывна сверху на $\pr_n X_0$, 
\end{enumerate}
а также выполнено ещё и условие
\begin{enumerate}[{\rm (iii)}]
\item\label{spk_iii} для любого не более чем счётного
подмножества $B\subset H$, ограниченного сверху
и удовлетворяющего условию $\inf q(B)>-\infty$,
найдётся $n_q\in \NN$, для которого при каждом 
$n\geqslant n_q$ проекция ${\pr}_n B\subset X_n$
секвенциально предкомпактна в $X_n$.
\end{enumerate}
Тогда для $L:=\lin^+\RR^{X_0}$ и супремальной функции  $\spf_{H,q}\in \spl^+\RR_{-\infty}^X$ справедливо заключение теоремы\/ {\rm \ref{thrdq}} вместе c  представлением \eqref{eq5_1} на $X_0$.
\end{theorem}
\begin{proof} 
В обозначениях  $L:=\lin \RR^{X_0}$ и  $L_n:= \lin \RR^{\pr_n X_0}$ так же, как и при доказательстве  
теоремы \eqref{thrdq}, убеждаемся, что выполнено 
\eqref{LnL}, функция $f:=\spf_{H,q}\in \spl^+\RR_{-\infty}^X$ возрастающая  и для $f_n$ из \eqref{spflln+}, или из условия \eqref{kii},  имеем $f_n \in \spl^+\RR_{-\infty}^{X_n}$.  Если выполнено условие минорирования \eqref{ki}, то $f_n \in \spl^+\RR^{X_n}$ и  по теореме Хана\,--\,Банаха об огибающей из введения, применённой в части импликации \ref{HBi}$\Rightarrow$\ref{HBii}, получаем условие \eqref{drqn} теоремы \ref{th3_1} на векторном подпространстве $\pr_n X_0$ при $n\in \NN$. В случае же выполнения  условия 
\eqref{kii} по теореме Хёрмандера об огибающей из введения, применённой в части импликации \ref{HI1}$\Rightarrow$\ref{HI2}, по-прежнему получаем условие \eqref{drqn} теоремы \ref{th3_1} при $n\in \NN$. Остаётся проверить условие  \eqref{stabii} 
теоремы \ref{th3_1} для убывающей стабилизирующейся к $x=(x_n)_{n\in \NN}\in X_0$ последовательности \eqref{eq5_2}. При этом по замечанию \ref{prl_subs} можем в условии (iii) положить $n_q=1$. 

Если для последовательности $(a_k)$, определённой в \eqref{eq5_3}, $\inf_{k\in \NN} a_k =-\infty$, то условие \eqref{eq3_7} выполнено очевидным образом. Пусть
\begin{equation}\label{eq6_1}
\inf_{k\in\NN} a_k \overset{\eqref{eqa0}}{=:}a_0 >-\infty.
\end{equation}
Тогда согласно \eqref{eq5_3} найдётся последовательность
$(h^{(k)})_{k\in \NN}\subset H$, для которой  выполнено \eqref{eq5_4}. При этом ввиду \eqref{eq6_1} для ограниченного сверху множества $B=\{ h^{(k)}\colon k\in \NN \}$ выполнено
условие (iii), значит все проекции ${\pr}_n B$ секвенциально
предкомпактны в $X_n$. Следовательно, в множестве $B$
можно выбрать последовательность $(h^{(1{,}n)})$ , сходящуюся в $X_1$; из последовательности $(h^{(1{,}n)})$ можно выбрать подпоследовательность $(h^{(2{,}n)})$, сходящуюся в $X_2$  и т.д. Диагональный процесс даёт последовательность $(h^{(n{,}n)})$, сходящуюся в каждом $X_n$, $n\in \NN$.
Из этого по определению топологии проективного предела вытекает
сходимость последовательности $(h^{(n{,}n)})$ в $X$.
Так как по построению $(h^{(n{,}n)})\subset H$,
то в силу секвенциальной замкнутости $H$ в $X$ предел $h$
последовательности $(h^{(n{,}n)})$ лежит в $H$.
Из включения
\begin{equation}\label{eq6_2}
(h^{(n{,}n)}) \subset (h^{(k)})
\end{equation}
и неравенств $h^{(k)}\leqslant x^{(k)}$ из \eqref{eq5_4} ввиду стабилизации $(x^{(k)})$ к $x=(x_n)$ имеем
\begin{equation*}
{\pr}_k\, h^{(n{,}n)}\leq_k {\pr }_k x^{(k)} =x_k,
\quad k\in \NN. 
\end{equation*}
Отсюда ввиду замкнутости конуса положительных векторов в каждом
$X_k$ и непрерывности проекций ${\pr}_k$ получаем
${\pr}_k h \leq_k {\pr}_k x$, $k\in \NN$. Последнее
означает, что $h\leqslant x$ в $X$.

Пусть $q_k$ --- функция, определённая в \eqref{q1qn}   на решётке Фреше $X_k$. Так как $q_k$ --- линейный положительный функционал на решётке Фреше, то $q_k$ --- непрерывный функционал. Отсюда 
\begin{equation*}
q(h)=q_k\bigl({\pr }_k h\bigr)=
q_k \bigl({\pr }_k \lim_n h^{(n{,}n)}\bigr)
=\lim_n q_k \bigl({\pr }_k  h^{(n{,}n)}\bigr)
=\lim_n q \bigl( h^{(n{,}n)}\bigr).
\end{equation*}
Отсюда в силу включения \eqref{eq6_2} и выполнения первого неравенства из \eqref{eq5_4} имеем $q(h)\geqslant a_0-\varepsilon$, и как показано выше, $h\in H$, $h\leqslant x$. Следовательно, $\spf_{H,q}(x) \geqslant a_0$. Вспоминая определение
 величины $a_0$ из \eqref{eq6_1} и \eqref{eq5_3}, а также \eqref{eq5_2}, видим, что  выполнено \eqref{eq5_9}. По теореме \ref{th3_1} и предложению \ref{pr_pf} функция $\spf_{H,q}$ допускает представление  вида \eqref{eq5_1} с $L:=\lin \RR^{X_0}$. Возможность замены здесь $L:=\lin \RR^{X_0}$ на $L:=\lin^+ \RR^{X_0}$ следует из леммы \ref{lemp}.
\end{proof}

Приведённая ниже версия теоремы \ref{th6_1} для выпуклого множества $H$ ранее не отмечалась.

\begin{theorem}[{\rm вариация \cite[теорема 6.1]{Kh01_1}}]\label{th6_2}
Пусть $X=\projlim_n X_np_n$ ---  приведённый правильный проективный предел  решёток Фреше $X_n$, $n\in \NN$, и $X_0$ --- векторное подпространство в $X$, $q_1\in \aff^+\RR^{X_1}$ и $q\overset{\eqref{q1qn}}{:=}q_1\circ \pr_1$, выпуклое подмножество $H\subset X$ секвенциально замкнуто в $X$.
Допустим, что при каждом $n\in \NN$ выполнено одно из двух условий \eqref{ki} или \eqref{kii} теоремы\/ {\rm \ref{th6_1}}, а также выполнено ещё и условие\/  {\rm (iii)} теоремы\/ {\rm \ref{th6_1}}. Тогда  для $L:=\aff \RR^{X_0}$ вогнутая 
функция  $\spf_{H,q}\in \conc \RR_{-\infty}^X$ допускает представление \eqref{eq5_1}. 
Если дополнительно $0\in H$ и $q_1(0)\geq 0$, то $\spf_{H,q}\in \conc^+\RR_{-\infty}^X$, а  $L=\aff \RR^{X_0}$ в \eqref{eq5_1} можно заменить на $L:=\aff^+\RR^{X_0}$ из \eqref{aff+}.
\end{theorem}
Доказательство теоремы \ref{th6_2} опускаем, поскольку оно представляет из себя симбиоз доказательств теорем  \ref{thrdv} и \ref{th6_1}. Отметим лишь, что при доказательстве условия 
\eqref{drqn} теоремы \ref{th3_1} на векторных подпространствах  $\pr_n X_0$ при $n\in \NN$ в случае выполнения  условия 
\eqref{kii} нужно использовать теорему Хёрмандера об огибающей из введения в части импликации \ref{HII1}$\Rightarrow$\ref{HII2}. 

\begin{remark}\label{rx0up} Это замечание относится ко всем четырём теоремам \ref{thrdq}--\ref{th6_2}. Заключение этих теорем о возможности представления \eqref{eq5_1} для векторов $x$ из векторного подпространства $X_0$ c помощью утверждений  типа  теорем о минимаксе можно распространить на векторы $x\in X$, представимые в виде точных верхних границ возрастающих последовательностей векторов из $X_0$ \cite{KKh09}, \cite[2.5]{Kh10}. Доказательство этого перехода довольно объёмно. Полное изложение возможности такого распространения представления \eqref{eq5_1} намечается обосновать в ином месте.
\end{remark}

\chapter{Применения в теории функций}\label{gh2}

\section{Введение. Где возникают огибающие?}\label{int_hf}
\setcounter{equation}{0}

Кратко остановимся на схемах применения двойственных описаний верхних огибающих из главы~\ref{ch_en} в комплексном анализе.
В основном следуем изложению таких схем в \cite{Kh01_1}--\cite{Kh10}.  

Всюду  во  этом разделе   $D$ --- {\it область\/} (открытое связное подмножество)  в комплексной плоскости $\CC$ или  $n$-мерном комплексном пространстве  $\CC^n$ со стандартной евклидовой топологией, $n \in \NN$. 
Через  $\Hol (D)$, $\Mer (D)$, $\har (D)$, $\phar (D)$, $\sbh (D)$ и $\psbh (D)$ обозначаем классы соотв. голоморфных, мероморфных, гармонических, плюригармонических, субгармонических и плюрисубгармонических  функций на $D$; последние два класса по определению содержат функцию $\boldsymbol{-\infty}$ --- тождественную $-\infty$;
$\sbh (D)=\psbh (D)$ и  $\har (D)=\phar (D)$ при $n=1$.

Для $S\subset \CC^n$ различаем   $C_{\RR}(S)$ и  $C_{\CC}(S)$ --- классы непрерывных вещественнозначных и соответственно комплекснозначных функций.

Всюду во введении $\sf Z$ ---  главное аналитическое множество в $D$, т.\,е. нулевое множество некоторой ненулевой функции $g_{\sf Z}\in \Hol (D)$, заданное вместе с функцией кратности нулей, или с дивизором нулей,  функции $g_{\sf Z}$. При этом и
дивизор нулей обозначаем той же буквой  $\sf Z$.

 В случае $n=1$ часто удобно мыслить $\sf Z$ как последовательность, вообще говоря, повторяющихся точек $\{{\sf z}_k\}_{k\in \NN}$, для которой порядок нумерации не важен, а имеет значение только число повторений точки в ней.  Если $\sf Z$ --- последовательность всех нулей (корней) функции $f\in \Hol (D)$, то  каждая точка ${\sf z} \in D$ повторяется $\sf Z$ в столько раз, какова кратность нуля (корня) функции $f$ в точке $\sf z$. При этом полагаем  $\Zero_f:={\sf Z}$. Отношения и операции для множеств (последовательностей) нулей понимаются как в \cite{Kh01_1}--\cite{Kha12}. 

Различные задачи комплексного анализа сводятся к построению или доказательству существования (точной) огибающей --- верхней или нижней --- из определенного класса функций на $D$ или на подмножестве $S \subset D$. Отметим некоторые из них в подходящей трактовке.

\subsection{Нетривиальность весового класса}\label{nontr7_1} 
По этому подразделу см. {\rm  \cite[\S~10]{Kh01_2}, \cite[1.1]{Kh10}}. При каких условиях на {\it функцию-мажоранту\/} ({\it весовую функцию, вес\/})\/ $M\colon D\to \RR_{\pm \infty}$ найдется  ненулевая функция  $f\in \Hol(D)$ с ограничением $\log |f|\leqslant M$ на $D$? Другими словами,  вопрос состоит в исследовании условий нетривиальности класса функций 
$\Hol(D, M):=\{f\in \Hol(D)\colon \log|f|\leq M +\const \text{ на }D\}$,
т.\,е. условий, при которых  $\Hol(D, M)\neq \{0\}$.
Зачастую достаточно убедиться в  нетривиальности выпуклого множества
\begin{equation}\label{env:1}
\{h\in \psbh (D)\colon h\leq M \text{ на }D\}, 
\end{equation}
т.\,е. доказать существование функции $h\neq\boldsymbol{-\infty}$ из этого класса, а затем для такой функции $h\leq M$ попытаться построить  голоморфную функцию $f\not\equiv 0$, удовлетворяющую оценке $\log |f| \leq \check h\leq \check M$, где $\check h$ и  $\check M$ --- некоторые незначительные увеличения соответственно функций $h$ и $M$. При этом существование такой функции $f$ можно доказать при помощи решений $\overline\partial$-  или $\partial\overline\partial$-задачи с оценками (см. \cite[\S~9]{Kh01_2}, \cite[IIIA]{Koo}, \cite[теорема 1]{BaiKha16}, \cite[теорема 3, следствия 1--3]{KhaBai16}) в духе Л.~Хёрмандера \cite[гл.~IV]{HorC} или же путем аппроксимации функции $h$ логарифмом модуля голоморфной функции. Правда, последний аппроксимационный способ представляется нам не соответствующим  постановке задачи и в данной тематике излишним, поскольку здесь требуется только минорирование функции $h$, а не ее аппроксимация.           

\subsection{Описание нулевых множеств}\label{zero7_2}
По этому подразделу см. {\rm \cite[\S~8]{Kh01_2}, \cite{Kha07}, \cite{KhCh}, \cite[1.2]{Kh10}.}  
Пусть ${\sf Z}$ --- нулевое множество некоторой функции $g_{\sf Z}\in \mathrm{Hol}\,(D)$. Если $D$ --- односвязная при $n=1$ или звёздная относительно точки при $n>1$ область, то  ${\sf Z}$ --- {\it нулевое множество для класса\/} $\Hol (D, M)$ тогда и только тогда, когда найдется функция $h\in \psbh (D)$, для которой выполнено неравенство $\log |g_{{\sf Z}}|+h\leq M$. Это связано с тем, что  $h\in \psbh (D)$ для односвязной при $n=1$ или звёздной относительно точки при $n>1$ области $D$ в том и только том случае, когда имеет место представление $h=\Re f$ для некоторой $f\in \Hol (D)$, или $h=\log |e^f|$ \cite[предложение 2.2.13]{Klimek}. Другими словами, ${\sf Z}$ --- нулевое множество для $\Hol (D, M)$, если и только если непуст класс 
\begin{equation}\label{env:2}
	\{ h\in \psbh (D)\colon h\leq M-\log |g_{\sf Z}| \text{ на }D\}. 
\end{equation}
Некоторые подходы подобного рода возможны  и для конечносвязных областей $D\subset \CC$ \cite{Kha07}, \cite{KhaN10}, но при этом приходиться заменить функции $h$ и $M$ на их весьма малые изменения $\check h$ и $\check M$.

\subsection{Описание нулевых подмножеств}\label{subz7_3}
По этому подразделу см. {\rm \cite[\S~11]{Kh01_2}--\cite{Kha12}, \cite{KhCh}, \cite{KKh09}}. В обозначениях предыдущего подраздела задача состоит в нахождении голоморфной функции $f\neq 0$, для которой $g_{\sf Z} f\in \Hol (D , M)$, т.\,е. справедливо ограничение $\log|g_{\sf Z} f|\leq M$, или $\log |f|\leq M-\log|g_{\sf Z}|$. Аналогично предыдущим пунктам, вопрос вновь сводится к нетривиальности класса
\begin{equation}\label{env:3}
\{ h\in \psbh (D)\colon h\leq M-\log |g_{\sf Z}| \text{ на }D\}.	
\end{equation}
Эта постановка имеет двойственные выходы на проблемы аппроксимации в пространствах функций, прежде всего экспоненциальными системами, на существование для  голоморфных функций голоморфных функций-мультипликаторов, <<погашающих>> их рост, и др.  \cite{Koo}, \cite[\S~10]{Kh01_2}, \cite{Kha12}.

\subsection{Представление мероморфных функций}\label{reprMer7_4}  
По этому подразделу см.  {\rm \cite[\S~12]{Kh01_2}, \cite{Kha07}, \cite[1.4]{Kh10}, \cite{Kha02_surm}}.
Пусть $Q=q_1/q_2$ --- мероморфная функция в области $D$ и $q_1,q_2 \in \Hol(D)$, $q_1,q_2 \not\equiv 0$. Задача состоит в возможности представления функции $Q\in \Mer  (D)$ в виде отношения двух функций из класса $\Hol (D, M)$, возможно {\it без общих нулей.\/} При этом решение ее  наиболее естественно искать в терминах  функций $M$ и  $U_Q:=\max\bigl\{\log|q_1|, \log |q_2|\bigr\}$ или $U_Q:=\log \sqrt{|q_1|^2+|q_2|^2}$ в связи с известными определениями различных вариантов характеристики Неванлинны функции $Q$ именно через $U_Q$. При этом задачу можно в несколько ослабленной форме переформулировать как поиск условий, при которых класс 
\begin{equation}\label{env:4}
\{ h\in \psbh (D)\colon h\leq M-U_Q \text{ на }D\}	
\end{equation}
или класс (при дополнительном требовании <<без общих нулей>>)
\begin{equation}\label{env:5}
\{ h\in \phar (D)\colon h\leq M-U_Q \text{ на }D\}	
\end{equation}
нетривиален.

\subsection{Комплексная теория потенциала}\label{ktp} По этому разделу см. \cite{Pol91}--\cite{Pol12}, \cite{Klimek}. Основные объекты, такие, как (плюри)гармонические меры, функции Грина и им подобные, (максимальные) решения задачи Дирихле и пр., на которые опираются применения этой теории, строятся как верхняя огибающая специальных, чаще всего выпуклых и ограниченных сверху некоторой функцией-мажорантой $M$, семейств (плюри)субгармонических функций.

\subsection{Теория равномерных алгебр} Многочисленные применения в этой теории нашла теорема  двойственности Эдвардса \cite{Edwards}, \cite[1.2]{Gamelin}.  Конкретнее, пусть $H$ --- некоторый выпуклый конус вещественнозначных  полунепрерывных сверху функций  на некотором компактном топологическом пространстве $K$  со значениями в $\RR_{-\infty}$ и $H$ содержит все константы. Через $J_a (H)$ обозначим класс  {\it мер Йенсена\/} в точке $a\in K$ относительно  $H$, а именно: положительных мер Радона $\mu$ на $K$, удовлетворяющих условию 
\begin{equation}\label{Jm}
h(a)\leq \int h\dd \mu \quad \text{ для всех }\; h\in	H.
\end{equation}
По одной из теорем Эдвардса из \cite[2]{Edwards} {\it для любой полунепрерывной снизу функции $x\in \RR_{+\infty}^K$}
\begin{equation}\label{dual:E}
\sup \{ h(a)\colon h\in H, \; h\leq x  \}	=\inf \left\{ \int x \dd \mu  \colon \mu \in J_a(H)\right\}.
\end{equation}
Отсюда, в частности, следует, что {\it множество 
$ \{ h\colon h\in H, \; h\leq x  \}$
непусто тогда и только тогда, когда хотя бы для одной точки $a\in K$ конечна правая часть в \eqref{dual:E}.}

\subsection{Связь с задачами \ref{pr1}--\ref{pr_3} из подраздела \ref{stpr}} Трактовки проблем из \ref{nontr7_1}--\ref{reprMer7_4}, касающиеся условий нетривиальности классов  \eqref{env:1}--\eqref{env:5}, можно сформулировать в следующей общей форме. 
Пусть $H$ --- выпуклый конус или, более общ\'о, выпуклое множество в векторной решетке $X$ с отношением порядка $\leq$. {\it Для каких $X_0\subset X$ при всех $x\in X_0$ множество  
$\{h\in H \colon h\leq x\}$	
непусто\/}? Если $q\colon X\to \RR$ --- некоторая функция на $X$, то {\it это множество  непустое тогда и только тогда, когда выполнено соотношение\footnote{Напоминаем, что  $\sup \varnothing :=-\infty$ и $\inf \varnothing :=+\infty$ для пустого  множества $\varnothing$.}}
\begin{equation}\label{spfH}
-\infty <\sup \bigl\{q(h) \colon h\in H , \; h\leq x\bigr\}\overset{\eqref{eq4_1}}{=}\spf_{H,q}(x).		
\end{equation}
Собственно, это и есть частный случай задачи \ref{pr_3} из подраздела \ref{stpr}.
Если удается получить равенство вида \eqref{dual:E}, где в роли $q$ выступает мера Дирака $\delta_a(h)=h(a)$ в точке $a$, то 
это решает в определённом смысле частные случаи задач
\ref{pr1} и \ref{pr2}. В частности, в таком случае соотношение \eqref{spfH} с $q:=\delta_a$ выполнено тогда и только тогда, когда 
\begin{equation*}
-\infty <\inf \left\{ \int x \dd \mu  \colon \mu \in J_a(H)\right\}, 
\end{equation*}
где $\inf$ в правой части берётся по действиям на $x\in X_0$ всевозможных выметаний меры Дирака в точке $a$, т.\,е. ростка $\Spr(\delta_a; H, X_0, \lin^+\RR^X)$, определённого в  подразделе \ref{Sbal}, определение \ref{balpb}.

\subsection{Некоторые дополнения к предшествующим  определениям и обозначениям}\label{adddef}

\paragraph{\bf Множества, топология, порядок}\label{1_1_1}
Для $n,m\in \NN$ аффинные пространства  $\CC^n$ и  $\RR^m$ соответственно над $\CC$ и $\RR$ наделяются стандартной евклидовой нормой-модулем $|\cdot|$.
Полагаем $\RR_{\infty}^m:=(\RR^m)_{\infty}$, $\CC_{\infty}^n:=(\CC^n)_{\infty}$, $\CC_{\infty}:=(\CC^1)_{\infty}$ --- одноточечные  компактификации Александрова;   $|\infty|:=+\infty$. При необходимости $\CC^n$ и 
 $\CC_{\infty}^n$ отождествляем соответственно с $\RR^{2n}$ и  $\RR_{\infty}^{2n}$ (над $\RR$). Далее, когда  возможно, обозначения вводятся и определения даются только для $\RR^{m}$ и   $\RR_{\infty}^{m}$.  Для подмножества  $S\subset \RR_{\infty}^m$ через $\clos S$, $\Int S$ и  $\partial S$  обозначаем соотв. {\it замыкание,\/} {\it внутренность\/} и {\it границу\/} $S$ в  $\RR^m_{\infty}$.  (Под)область в $\RR_{\infty}^m$ --- открытое связное  подмножество в $\RR_{\infty}^m$. 
 Для $S_0 \subset S\subset \RR^m_{\infty}$ пишем $S_0\Subset S$, если $\clos S_0$ --- {\it компактное подмножество\/} в $S$ в топологии, индуцированной с $\RR_{\infty}^m$ на $S$. Для   
 $r\in \RR_{+\infty}^+$ и $x\in \RR^m$ полагаем $B(x,r):=\{x' \in \RR^m \colon |x'-x|<r\}$ --- открытый шар радиуса $r$ с центром  $x$, $B(r):=B(0,r)$ и $B(x,+\infty)= \RR^m$;  $B(\infty,r):= \{x\in \RR^m \colon |x|>1/r\}$ и $B(\infty, +\infty):=\RR^m_{\infty}\setminus \{0\}$. При $r\leq 0$ естественно  $B(0,r):=\varnothing$. При $r>0$ полагаем $\overline B(\cdot,r):=\clos B(\cdot,r)$ --- замкнутые шары, но 
$\overline B(x,0):=\{x\}$ и при $r<0$ естественно $\overline B(\cdot,r):=\varnothing$. {\it Положительность\/} всюду понимается, в соответствии с контекстом, как $\geq 0$, а $>0$ --- {\it строгая положительность}. Аналогично для {\it отрицательности.\/} Открытые (замкнутые с непустой внутренностью) шары строго положительного  радиуса с центром 
$x\in \RR^m_{\infty}$ образуют открытую (соответственно замкнутую) базу окреcтностей точки $x\in  \RR^m_{\infty}$.

\paragraph{\bf Функции}\label{1_1_2}  Для произвольной функции $f\colon X\to Y$ 
 допускаем, что не для всех  $x\in X$  определено значение $f(x)\in Y$.  {\it Функция $f$ расширенная числовая,\/} если ее образ --- подмножество в $\RR_{\pm \infty}$. 

Для  открытого подмножества   $\mathcal O \subset \RR^m$ через  $\har ({\mathcal O})$ обозначаем  векторное  пространство над $\RR$ {\it гармонических\/} (аффинных  при $m=1$) в ${\mathcal O}$ функций; $\sbh ({\mathcal O})$ --- выпуклый конус над $\RR^+$ {\it субгармонических\/} (выпуклых при $m=1$) в ${\mathcal O}$ функций. 		 Функцию, {\it тождественно равную\/} $-\infty$ на ${\mathcal O}$, обозначаем символом $\boldsymbol{-\infty}\in \sbh ({\mathcal O})$; $\sbh_*(\mathcal O ):=\sbh (\mathcal O)\setminus \{\boldsymbol{-\infty}\}$. 	Для  $\mathcal O \subset \CC^n$ через  $\Hol ({\mathcal O})$  обозначаем  векторное  пространство над $\CC$ {\it голоморфных\/}  функций на ${\mathcal O}$.   

При $p\in \NN\cup \{\infty\}$ класс $C^p$ состоит из числовых функций на открытых подмножествах с непрерывными (частными) производными  до порядка $p$ включительно.

Через $\dist (\cdot , \cdot) $ обозначаем {\it функции евклидова расстояния\/} между парами точек, между  точкой и множеством, 
между множествами в $\RR^m_{\infty}$. По определению $\dist (\cdot , \varnothing):=:\dist (\varnothing,\cdot ):=\inf \varnothing:=+\infty=:\dist (x,\infty):=:\dist (\infty, x)$ при  $x\in \RR^m_{\infty}$. 

\paragraph{\bf Меры}  Далее ${\Meas} (S)$ --- класс {\it борелевских вещественных мер\/} на подмножествах из $S\in \mathcal B ( \RR^m_{\infty})$ со значениями в $\RR_{\pm\infty}$,       иначе называемых  {\it зарядами;\/}   $ {\Meas}_{\comp}(S)$  --- подкласс мер в ${\Meas} (S)$ с компактным {\it носителем\/} $\supp \nu\Subset S$, ${\Meas}^+ (S):=\bigl({\Meas} (S)\bigr)^+$. 
Для $x\in \RR_{\infty}^m$ и $0<r\in \RR_+$ полагаем $\mu(x,r):=\mu\bigl(B(x,r)\bigr)$.
Меру Рисса функции ${u}\in \sbh ({\mathcal O})$, $\mathcal O\subset \RR^m$,  чаще всего обозначаем как 
\begin{equation}\label{df:cm}
\nu_u:= \frac{1}{s_{m-1}} \,\Delta u\in {\Meas}^+(\mathcal O),  \quad\text{где } s_{m-1}:=\frac{2\pi^{m/2}\max\bigl\{1, (m-2)\bigr\}}{\Gamma(m/2)}
\end{equation}
--- площадь $(m-1)$-мерной единичной сферы $\partial B(1)$  в $\RR^m$, $\Delta$ --- {\it оператор Лапласа,\/}  действующий в смысле теории распределений, или обобщённых функций, а $\Gamma$ --- гамма-функция. Такие меры $\nu_u$ --- меры Радона, т.\,е. определяют линейный положительный непрерывный и ограниченный  функционал на пространстве $C_0(\mathcal O)$ непрерывных финитных функций на $\mathcal O$. В частности, $\nu_u(S)<+\infty$ для каждого измеримого  по $\nu_u$ подмножества $S\Subset \mathcal O$. При $u=\boldsymbol{-\infty}\in \sbh(\mathcal O)$ по определению   $\nu_{\boldsymbol{-\infty}}(S):=+\infty$ для всех $S\subset \mathcal O$.
 
Через $\lambda_m\in \mathcal{M}^+ (S)$ обозначаем сужения {\it меры Лебега\/} на собственные  борелевские подмножества  $S\subset \RR^m$; $\delta_x\in \mathcal{M}^+ (S)$ --- {\it мера Дирака в точке\/} $x\in S$, т.\,е. $\supp \delta_x=\{x\}$ и $\delta_x\bigl(\{x\}\bigr)=1$. В обозначении меры Лебега индекс $m$ часто будем опускать. Во включении $\Meas_{\infty}(\mathcal O)\subset \Meas(\mathcal O)$ с открытым множеством $\mathcal O$, первое множество $\Meas_{\infty}(S)$ состоит из  всех мер $\mu$ с 
бесконечно дифференцируемой плотностью, т.\,е. $\dd \mu= f\dd \lambda$, где $f$ из класса $C^{\infty}$ на $\mathcal O$.

\begin{definition}[вариация  определения \ref{balpb}]\label{bal_meas}
Пусть $D$ ---  подобласть в $\RR^m$,
$H \subset \sbh (D)$, $\nu\in \Meas^+(D)$,
$\mathcal M \subset {\Meas}^+ (D)$.
Говорим, что мера $\mu \in  \mathcal M$ ---
{\it линейное выметание  меры $\nu$ относительно $H$ в $\mathcal M$,}\/
и пишем $\nu {\prec}_{H, \mathcal M} \mu$,  если
\begin{equation}\label{eq7_1}
\int h \dd \nu \leqslant \int h \dd\mu\quad \text{для всех функций $h\in H$}.
\end{equation}
 При $\nu=\delta_x$, $x\in D$, линейные выметания $\mu \in  \mathcal M$ меры Дирака $\delta_x$ относительно $H$ в $\mathcal M$ называют {\it мерами Йенсена для точки $x\in D$,\/} если выбрано   $H:=\sbh_*(D)$ и  $\mathcal M=\Meas_{\comp}^+(D)$ (см. и \eqref{Jm}).

Говорим, что аффинная функция  $\mu+c:=\mu+c\cdot {\boldsymbol{1}}  \in  \mathcal M+\RR \boldsymbol{1}$ ---
{\it аффинное выметание  меры $\nu$ относительно $H$ в $\mathcal M+\RR \boldsymbol{1}$,}\/
и пишем $\nu {\curlyeqprec}_{H, \mathcal M} \mu+c$,  если
\begin{equation}\label{eq7_1a}
\int h \dd \nu \leqslant \int h \dd\mu+c\quad \text{для всех функций $h\in H$}.
\end{equation}
Нижние индексы ${H, \mathcal M}$ в ${\prec}_{H, \mathcal M}$
и ${\curlyeqprec}_{H, \mathcal M}$ опускаем и пишем соотв. просто  
${\prec}$ и ${\curlyeqprec}$, если из контекста ясно, какие $H$ и $\mathcal M$ подразумеваются. Очевидно, если 
$\nu {\prec}\mu$, то  $\nu {\curlyeqprec} \mu+c$ для любого  $c\in \RR^+$ в \eqref{eq7_1a}. 
\end{definition}

\paragraph{{\bf Нули  голоморфных функций} {\rm \cite[\S~11]{Ron77}, \cite{Ron86}, \cite{Ron92}, \cite{GL}, \cite{Kha12},  \cite[гл.~1]{Chi}}}\label{zf}
Пусть $D$ ---   подобласть в $\CC_{\infty}^n$, $0\neq f\in \Hol(D)$. {\it Дивизором нулей} функции $f$ называем функцию  $\Zero_f\colon D \to \NN_0$, равную кратности нуля функции $f$ в каждой точке $z\in D$. 
  Для $f=0\in \Hol (D)$ по определению $\Zero_0\equiv +\infty$ на  $D$. В главе \ref{gh2} 	далее всюду 
\begin{equation*}
\boxed{D\neq \varnothing \text{ {\sffamily{ --- область в $\RR^m$ или в $\CC^n$.}}}}
\end{equation*}

\section{Существование (плюри)субгармонической нижней огибающей}\label{exse}
\setcounter{equation}{0}

\subsection{Основной результат для выпуклых подмножеств субгармонических функций}\label{main}
\begin{propos}[{\rm частный случай предложения \ref{pr7_1}, см. и \cite[предложение 7.1]{Kh01_2}}]\label{pr7_1} 
Пусть $H$ --- непустое подмножество в $\sbh_* (D)$ и ненулевая  
 мера $\nu \in {\Meas}_{\comp}^+(D)$ такова, что
\begin{equation}\label{cnu}
-\infty < \int h \dd \nu\quad \text{для всех функций $h\in H$}. 
\end{equation}
Пусть $\mathcal M\subset \Meas^+_{\comp}(D)$ и  расширенная числовая функция $F$ на $D$ локально универсально измеримая для $\mathcal M$, т.\,е. для любого компакта $K\subset D$ сужение $F\bigm|_K$   является $\mu\bigm|_K$-измеримой функцией для всех $\mu \in \mathcal M$. Если существует функция $h\in H$,  
для которой  $h\leqslant F$ на $D$, то  
\begin{subequations}\label{exhF}
\begin{align}
-\infty  &\overset{\eqref{eq7_1}}{<}
\inf \left\{ \, \int F\, d\mu \colon
\nu {\prec } \mu  \right\}.
\tag{\ref{exhF}l}\label{exhFl}\\
-\infty  &\overset{\eqref{eq7_1a}}{<} \inf \left\{ \, \int F\, d\mu +c\colon
\nu {\curlyeqprec} \mu +c \right\}, 
\tag{\ref{exhF}a}\label{exhFa}
\end{align}
\end{subequations}
\end{propos}

Приведённая  ниже теорема \ref{th7_1}, в определенной степени обратная к предложению \ref{pr7_1}, и является основным результатом раздела \ref{exse}.

Пусть $(u_k)_{k\in \NN}$ --- последовательность субгармонических в $D$ функций, (равномерно по $k$) локально  ограниченных сверху (т.е. на компактах) в $D$.
Через $\limsup_k^* u_k$ будем обозначать полунепрерывную сверху
регуляризацию (поточечного) верхнего предела последовательности
$(u_k)$, которая также является субгармонической
функцией. При этом $\limsup_k^* u_k\in \sbh_*(D)$,  если  $\limsup_n u_n(x)\not\equiv -\infty$ 
 на $D$, и $\lambda_m$-п.\,в. совпадает с $\limsup_k u_k$  \cite[Приложение I]{GL}. Поэтому $\limsup_k u_k$ и 
$\limsup_k^* u_k$ в $L^1_{\loc}(D,\lambda_m)$ для последовательностей $(u_k)_{k\in \NN}\subset \sbh_*(D)$,  локально  ограниченных сверху, можно не различать.

Функция $f$ называется {\it  строго положительной\/}
(соотв. {\it  строго отрицательной\/}) на $D$,
если $f(x)> 0$ (соотв. $f(x) < 0$) для всех $x\in D$.
Положительная функция $f\colon D\to \RR^+$ {\it локально отделена от нуля,\/} если $\inf_{x\in K} f(x)>0$ для любого компакта $K\Subset D$.
Очевидно, полунепрерывная снизу строго положительная функция на $D$ локально отделена от нуля.
 
Для сглаживания функций и мер  операцией свертки $*$ и её <<скользящей>> версией будет использована  аппроксимационная единица $a\colon \RR^m\to \RR^+$, зависящая только от $|x|$, со свойствами   
\begin{equation}\label{ahat}
 a(x)\equiv 0 \text{ при $|x|\geq 1$}, \quad  \int_{\RR^m} a(x) \dd \lambda_m(x)=1, \quad a\in C^{\infty}.
\end{equation}
 В шаре $B(r)$, $r>0$, эта функция задаётся как 
\begin{equation}\label{ahate}
a_{r}(x)\overset{\eqref{ahat}}{:=}
c_{r}a\Bigl(\frac{1}{r} \,x\Bigr), \quad c_{r}\in \RR^+,
\quad \int_{\RR^m}  a_{r}\dd \lambda_m=1,
\end{equation}  
а соответствующая ей мера $\alpha^{(r)}\in \Meas^+(\RR^m)$
через её плотность как
\begin{equation}\label{ahatm}
\dd \alpha^{(r)}\overset{\eqref{ahate}}{:=}
a_{r}\dd \lambda_m, \quad \supp \alpha_{r}\subset 
\overline B(r).
\end{equation}
Рассмотрим строго положительную отделенную от нуля функцию  
\begin{equation}\label{eq7_3}
r\colon D\to \RR^+, \quad r(x) < \dist (x,\partial D) \quad \text{при всех $x\in D$}.
\end{equation}
Тогда для каждой точки $x\in D$ определена мера 
$\alpha_x^{(r(x))}\in \Meas^+\bigl(B(x, r(x))\bigr)$ с носителем  $\supp \alpha_x^{(r(x))}\subset D$,полученная сдвигом меры $\alpha^{(r(x))}$ в точку $x$: $\alpha_x^{(r(x))}(S):=\alpha^{(r(x))}(S-x)$. 

\begin{theorem}\label{th7_1}  Пусть  $\varnothing \neq H\subset \sbh_* (D)$, $0\neq \nu\in {\Meas}_{\comp}^+(D)$, 
$F$ --- расширенная числовая функция из $L^1_{\loc}(D):=L^1_{\loc}(D,\lambda)$. 
Допустим, что 
\begin{enumerate}
\item[{\rm (c)}]\label{c} для любого компакта $K\subset D$ и любой постоянной $c$ существует функция $h\in H$, удовлетворяющая неравенству $h\leqslant c$  на $K$,
\end{enumerate} 
и выполнено одно из следующих двух условий:
\begin{enumerate}[{\rm (a)}]
\item\label{sla} 
для любой локально ограниченной сверху последовательности функций $( h_k) \subset H$ имеем  ${\limsup}_k^* h_k \in H$, если только $\limsup_k h_k(x)\not\equiv -\infty$ на $D$;
\item\label{slb} 
множество $H$ секвенциально замкнуто в $L^1_{\loc}(D)$.
\end{enumerate}
Положим \begin{equation}\label{M0}
\mathcal M:={\Meas}_{\infty}^+(D) 
\cap {\Meas}_{\comp}^+(D\setminus U_0),
\end{equation} где область 
$U_0$ такова, что 
\begin{equation}\label{U01}
\supp \nu \subset U_0\Subset D. 
\end{equation}
\begin{enumerate}[{\rm I.}]
\item\label{Ha} Если  $H$ --- выпуклый конус, содержащий отрицательную функцию, и 
\begin{equation}\label{eq7_2}
-\infty \overset{\eqref{eq7_1}}{<} \inf \left\{ \int F \dd \mu \colon \nu {\prec } \mu , \; \mu \in \mathcal M \right\} ,
\end{equation}
то для любой положительной  локально отделённой от нуля функции
\eqref{eq7_3}
найдутся строго положительная  функция $\widehat{r}\leq r$ класса $C^{\infty}$ на $D$, постоянная $C\in \RR$ и  функция $h\in H$, для которых 
\begin{equation}\label{eq7_4}
h(x)\leqslant F^{*\widehat{r}}(x)+C
\quad\text{при всех $x\in D$, где}\quad F^{*\widehat{r}}(x) :=
\int\limits_{B(\widehat{r}(x))}F(x+y)\, \dd \alpha^{(\widehat{r}(x))}(y).
\end{equation}
\item\label{Hb}  Если $H$ --- выпуклое множество и $0\in H$, то при условии 
\begin{equation}\label{eq7_2_1}
-\infty  \overset{\eqref{eq7_1a}}{<}\inf \left\{ 
\int F \dd \mu +c\colon
\nu {\curlyeqprec} \mu+c , \; \mu \in \mathcal M,\; c \in \RR^+
\right\} 
\end{equation}
для любой положительной  локально отделённой от нуля функции
\eqref{eq7_3} найдутся такая же, как в п.~{\rm \ref{Ha},} 
функция $\widehat{r}\leq r$ на $D$, постоянная $C\in \RR$ и функция $h\in H$, для которых выполнено \eqref{eq7_4}.
\end{enumerate}
\end{theorem}
\begin{proof} По условию \eqref{U01} согласно выбору  $U_0$ всегда найдутся <<промежуточная>> регулярная для задачи Дирихле область $U_1$ в $D$ с кусочно гладкой границей $\partial U_1$, для которой 
\begin{equation}\label{UD}
\varnothing \neq \supp \nu \subset U_0\Subset  U_1\Subset D. 
\end{equation}

  Ранее нами уже отмечалась в частном случае 
\begin{lemma}\label{lemotz} Для  положительной  локально отделённой от нуля функции \eqref{eq7_3} найдётся строго положительная класса $C^{\infty}$, а значит, и отделённая от нуля функция $\widehat{r}\leq r$  на  $D$ с
\begin{equation}\label{unBxr}
 \bigcup_{x\in S}  B\bigl(x, \widehat{r}(x)\bigr)\Subset D\quad \text{для любого $S\Subset D$}, \quad 
\left(\bigcup_{x\in D\setminus U_1}  B\bigl(x, \widehat{r}(x)\bigr)\right)\bigcap U_0=\varnothing.
\end{equation}
\end{lemma}
Элементарное, но довольно трудоемкое доказательство этой леммы 
опускаем.  По лемме \ref{lemotz} сразу можем считать, что уже изначально функция $r$ обладает свойствами функции $\widehat{r}$, указанными в лемме \ref{lemotz} и теореме \ref{th7_1}, п.~\ref{Ha}. 
При этом соглашении для любой меры  
\begin{equation}\label{mus}
\mu\in \Meas^+_{\comp} (D),  \quad \supp \mu \subset D\setminus U_1,
\end{equation} 
определён интеграл семейства мер $\bigl\{\alpha_x^{(r(x))}\bigr\}_{x\in D}$ по мере $\mu$  \cite{Bour}, записываемый далее как
\begin{equation}\label{murint}
{\mu }^{*r} :=\int \alpha_x^{(r(x))} \dd \mu (x)\in \Meas_{\infty}^+(D)\cap \Meas_{\comp}^+(D\setminus U_0)\subset \mathcal M,  
\end{equation}
где  принадлежность классу $\Meas_{\infty}^+(D)$ обеспечивается построением мер $\alpha_x^{(r(x))}$ после \eqref{ahate}--\eqref{ahatm} через аппроксимативную единицу из \eqref{ahat}, а также ввиду соглашения о принадлежности функции $r$ классу $C^{\infty}$, а принадлежность классу $\Meas_{\comp}^+(D\setminus U_0)$ основана на \eqref{UD} и соотношениях \eqref{unBxr}  леммы \ref{lemotz}.  При этом мера  ${\mu }^{*r}$ по её определению \eqref{murint} действует на  $F\in L^1_{\loc}(D)$ по правилу
\begin{equation}\label{intFa}
\int F \dd{\mu }^{*r} =
\int \int\limits_{|y|<r(x)} F(x+y)\dd \alpha^{(r(x))}(y)\dd \mu (x)\overset{\eqref{eq7_2}}{=}
\int F^{*r}\dd \mu .
\end{equation}

\begin{lemma}\label{lem7_3} 
Пусть  $\mu\in \mathcal M_1=\Meas_{\comp}^+(D\setminus U_1)$ --- линейное (соотв. $\mu+c\in \mathcal M_1+c\boldsymbol{1}$ --- аффинное) выметание меры $\nu$ относительно  $H\subset \sbh_*(D)$ в $\mathcal M_1$ (соотв. в $\mathcal M+\RR^+\boldsymbol{1}$) в смысле определения\/  {\rm \ref{bal_meas}} 
и соотношения \eqref{eq7_1} (соотв. \eqref{eq7_1a}), т.\,е. $\nu \prec \mu$ (соотв. $\nu \curlyeqprec \mu+c$) в $\mathcal M_1$
(соотв. в $\mathcal M_1+\RR^+\boldsymbol{1}$). Тогда в обозначении \eqref{murint} имеем $\nu \prec \mu^{*r}$ (соотв.\/ $\nu \curlyeqprec \mu^{*r}+c$) в $\mathcal M$ (соотв. в $\mathcal M+\RR^+\boldsymbol{1}$), где $\mathcal M$ из \eqref{M0}.
\end{lemma}
\begin{proof}[Доказательство леммы\/ {\rm \ref{lem7_3}}] 
Каждая мера  $\alpha^{(r)}$ из \eqref{ahatm} --- это мера Йенсена для 0  в $\sbh_*\bigl(\overline B(r')\bigr)$ при любом $r'>r$. 
Отсюда при условии $\nu \prec \mu$  из леммы  имеем
\begin{equation*}
\int h\dd \nu \leq \int h(x) \dd \mu (x)\leq 
\int \int\limits_{|y|<r(x)} h(x+y)\dd \alpha^{(r(x))} (x) \dd \mu  (x) \overset{\eqref{intFa}}{=} \int h \dd \mu^{*r},
\end{equation*}
что даёт $\nu \prec \mu^{*r}$ в $\mathcal M$ ввиду \eqref{mus}--\eqref{murint}. При условии $\nu \curlyeqprec \mu+c$  для  меры Йенсена $\alpha^{(r(x))}$ 
для нуля имеем  
\begin{equation*}
\int h\dd \nu \leq \int h(x) \dd \mu (x)+c\leq 
\int \int\limits_{|y|<r(x)} h(x+y)\dd \alpha^{(r(x))} (x) \dd \mu  (x)
+c\overset{\eqref{intFa}}{=}
\int h \dd \mu^{*r}+c,
\end{equation*}
что даёт $\nu \curlyeqprec \mu^{*r}+c$ в $\mathcal M+\RR^+\boldsymbol{1}$ ввиду \eqref{mus}--\eqref{murint}. 
\end{proof}
В условиях п.~\ref{Ha} из \eqref{eq7_2} ввиду \eqref{intFa} согласно лемме \ref{lem7_3} следует
\begin{multline}\label{Imu}
-\infty < \inf \left\{ \int F \dd \mu \colon \nu {\prec } \mu , \; \mu \in \mathcal M \right\}\leq \inf \left\{ \int F \dd \mu^{*r} \colon \nu {\prec } \mu , \; \mu \in \mathcal M_1 \right\}
\\
\overset{\eqref{intFa}}{=}
\inf \left\{ \int F^{*r} \dd \mu \colon \nu {\prec } \mu , \; \mu \in \mathcal M_1 \right\},\quad\text{где $\mathcal M_1=\Meas_{\comp}^+(D\setminus U_1)$}.
\end{multline}
В условиях п.~\ref{Hb} из \eqref{eq7_2_1} ввиду \eqref{intFa} согласно лемме \ref{lem7_3} следует 
\begin{multline}\label{Imub}
-\infty < \inf \left\{ \int F \dd \mu+c \colon \nu {\curlyeqprec} \mu +c, \; \mu \in \mathcal M,\; c\in \RR^+ \right\}\\
\leq \inf \left\{ \int F \dd \mu^{*r} +c\colon \nu {\curlyeqprec} \mu +c , \; \mu \in \mathcal M_1,\; c\in \RR^+ \right\}
\\
\overset{\eqref{intFa}}{=}
\inf \left\{ \int F^{*r}\dd \mu+c \colon \nu {\curlyeqprec} \mu +c, \; \mu \in \mathcal M_1,\; c\in \RR^+ \right\},\quad\text{где $\mathcal M_1=\Meas_{\comp}^+(D\setminus U_1)$}.
\end{multline}
Функция $F^{*r}$ в правых частях \eqref{Imu} и \eqref{Imub}
непрерывна (и даже из класса $C^{\infty}$), поскольку $F\in L^1_{\loc}(D)$ и по соглашению $r$ из класса $C^{\infty}$. 
Временно заменим функцию $F^{*r}$ на функцию
\begin{equation}\label{Fbal}
F_{\bal}^{*r}:=\begin{cases} F^{*r}\quad\text{на $D\setminus U_1$},\\
\text{гармоническое продолжение  $F^{*r}$ с $\partial U_1$ внутрь $U_1$ на $U_1$}.
\end{cases}
\end{equation}
При выполнении \eqref{Imu} или  \eqref{Imub} ввиду $\mu \in \mathcal M_1=\Meas_{\comp}^+(D\setminus U_1)$  имеем соотв.
\begin{subequations}\label{Fbc}
\begin{align}
-\infty&<\inf \left\{ \int F_{\bal}^{*r} \dd \mu \colon \nu {\prec} \mu , \; \mu \in \mathcal M_1 \right\}
\tag{\ref{Fbc}l}\label{Fbcl}\\
\intertext{или} 
-\infty&<\inf \left\{ \int F_{\bal}^{*r} \dd \mu+c \colon \nu {\curlyeqprec} \mu+c , \; \mu \in \mathcal M_1,\; c\in \RR^+ \right\}.
\tag{\ref{Fbc}a}\label{Fbca}
\end{align}
\end{subequations} 
Пусть теперь $\mu \in \Meas_{\comp}^+(D)$. Рассмотрим классическое выметание $\mu^{\bal}$ меры  $\mu$ из $U_1$ \cite{L}:
\begin{equation}\label{mub}
\mu^{\bal}:=\begin{cases}
\mu\bigm|_{D\setminus \clos U_1}\quad\text{на $D\setminus \clos U_1$,}\\
\text{\it выметание из $\clos U_1$ меры $\mu\bigm|_{\clos U_1}$ на границу $\partial U_1$}\quad\text{на $\clos U_1$.}   
\end{cases}
\end{equation} 
Из определения \eqref{mub} следует $\supp \mu^{\bal}\Subset D\setminus U_1$ и для непрерывной функции $F_{\bal}^{*r}$ из \eqref{Fbal} в силу её гармоничности в $U_1$ по определению классического выметания меры из \cite{L} имеем 
\begin{equation}\label{IFbal}
\int F_{\bal}^{*r} \dd \mu =\int F^{*r} \dd \mu^{\bal}.
\end{equation}
При этом если $\nu \prec \mu$ (соотв.  $\nu \curlyeqprec \mu+c$) в $\Meas_{\comp}^+(D)$ (соотв. в $\Meas_{\comp}^+(D)+\RR^+\boldsymbol{1}$), то $\nu\prec \mu^{\bal}$ (соотв. $\nu \curlyeqprec \mu^{\bal}+c$)
в $\Meas_{\comp}^+(D\setminus U_1)$ (соотв. в $\Meas_{\comp}^+(D\setminus U_1)+\RR^+\boldsymbol{1}$). Отсюда и из \eqref{Fbc} в силу  \eqref{IFbal}  имеем соотв.
\begin{subequations}\label{Fbcb}
\begin{align}
-\infty&\overset{\eqref{Fbcl}}{<}\inf \left\{ \int F_{\bal}^{*r} \dd \mu \colon \nu {\prec} \mu , \; \mu \in \Meas_{\comp}^+(D) \right\}
\tag{\ref{Fbcb}l}\label{Fbcba}\\
\intertext{или} 
-\infty&\overset{\eqref{Fbca}}{<}\inf \left\{ \int F_{\bal}^{*r} \dd \mu +c\colon \nu {\curlyeqprec} \mu+c , \; \mu \in \Meas_{\comp}^+(D),\; c\in \RR^+ \right\}.
\tag{\ref{Fbcb}a}\label{Fbcbb}
\end{align}
\end{subequations}

Положим теперь 
\begin{equation}\label{stopr}
X:=L_{\loc}^1(D), \quad X_0:=C_{\RR}(D)\subset X, \quad q_1=q:=\nu 
\end{equation}
и соотв.
\begin{equation}\label{stoprr}
\text{$H\subset \sbh_*(D)$ --- выпуклые конус или множество};
\quad \text{$L=\lin^+\RR^X$ или $L=\aff^+\RR^X$}.  
\end{equation} 
При этом $X$ можно трактовать как приведённый правильный проективный предел векторных решёток (Фреше) $L^1(D_n)$ с естественным отношением порядка $\leq$ поточечно п.\,в. --- подраздел \ref{examp}, пример 4. Здесь исчерпание области $D$ областями $D_n$ можно начинать с области $D_1\supset \supp \nu$ (см. замечание \ref{prl_subs}). 

Для удобства ссылок случаи выпуклых конуса $H$ и множества $H$ рассмотрим раздельно.

\paragraph{\ref{Ha}. {\bf $H$ --- выпуклый конус}} 
При условии \eqref{sla} теоремы \ref{th7_1} используем следствие \ref{cor_RR} из теоремы \ref{thrdq}, а при условии \eqref{slb} --- топологическую теорему \ref{th6_1}. Условие (iii-iv) следствия  
\ref{cor_RR} составляет содержание условия \eqref{sla} теоремы \ref{th7_1}, а условие \eqref{slb} теоремы \ref{th7_1} требуется в теореме  \ref{th6_1}. Условие (c) теоремы  \ref{th7_1} влечёт за собой выполнение как условия \eqref{qHi} теоремы \ref{thrdq}, так  и условия 
\eqref{ki} теоремы \ref{th6_1}. 
В рамках условия \eqref{sla} теоремы \ref{th7_1} для $H\subset \sbh_*(D)$ как условие \eqref{qHv} теоремы \ref{thrdq}, так и условие (iii) теоремы \ref{th6_1} входят в перечень основных свойств субгармонических функций. То, что $H$ содержит отрицательную функцию означает, что $H\cap (-X^+)\neq \varnothing$. В частности, выполнено условие \eqref{qHii} теоремы \ref{thrdq}. Таким образом, при выборе условия \eqref{sla} в рамках теоремы \ref{th7_1} выполнены все условия следствия \ref{cor_RR} из теоремы \ref{thrdq}, а при выборе  условия \eqref{slb} в рамках теоремы \ref{th7_1} выполнены все условия 
топологической теоремы \ref{th6_1} в дополненной усиленной  версии с $H\cap (-X^+)\neq \varnothing$. Следовательно, выполнено заключение    
\eqref{eq5_1}, которое в рассматриваемой конкретной ситуации может быть записано для любой непрерывной функции $f\in  X_0\overset{\eqref{stopr}}{=}C_{\RR}(D)$ как
\begin{multline}\label{eq5_1_cc}
\spf_{H,q}(f)\overset{\eqref{eq4_1}}{:=}
\sup \left\{\, \int h \dd \nu  \colon H\ni h\leqslant f\right\}
\\=\inf \left\{l(f) \colon l\in \lin^+\RR^X; \;  
\int h \dd \nu \leq l(g) \text{ при всех $h\in H$, $h\leq g$, $g\in C_{\RR}(D)$}\right\}
\\
=\inf \bigl\{l(f) \colon l\in \lin^+\RR^X, \;  
q {\prec_H^{X_0}} l\bigr\}
\overset{\eqref{sprqb}}{=:}\inf
\bigl(\Spr(\nu; H,X_0,\lin^+\RR^X) (f)\bigr)
\quad \text{для всех $f\in C_{\RR}(D)$}.
\end{multline}
Но по следствию \ref{Tf+c} можем $\lin^+\RR^X$ отождествить с 
$\Meas_{\comp}^+(D)$, что, впрочем, в рамках \eqref{stopr} и \eqref{stoprr} давно известно и ранее. Кроме того, отметим, что по определению \ref{bal_meas} линейного выметания меры и определению \ref{balpb} абстрактного выметания и ростка в \eqref{sprqb} в данной конкретной ситуации
\begin{equation*}
\Spr(\nu; H,X_0,\lin^+\RR^X)=\{\mu \in \Meas_{\comp}^+ (D)\colon
\nu \prec \mu  \}.
\end{equation*} 
Отсюда  и из \eqref{eq5_1_cc} 
\begin{equation}\label{bmu}
\sup \left\{\, \int h \dd \nu  \colon H\ni h\leqslant f\right\}=
\inf \left\{\, \int f \dd \mu \colon
 \nu \prec \mu \in \Meas_{\comp}^+ (D) \right\}
\quad\text{для любой $f\in C_{\RR}(D)$.}
\end{equation} 
Применяя \eqref{bmu} к непрерывной функции $F_{\bal}^{*r}$
из \eqref{Fbal} согласно \eqref{Fbcba} получаем
\begin{equation*}
\sup \left\{\, \int h \dd \nu  \colon H\ni h\leqslant F_{\bal}^{*r} \right\}=
\inf \left\{\, \int F_{\bal}^{*r} \dd \mu \colon
 \nu \prec \mu \in \Meas_{\comp}^+ (D) \right\}\overset{\eqref{Fbcba}}{>}-\infty
\end{equation*}
Следовательно, существует функция $h\in H$, которая удовлетовряет 
неравенству $h\leq F_{\bal}^{*r}$ на $D$. Так как функции $F_{\bal}^{*r}\in C_{\RR}(D)$ и $F^{*r}\in C_{\RR}(D)$ совпадают в $D\setminus U_1$, то для некоторой достаточно большой постоянной $C\in \RR$ имеем \begin{equation}\label{FFusr}
F_{\bal}^{*r}\leq F^{*r}+C\quad\text{на $D$}.
\end{equation} Это завершает доказательство  для п.~\ref{Ha}.   

\paragraph{\ref{Hb}. {\bf $H$ --- выпуклое множество}} 
При условии \eqref{sla} теоремы \ref{th7_1} используем теорему \ref{thrdv}, а при условии \eqref{slb} --- топологическую теорему \ref{th6_2}. Условие (iii-iv) следствия  
\ref{cor_RR}, требуемое в теореме \ref{thrdv}, составляет содержание условия \eqref{sla} теоремы \ref{th7_1}, а условие \eqref{slb} теоремы \ref{th7_1} требуется в теореме  \ref{th6_2}.
Условие (c) теоремы  \ref{th7_1} влечёт за собой выполнение как условия \eqref{qHi} теоремы \ref{thrdq}, требуемое в теореме \ref{thrdv}, так  и условия \eqref{ki} теоремы \ref{th6_1}, требуемое в теореме \ref{th6_2}.
В рамках условия \eqref{sla} теоремы \ref{th7_1} для $H\subset \sbh_*(D)$ как условие \eqref{qHv} теоремы \ref{thrdq}, требуемое в 
теореме \ref{thrdv}, так и условие (iii) теоремы \ref{th6_1}, требуемое в теореме \ref{th6_2},  входят в перечень основных свойств субгармонических функций. При всём этом $0\in H$. Таким образом, при выборе условия \eqref{sla} в рамках теоремы \ref{th7_1} выполнены все условия  теоремы \ref{thrdv}, а при выборе  условия \eqref{slb} в рамках теоремы \ref{th7_1} выполнены все условия 
топологической теоремы \ref{th6_2} в дополненной усиленной  версии с $0\in H$ и $q_1(0)=\int 0 \dd \nu=0$. Следовательно, выполнено заключение  \eqref{eq5_1}, которое в рассматриваемой конкретной ситуации может быть записано для любой непрерывной функции $f\in  X_0\overset{\eqref{stopr}}{=}C_{\RR}(D)$ как
\begin{multline}\label{eq5_1_conc}
\spf_{H,q}(f)\overset{\eqref{eq4_1}}{:=}
\sup \left\{\, \int h \dd \nu  \colon H\ni h\leqslant f\right\}
\\=\inf \left\{a(f) \colon a\in \aff^+\RR^X; \;  
\int h \dd \nu \leq a(g) \text{ при всех $h\in H$, $h\leq g$, $g\in C_{\RR}(D)$}\right\}
\\
=\inf \bigl\{a(f) \colon a\in \aff^+\RR^X, \;  
q {\prec_H^{X_0}} a\bigr\}
\overset{\eqref{sprqb}}{=:}\inf
\bigl(\Spr(\nu; H,X_0,\aff^+\RR^X) (f)\bigr)
\quad \text{для всех $f\in C_{\RR}(D)$}.
\end{multline}
Но по следствию \ref{cor_aff} можем $\aff^+\RR^X$ отождествить с 
$\Meas_{\comp}^+(D)+\boldsymbol{1}\cdot \RR^+$ в рамках \eqref{stopr} и \eqref{stoprr}. 
Кроме того, отметим, что по определению \ref{bal_meas} аффинного  выметания меры и определению \ref{balpb} абстрактного выметания и ростка в \eqref{sprqb} в данной конкретной ситуации
\begin{equation*}
\Spr(\nu; H,X_0,\aff^+\RR^X)\overset{\eqref{aff+}}{=}\{\mu+c\cdot \boldsymbol{1}\colon
\mu \in \Meas_{\comp}^+ (D), \; c \in \RR^+,\;
\nu \curlyeqprec \mu +c \}.
\end{equation*} 
Отсюда  и из \eqref{eq5_1_conc} для всех $f\in C_{\RR}(D)$
\begin{equation}\label{bmu_conc}
\sup \left\{\, \int h \dd \nu  \colon H\ni h\leqslant f\right\}=
\inf \left\{\, \int f \dd \mu+c \colon
 \nu \curlyeqprec \mu +c,\; \mu\in \Meas_{\comp}^+ (D),\; c\in\RR^+ \right\}.
\end{equation} 
Применяя \eqref{bmu_conc} к непрерывной функции $F_{\bal}^{*r}$
из \eqref{Fbal} согласно \eqref{Fbcbb} получаем
\begin{equation*}
\sup \left\{\, \int h \dd \nu  \colon H\ni h\leqslant F_{\bal}^{*r} \right\}=
\inf \left\{\, \int F_{\bal}^{*r} \dd \mu +c\colon
 \nu \curlyeqprec \mu +c, \; \mu\in \Meas_{\comp}^+ (D), \;
c\in \RR^+ \right\}\overset{\eqref{Fbcbb}}{>}-\infty
\end{equation*}
Следовательно, существует функция $h\in H$, которая удовлетворяет 
неравенству $h\leq F_{\bal}^{*r}$ на $D$. Так как  для некоторого  $C\in \RR$ имеем \eqref{FFusr}, это завершает доказательство  для п.~\ref{Hb}.   
\end{proof}
\begin{remark}
Если  $H$ из теоремы \ref{th7_1} --- {\it конус,\/} содержащий строго отрицательную полунепрерывную сверху функцию на $D$, то условие (c) теоремы \ref{th7_1} выполнено автоматически.
\end{remark}

В заключение подраздела   приведём простые примеры выпуклых
конусов и множеств в $\sbh(D)$, $D\subset \RR^m$  или  $D\subset 
\CC^n$, всем условиям теоремы \ref{th7_1}. Эти примеры легко получаются на  основе известных свойств  выпуклых, (плюри)субгармонических и (плюри)гармонических функций (см. также \cite[примеры 7.1--7.4]{Kh01_2}).

\begin{examples} {\bf 1.} Для произвольной области $D\subset \RR^m$ --- это выпуклые конусы $H$, равные $\sbh_*(D)$, $\har (D)$, $\conv \RR^D$ в случае выпуклости $D$. 

{\bf 2.} Для области $D\subset \CC^n$ --- это выпуклые конусы  $H$, равные $\psbh_*(D):=\psbh(D) \setminus \{\boldsymbol{-\infty}\}$, $\phar(D)$.

{\bf 3.} Пусть  $C$ --- некоторый выпуклый конус  непрерывных положительных функций на $D\subset \RR^m$ или на $D\subset \CC^n$. Тогда для любого конуса $H$ из предыдущих пунктов
 1 и 2 их подмножества
\begin{equation}\label{ex3}
\{h\in H\colon \exists f\in C, h\leq f \text{ на } D\}
\end{equation}
--- выпуклые конусы, удовлетворяющие условиям части \ref{Ha} теоремы  \ref{th7_1}. 

{\bf 4.} Пусть  $C$ --- некоторое  выпуклое множество непрерывных положительных функций на $D\subset \RR^m$ или на $D\subset \CC^n$. Тогда \eqref{ex3} --- выпуклые подмножества, удовлетворяющие условиям части \ref{Hb} теоремы  \ref{th7_1}. 
\end{examples}

\begin{remark} Отметим, что как в п. 4, так и в п. 3 можно рассматривать классы функций $C$, определённых только на части $S\subset D$, с изменениями <<{\it $h\leq f$ на $S$}>> в определениях классов \eqref{ex3}, но при этом придётся накладывать на $S$ некоторые условия, связанные, например, с принципом максимума, чтобы верхние регуляризации верхних пределов последовательностей функций из $H$, ограниченных сверху на $S$ функциями из $C$, не принимали значений $+\infty$ и принадлежали $H$. Такого рода выпуклые конусы использовались, например, в 
\cite{KhCV98}.   

Кроме того, многие инвариантные относительно определённых преобразований области $D$ выпуклые подконусы и выпуклые подмножества   конусов $H$ из пунктов 1 и 2 могут удовлетворять условиям теоремы \ref{th7_1} (чётные, периодические, инвариантные относительно конформных и биголоморфных автоморфизмов псевдовыпуклой области $D\subset \CC^n$ и т.\,д., и т.\,п.). Это тема отдельного исследования и предполагается рассмотреть её в ином месте.
\end{remark}

\subsection{Одна версия теоремы \ref{th7_1} для применений}\label{ss_ver}
Для применений к проблемам \ref{nontr7_1}--\ref{reprMer7_4} из  раздела \ref{int_hf} будет отдельно рассмотрен случай функции вида $F:=M-u$, где $u$ --- некоторая (плюри)субгармоническая функция. Следующее следствие несколько уточняет и обобщает результаты с подобной функцией $F=M-u$ из наших работ \cite{KKh09_1}, \cite{Kha93_1}--\cite{Kha12}, \cite{Kh01_0}. 
 
\begin{corollary}\label{corm} Пусть $H\subset \sbh_*(D)$, $u\in \sbh_*(D)$, $M \in L^1_{\loc} (D)$, $M(D)\subset \RR_{\pm\infty}$ и фиксированы  
\begin{equation}\label{nuU}
0\neq \nu \in \Meas^+_{\comp}(D), \quad \supp \nu \overset{\eqref{U01}}{\subset} U_0\Subset D, \quad \text{где $U_0$
--- область}. 
\end{equation}

\begin{enumerate}[{\rm I.}]
\item\label{uMI} 
Пусть выполнено \eqref{cnu},  функция $M$ локально универсально 
измеримая для множества мер\footnote{Для этого достаточно, например, полунепрерывности (сверху или снизу) функции $M$.} $\mathcal M\subset \Meas^+_{\comp}(D)$. Если  существует функция $h\in H$, удовлетворяющая неравенству
\begin{equation}\label{esthM}
u+h\leq M\text{ на $D$},
\end{equation}
то найдется постоянная $C\in \RR$, для которой 
\begin{subequations}\label{estM}
\begin{align}
\int u \dd \mu &\leq \int M\dd \mu +C\quad \text{при всех $\nu \prec \mu\in \mathcal M$},  \tag{\ref{estM}l}\label{estMl}\\
\int u \dd \mu &\leq \int M\dd \mu+c +C\quad \text{при всех $\nu \curlyeqprec \mu+c$, $\mu\in \mathcal M$, $c\in \RR^+$.} 
\tag{\ref{estM}a}\label{estMa} 
\end{align}
\end{subequations}
\item\label{uMII} Пусть выполнено условие\/ {\rm (c)} и одно из условий \eqref{sla} или \eqref{slb} из теоремы\/ {\rm \ref{th7_1}}, $\mathcal M\overset{\eqref{M0}}{:=}{\Meas}_{\infty}^+(D) 
\cap {\Meas}_{\comp}^+(D\setminus U_0)$, $r\geq 0$ --- локально отделённая от нуля функция из \eqref{eq7_3}.   

\begin{enumerate}[{\rm 1.}]
\item\label{uMII1} Пусть $H$ --- выпуклый конус, содержащий отрицательную функцию, и для некоторой постоянной  $C\in \RR$ имеет место  \eqref{estMl}. Тогда найдутся строго положительная  функция $\widehat{r}\leq r$ класса $C^{\infty}$ на $D$, постоянная $\const \in \RR$ и  функция $h\in H$, для которых 
\begin{equation}\label{eq7_4_l}
u+h\leqslant M^{*\widehat{r}}+\const\quad\text{на $D$}\qquad\text{\rm 
(ср. с \eqref{esthM})}.
\end{equation}
Кроме того, часто можно избавиться от  участия функции $\widehat{r}$ в правой части \eqref{eq7_4_l}: 
\begin{enumerate}[{\rm (i)}]
\item\label{Mi} если $M\in C(D)$, то в  \eqref{eq7_4_l} можно заменить $M^{*\widehat{r}}$ на $M$, что отличается от  \eqref{esthM} лишь дополнительным слагаемым $+\const$ в правой части; 
\item\label{Mii} если $M\in \sbh_*(D)$, то в 
\eqref{eq7_4_l} можно заменить $M^{*\widehat{r}}$ на усреднение $M$ по сферам
\begin{equation}\label{S}
S_M\bigl(x, r(x)\bigr):=\frac{1}{s_{m-1}r^{m-1}(x)}\int\limits_{\partial B(x,r(x))} M(y)\dd \sigma(y), \quad x\in D,
\end{equation} 
где $s_{m-1}$ из \eqref{df:cm}, а $\dd \sigma$ --- элемент площади поверхности сферы $\partial B\bigl(x, r(x)\bigr)$.

\end{enumerate}
\item\label{uMII2} Пусть $H$ --- выпуклое множество, $0\in H$,
и для некоторой постоянной  $C\in \RR$ имеет место  \eqref{estMa}. Тогда найдутся строго положительная  функция $\widehat{r}\leq r$ класса $C^{\infty}$ на $D$, постоянная $\const \in \RR$ и  функция $h\in H$, для которых выполнено \eqref{eq7_4_l}.
Кроме того, можно избавиться от  участия функции $\widehat{r}$ в правой части \eqref{eq7_4_l} как в\/ {\rm \ref{Mi}} и\/ {\rm \ref{Mii}}. 
\end{enumerate}
\end{enumerate}
\end{corollary}
\begin{proof} В доказательстве полагаем $F:=M-u$.

\ref{uMI}. Функция $F$ локально универсально измеримая для $\mathcal M$, поскольку функция $u\in \sbh_*(D)$ полунепрерывна сверху. По предложению \ref{pr7_1} из \eqref{exhFl} и \eqref{exhFa} получаем соотв. \eqref{estMl} и \eqref{estMa}. 

\ref{uMII1}. Выполнены все условия теоремы \ref{th7_1} в части  
\ref{Ha}, откуда ввиду $F^{*\widehat{r}}=M^{*\widehat{r}}-u^{*\widehat{r}}$ получаем 
\begin{equation}\label{uhM}
u^{*\widehat{r}}+h\overset{\eqref{eq7_4}}{\leq} M^{*\widehat{r}}+\const\quad\text{на $D$ для некоторой функции $h\in H$}.
\end{equation}
Мера $\alpha^{(r(x))}$ --- это мера Йенсена, вследствие чего 
$u\leq u^{*\widehat{r}}$ на $D$ и получаем  неравенство \eqref{eq7_4_l}. 

Если $M\in C(D)$, то $M$ равномерно непрерывна на компактах из $D$. Поэтому изначально можно выбрать
локально отделённую от нуля функцию $r>0$ из \eqref{eq7_3} столь малой, что $\sup\limits_{|y-x|\leq r(x)} M(y)\leq M(x)+1$ для всех $x\in D$. Отсюда $M^{*\widehat{r}} \leq M+1$ на $D$ и     
имеет место усиление \ref{Mi}. 

Если $M\in \sbh_*(D)$, то значение усреднения $M^{*\widehat{r}}(x)$ не превышает (см.\footnote{Это утверждение доказано для субгармонических функция на $\CC$, но доказательство без труда переносится на субгармонические функции в окрестности $\overline B\bigl(x,\widehat{r}(x)\bigr)\subset \RR^m$.} \cite[предложение 3]{BTKh16}) усреднения $S_M\bigl(x,\widehat{r}(x)\bigr)\leq S_M\bigl(x,r(x)\bigr)$ ввиду  возрастания по радиусу усреднений \eqref{S}.  

\ref{uMII2}. Выполнены все условия теоремы \ref{th7_1} в части  
\ref{Hb}, откуда получаем \eqref{uhM}. Дальнейшие рассуждения дословно те же, что и при доказательстве \ref{uMII1} после \eqref{uhM}.
\end{proof}

\section{Применения к голоморфным функциям}\label{int_hfh}
\setcounter{equation}{0}

Для области  $D\subset \CC^n$ и {\it расширенной числовой функции\/}
$M$ на $D$ рассматриваем весовые классы
\begin{equation}\label{HolDM}
\Hol(D,M):=\left\{ f\in \Hol (D)\colon \sup_{z\in D}\frac{\bigl|f(z)\bigr|}{\exp M(z)}<+\infty\right\}.
\end{equation}
Далее всюду $\mathcal M\overset{\eqref{M0}}{:=}{\Meas}_{\infty}^+(D) 
\cap {\Meas}_{\comp}^+(D\setminus U_0)$, $r\geq 0$ --- произвольная локально отделённая от нуля функция из \eqref{eq7_3},  а {\it фиксированные\/}  $\nu\in \Meas^+_{\comp}(D)$ и $U_0$ удовлетворяют \eqref{nuU}, $M\in L^1_{\loc}(D)$.

\subsection{К нетривиальности весовых классов голоморфных функций}\label{nontr7_1+}

\begin{theorem}\label{trM} Пусть $D$ --- псевдовыпуклая область,  $H=\psbh_*(D)$. Если 
\begin{equation}\label{1M}
-\infty \overset{\eqref{estMl}}{<}\inf\left\{\int M \dd \mu \colon
\nu\prec \mu\in \mathcal M \right\}
\end{equation}
то найдётся строго положительная  функция $\widehat{r}\leq r$ класса $C^{\infty}$ на $D$,  для которой при любом значении числа $a>0$ для весовой функции 
\begin{equation}\label{Mrd}
\widetilde M(z):= \inf_{0<d<\min\bigl\{1,\dist(z,\partial D)\bigr\}} \left(B_{M^{*\widehat{r}}} (z,d)+n\log \frac{1}{d}\right)
+(n+a) \log \bigl(2+|z|\bigr),
\end{equation}
зависящей только от $M$,  $\widehat{r}$ и $a$, где 
\begin{equation}\label{avB}
 B_M(z,d):=\frac{n!}{\pi^nd^{2n}} \int\limits_{B(z,d)}  M \dd \lambda\quad\text{--- усреднение функции $M$ по шару $B(z,d)$},
\end{equation}
класс $\Hol \bigl(D,\widetilde{M}\,\bigr)$ не тривиален, т.\,е. содержит ненулевую функцию.

Если $M\in C(D)$, то функцию $\widehat{r}$  можно убрать, заменив  $M^{*\widehat{r}}$ в правой части \eqref{Mrd} на $M$. 

Если $M\in \sbh_*(D)$, то для любой локально отделённой от нуля функции $d\colon D\to \RR_*^+$ с 
\begin{equation}\label{dc}
d(z)<\min\bigl\{1, \dist(z,\partial D)\bigr\},\;z\in D, \quad \sup_{z\in D} \left(\frac1{d(z)}
\sup\bigl\{d(z')\colon |z-z'|\leq d(z)\bigr\}\right)\leq A<+\infty
\end{equation}
вместо функции $\widetilde{M}$ из  \eqref{Mrd} можно использовать
функцию 
\begin{equation}\label{Msbh}
\widetilde M(z):= B_M \bigl(z,d(z)\bigr)+n\log \frac{1}{d(z)}
+(n+a) \log \bigl(2+|z|\bigr), \quad z\in D.
\end{equation}
\end{theorem}
\begin{proof} По следствию \ref{corm}, часть \ref{uMII1}, для $u=0$
из условия \eqref{1M}, соответствующего \eqref{estMl}, найдется функция $h\in \psbh_*(D)$, с которой выполнено  \eqref{eq7_4_l}, т.\,е.
\begin{equation}\label{eq7_4_l+}
h\leqslant M^{*\widehat{r}}+\const\quad\text{на $D$}.
\end{equation}
\begin{theoremA}[{\rm \cite[теорема 1]{BaiKha16}}] 
Для псевдовыпуклой области $D\subset \CC^n$ и $h \in \psbh_*(D)$ для любого числа $a>0$ найдется ненулевая функция $f\in \Hol (D)$, удовлетворяющая  оценке  
\begin{equation}\label{est:main}
\log \bigl|f(z)\bigr|\leq B_h(z,d)+ n\log \frac{1}{d}+(n+a)\log\bigl(1+|z|+d\bigr) 
\quad\text{для всех $z\in D$ и\/ $d\in \bigl(0,\dist (z, \partial D)\bigr)$.}
\end{equation}
\end{theoremA}
Применим усреднение по шарам $B(z,d)$ к обеим частям неравенства 
\eqref{eq7_4_l+} с добавками двух последних слагаемых из правой части \eqref{est:main} при дополнительном ограничении $d<1$. Тогда  
по теореме A существует ненулевая функция $f\in \Hol(D,\widetilde M)$ с функцией $\widetilde M$ из \eqref{Mrd}. Вариация теоремы \ref{trM} для $M\in C(D)$ следует из части \ref{Mi} следствия \ref{corm}. В случае $M\in \sbh_*(D)$ следует воспользоваться частью \ref{Mii} следствия \ref{corm} при $r=\frac{1}{A}d$ для функции $d$ из \eqref{dc}. При таком выборе некоторые технические 
выкладки позвляют выбрать $\widetilde{M}$ как в  \eqref{Msbh}.
\end{proof}

\subsection{К описанию нулевых множеств}\label{zero7_1+}

\begin{theorem}\label{th_zero} Пусть область $D\subset \CC^n$ односвязная при $n=1$ и звёздная относительно некоторой точки $z_0\in D$ при $n>1$, т.\,е. для любой точки $z\in D$ отрезок $[z_0,z]$ лежит в $D$. Пусть $H=\phar (D)$, $\nu:=\delta_{z_0}$ --- мера Дирака в точке $z_0\in D$, ${\sf Z}$ --- дивизор нулей некоторой ненулевой 
функции  $f_{\sf Z}\in \Hol(D)$, т.\,е. $\Zero_{f_{\sf Z}}={\sf Z}$.  Если для некоторой постоянной $C\in \RR$ имеем
\begin{equation}\label{zM}
\int \log|f_{\sf Z}| \dd \mu \leq \int M\dd \mu +C\quad \text{при всех $\delta_{z_0} \prec \mu\in \mathcal M$},
\end{equation}
то найдутся строго положительная  функция $\widehat{r}\leq r$ класса $C^{\infty}$ на $D$ и функция $f\in \Hol (D,M^{*\widehat{r}})$ с дивизором нулей $\Zero_f={\sf Z}$. Если дополнительно $M\in C(D)$, то можно выбрать такую функцию $f$ с $\Zero_f={\sf Z}$ из $\Hol (D,M)$. Если $M\in \sbh_*(D)$, то такую $f$ можно выбрать из 
$\Hol (D,M^{*r})$.
\end{theorem}
\begin{proof}
По следствию \ref{corm}, часть \ref{uMII1}, для $u=\log |f_{\sf Z}|$
из условия \eqref{zM}, соответствующего \eqref{estMl}, найдется функция $h\in \phar_*(D)$, с которой выполнено  \eqref{eq7_4_l}, т.\,е.
\begin{equation}\label{eq7_4_l+z}
\log |f_{\sf Z}|+h\leqslant M^{*\widehat{r}}+\const\quad\text{на $D$}.
\end{equation}
Для плюригармонической функции $h$ в $D$ найдётся функция $g\in \Hol(D)$ \cite[предложение 2.2.13]{Klimek}, для которой $\Re g=h$. 
Тогда функция $f:=f_{\sf Z}e^g$ ввиду \eqref{eq7_4_l+z} искомая. 
Уточнения/упрощения для $M\in C(D)$ и $M\in \sbh_*(D)$ следуют соотв. из части \ref{Mi} следствия \ref{corm} и из $M^{*\widehat{r}}\leq M^{*r}$.
\end{proof}

\subsection{К описанию нулевых подмножеств}\label{szero7_1+}

\begin{theorem}\label{th_subz} Пусть область $D\subset \CC^n$ псевдовыпуклая, $H=\psbh_*(D)$, положительная функция 
${\sf Z}\leq \Zero_{f_0}$ на $D$, где $\Zero_{f_0}$ --- дивизор нулей некоторой ненулевой функции  $f_0\in \Hol(D)$.  Если для некоторой постоянной $C\in \RR$ имеем
\begin{equation}\label{zM0}
\int \log|f_0| \dd \mu \leq \int M\dd \mu +C\quad \text{при всех $\nu \prec \mu\in \mathcal M$},
\end{equation}
то найдутся строго положительная  функция $\widehat{r}\leq r$ класса $C^{\infty}$ на $D$ и ненулевая  $g\in \Hol \bigl(D,\widetilde{M}\,\bigr)$, где $\widetilde{M}$ из \eqref{Mrd}, 
с дивизором нулей $\Zero_g\geq {\sf Z}$ на $D$. Если дополнительно $M\in C(D)$, то можно выбрать такую функцию $g\neq 0$ с $\Zero_g\geq {\sf Z}$ из $\Hol (D,M)$. Если $M\in \sbh_*(D)$, то для любой локально отделённой от нуля функции $d\colon D\to \RR_*^+$, удовлетворяющей условиям \eqref{dc}, можно выбрать такую ненулевую $g\in \Hol \bigl(D,\widetilde{M}\,\bigr)$ с $\widetilde{M}$ из \eqref{Msbh} и $\Zero_g\geq {\sf Z}$. 
\end{theorem}
\begin{proof}
По следствию \ref{corm}, часть \ref{uMII1}, для субгармонической $u:=\log |f_0|$ из условия \eqref{zM0}, соответствующего \eqref{estMl}, найдется функция $h\in \psbh_*(D)$, с которой выполнено  \eqref{eq7_4_l}:
\begin{equation}\label{eq7_4_l+0}
\log |f_0|+h\leqslant M^{*\widehat{r}}+\const\quad\text{на $D$}.
\end{equation}
По теореме A существует ненулевая функция $f\in \Hol(D)$, удовлетворяющая \eqref{est:main}. Применим усреднение по шарам $B(z,d)$ к обеим частям неравенства  \eqref{eq7_4_l+0} с добавками двух последних слагаемых из правой части \eqref{est:main} при дополнительном ограничении $d<1$. Тогда ввиду субгармоничности 
$\log \bigl|f_0(z)\bigr|\leq B_{\log |f_0|}(z,d)$ при всех $z\in D$
и функция $g:=f_0f$ --- требуемая. Вариация теоремы \ref{th_subz} для $M\in C(D)$ следует из части \ref{Mi} следствия \ref{corm}. В случае $M\in \sbh_*(D)$ следует воспользоваться частью \ref{Mii} следствия \ref{corm} при $r=\frac{1}{A}d$ для функции $d$ из \eqref{dc}. При таком выборе некоторые технические 
выкладки позволяют выбрать $\widetilde{M}$ как в  \eqref{Msbh}.
\end{proof}

\subsection{К представлению мероморфных функций}\label{repM7_4}

\begin{theorem} Пусть область $D\subset \CC^n$ псевдовыпуклая,
$H=\psbh_*(D)$, 
\begin{equation}\label{Merf}
Q=\frac{q_1}{q_2}\in \Mer (D), \quad q_1, q_2\in \Hol(D)\setminus \{0\}, \quad U_Q:=\left[\begin{gathered} 
        \max\bigl\{\log|q_1|, \log |q_2|\bigr\},\\ 
         \log \sqrt{|q_1|^2+|q_2|^2}. 
      \end{gathered} \right.
\end{equation}
Если для некоторой постоянной $C\in \RR$ имеем
\begin{equation}\label{zMUQ}
\int U_Q \dd \mu \leq \int M\dd \mu +C\quad \text{при всех $\nu \prec \mu\in \mathcal M$},
\end{equation}
то найдутся строго положительная  функция $\widehat{r}\leq r$ класса $C^{\infty}$ на $D$ и ненулевые  функции $g_1,g_2\in \Hol \bigl(D,\widetilde{M}\,\bigr)$, где $\widetilde{M}$ из \eqref{Mrd},
представляющие $Q=g_1/g_2$. Если $M\in C(D)$, то можно выбрать эти голоморфные функции $g_1,g_2\neq 0$, представляющие $Q=g_1/g_2$, из $\Hol (D,M)$. Если $M\in \sbh_*(D)$, то для любой локально отделённой от нуля  $d\colon D\to \RR_*^+$, удовлетворяющей  \eqref{dc}, можно выбрать такую  пару ненулевых функций  $g_1,g_2\in \Hol \bigl(D,\widetilde{M}\,\bigr)$ с $\widetilde{M}$ из \eqref{Msbh} и $Q=g_1/g_2$. 
\end{theorem}
\begin{proof} По следствию \ref{corm}, часть \ref{uMII1}, для субгармонической функции  $u:=U_Q$
из условия \eqref{zMUQ}, соответствующего \eqref{estMl}, найдется функция $h\in \psbh_*(D)$, с которой выполнено  \eqref{eq7_4_l}: 
\begin{equation}\label{eq7_4_l+r}
U_Q+h\leqslant M^{*\widehat{r}}+\const\quad\text{на $D$}.
\end{equation}
По теореме A существует ненулевая функция $f\in \Hol(D)$, удовлетворяющая \eqref{est:main}. Применим усреднение по шарам $B(z,d)$ к обеим частям неравенства  \eqref{eq7_4_l+r} с добавками двух последних слагаемых из правой части \eqref{est:main} при дополнительном ограничении $d<1$. Тогда ввиду субгармоничности 
$U_Q(z)\leq B_{U_Q}(z,d)$ при всех $z\in D$ и 
получаем неравенство 
$U_Q+\log |f|\leq \widetilde{M}+\const$ на $D$
для функции $\widetilde{M}$ из \eqref{Mrd}.
Отсюда при любом из двух выборов функции $U_Q$ в \eqref{Merf}
имеем 
\begin{equation*}
\left.\begin{gathered} 
\max\bigl\{\log|q_1|, \log |q_2|\bigr\}+\log |f|\\ 
\log \sqrt{|q_1|^2+|q_2|^2}+\log |f|
      \end{gathered} \right]
\leq \widetilde{M}+\const\quad\text{на $D$}.
\end{equation*}
и функции
\begin{equation}\label{g12}
\begin{cases}g_1&:=q_1f,\\ g_2&:=q_2f, \end{cases} \qquad \frac{g_1}{g_2}=\frac{q_1f}{q_2f}=\frac{q_1}{q_2}=Q,
\end{equation} 
--- искомые. Вариация теоремы \ref{trM} для $M\in C(D)$ следует из части \ref{Mi} следствия \ref{corm}. В случае $M\in \sbh_*(D)$ следует воспользоваться частью \ref{Mii} следствия \ref{corm} при $r=\frac{1}{A}d$ для функции $d$ из \eqref{dc}. При таком выборе некоторые технические 
выкладки позволяют выбрать $\widetilde{M}$ как в  \eqref{Msbh}.
\end{proof}

\begin{theorem}\label{th_Mer0} Пусть область $D\subset \CC^n$ односвязная при $n=1$ и звёздная относительно некоторой точки $z_0\in D$ при $n>1$. Пусть $H=\phar (D)$, $\nu:=\delta_{z_0}$, $Q$ и $U_Q$ --- функции из \eqref{Merf}. Если для некоторой постоянной $C\in \RR$ имеем
\begin{equation}\label{zMU0}
\int U_Q \dd \mu \leq \int M\dd \mu +C\quad \text{при всех $\delta_{z_0} \prec \mu\in \mathcal M$},
\end{equation}
то найдутся строго положительная  функция $\widehat{r}\leq r$ класса $C^{\infty}$ на $D$ и функция $f\in \Hol (D)$, не обращающаяся в нуль на $D$, для которой выполнено \eqref{g12}, 
$g_1,g_2\in \Hol (D,M^{*\widehat{r}})$ и, в частности, 
$\Zero_{g_1}=\Zero_{q_1}$, $\Zero_{g_2}=\Zero_{q_2}$. Если дополнительно $M\in C(D)$, то можно выбрать такую пару 
функций $g_1,g_2$  из $\Hol (D,M)$. Если $M\in \sbh_*(D)$, то такие $g_1,g_2$ можно выбрать из $\Hol (D,M^{*r})$.
\end{theorem}
\begin{proof}
По следствию \ref{corm}, часть \ref{uMII1}, для субгармонической $u=U_Q$ из условия \eqref{zMU0}, соответствующего \eqref{estMl}, найдется функция $h\in \phar_*(D)$, с которой выполнено  \eqref{eq7_4_l}, т.\,е. \eqref{eq7_4_l+r}. 
Для плюригармонической функции $h$ в $D$ найдётся функция $g\in \Hol(D)$ \cite[предложение 2.2.13]{Klimek}, для которой $\Re g=h$. 
Тогда для функции $f:=e^g$ функции $g_1,g_2$ вида  \eqref{g12} согласно \eqref{eq7_4_l+r} искомые. 
Уточнения/упрощения для $M\in C(D)$ и $M\in \sbh_*(D)$ следуют соотв. из части \ref{Mi} следствия \ref{corm} и из $M^{*\widehat{r}}\leq M^{*r}$.
\end{proof}
\subsection{Заключительные замечания} Мы не использовали здесь применительно к голоморфным функциям следствие \ref{corm}
в части \ref{uMII2} для выпуклого множества $H$, а также её часть 
\ref{uMI}. Содержательное применение последней части дано нами в
\cite{Kha12}, \cite{KhaRoz18} (см. также библиографию в них).
Здесь мы также не применяли аппарат потенциалов Йенсена, двойственных к мерам Йенсена, теории голоморфных потоков Е.\,М. Полецкого вкупе с голоморфными и полиномиальными дисками. Эти и другие подходы будут рассмотрены в последующих работах.

\end{document}